\documentclass[a4paper]{article}

\usepackage[a4paper,margin=2cm]{geometry}

\usepackage{amsmath,amssymb,amsfonts,amsthm}
\usepackage{graphicx}
\usepackage{subcaption}
\usepackage{float}
\usepackage{xcolor}
\usepackage[numbers,square,sort&compress]{natbib}
\usepackage{hyperref}
\hypersetup{colorlinks,citecolor=blue,linkcolor=blue,breaklinks=true}
\usepackage{booktabs}
\usepackage{colortbl}
\usepackage{enumitem}
\usepackage{epstopdf}
\usepackage{array}
\usepackage{bbding}
\usepackage{fancyvrb}
\usepackage{longtable}
\usepackage{listings}
\usepackage{rotating,rotfloat}
\usepackage{yhmath}
\usepackage{bm}
\usepackage[ruled,vlined]{algorithm2e}
\usepackage{indentfirst,csquotes}

\usepackage{multirow}

\usepackage{placeins}

\newcommand{\dx}{\,\mathrm{d}{x}}
\newcommand{\ds}{\,\mathrm{d}s}

\numberwithin{equation}{section}

\def\b{\color{blue}}
\def\r{\color{red}}

\allowdisplaybreaks

\newtheorem{theorem}{Theorem}[section]
\newtheorem{lemma}{Lemma}[section]

\newcommand{\bbm}{\begin{bmatrix}}
\newcommand{\ebm}{\end{bmatrix}}

\AtBeginDocument{%
  \setlength{\abovedisplayskip}{7pt plus 2pt minus 2pt}
  \setlength{\belowdisplayskip}{7pt plus 2pt minus 2pt}
  \setlength{\abovedisplayshortskip}{3pt plus 1pt minus 1pt}
  \setlength{\belowdisplayshortskip}{5pt plus 2pt minus 2pt}
}

\setlist[itemize,enumerate]{
  topsep=4pt,
  itemsep=2pt,
  parsep=0pt,
  partopsep=0pt
}

\begin{document}

\title{A level-set method for shape and topology optimization of Stokes--Darcy flows with Beavers--Joseph--Saffman coupling\thanks{This work was supported in part by the National Key Basic Research Program (2022YFA1004402), the National Natural Science Foundation of China under grants (No. 12401534 and No. 12471377), and the Science and Technology Commission of Shanghai Municipality (No. 22DZ2229014).}}
\author{Mingbo Zhou\thanks{School of Mathematical Sciences, East China Normal University, Shanghai 200241, China. E-mail: \texttt{52285500015@stu.ecnu.edu.cn}}
\and Jiajie Li\thanks{Department of Applied Mathematics, The Hong Kong Polytechnic University, Kowloon, Hong Kong. Corresponding author. E-mail: \texttt{jiajienumermath.li@polyu.edu.hk}}
\and Shengfeng Zhu\thanks{School of Mathematical Sciences, East China Normal University, Shanghai 200241, China. Key Laboratory of Ministry of Education \& Shanghai Key Laboratory of Pure Mathematics and Mathematical Practice, East China Normal University, Shanghai 200241, China. Corresponding author. E-mail: \texttt{sfzhu@math.ecnu.edu.cn}}}
\date{}
\maketitle

\begin{abstract}
We propose a level-set method for shape and topology optimization
of Stokes--Darcy flows with Beavers--Joseph--Saffman interface coupling.
We formulate the Stokes--Darcy system with BJS interface conditions
using a level-set representation and establish existence and uniqueness
of its solution. We further prove existence of minimizers for the
viscous-energy minimization and least-squares inverse problems
within this framework. Continuous Fr\'echet differentiability of the
control-to-state map is proved, and objective derivatives are derived
using an adjoint approach. A reaction--diffusion evolution scheme combined
with mixed finite elements is developed to find the optimized interface
and accommodate topology changes. Two- and three-dimensional numerical experiments verify the derived
sensitivities and demonstrate the effectiveness of the method in
minimizing viscous energy dissipation and reconstructing interfaces
from partial and noisy velocity observations.
\end{abstract}

\noindent\textbf{Keywords:} Topology optimization, Stokes--Darcy systems, level-set method, reaction--diffusion evolution, viscous energy minimization, interface reconstruction

\section{Introduction}
\label{sec:introduction}

Coupled free-flow and porous-medium flows arise in membrane
filtration~\cite{Hanspal2009}, resin
infusion~\cite{Haji2015}, and surface-water--groundwater
interaction and biological transport~\cite{EggenweilerRybak2020}. In these settings, fluid motion in open channels and
transport through the surrounding porous material occur on different
scales and are commonly described by coupled Stokes--Darcy
models~\cite{DiscacciatiQuarteroni2009}. The location and connectivity of
the free-flow region influence both the distribution of the flow and its
viscous losses. Designing this region therefore requires more than
adjusting channel widths: connecting separated passages or removing
unnecessary branches may also be important. An essential modeling issue
is the transmission law at the fluid--porous interface. The experiments
of Beavers and Joseph~\cite{BeaversJoseph1967}, Saffman's
analysis~\cite{Saffman1971}, and the homogenization results of J\"ager and
Mikeli\'c~\cite{AMikelic_WJager_2000} provide the basis for the
Beavers--Joseph--Saffman (BJS) description of tangential slip, complementing
normal-flux continuity and normal-stress balance.

The analysis and approximation of coupled free-flow and porous-medium
problems have been developed extensively. Representative contributions
include the variational analysis and finite element method of Layton,
Schieweck, and Yotov~\cite{WJLayton_FSchieweck_IYotov_2021}, the
discontinuous Galerkin formulation of Rivi\`ere~\cite{BRiviere_2005}, and
finite element error estimates for the Beavers--Joseph coupling by Cao
et al.~\cite{YCao_MGunzburger_XHu_FHua_XWang_WZhao_2010}. Of particular
relevance here, Pacquaut et al.~\cite{Pacquaut2012} combined a level-set
description with a stabilized mixed formulation for Stokes--Darcy flows,
including interfacial slip. Their work provides an important precedent
for representing the coupling on a fixed computational domain. Turning
such a state formulation into an optimization model introduces further
questions: how the interface-dependent operator varies with the design,
which design classes provide compactness, and how topology changes can
be incorporated into the numerical update.

Topology optimization offers a systematic way to determine flow layouts
without prescribing their connectivity. The pioneering Stokes-flow
formulation of Borrvall and Petersson~\cite{BorrvallPetersson2003} uses a
design-dependent resistance to distinguish fluid and nearly impermeable
regions. Related developments include the finite element approach of
Guest and Pr\'evost~\cite{GuestPrevost2006} for coupled
Stokes--Darcy flows and the optimization of Stokes and Darcy
flow regions by Wiker, Klarbring, and
Borrvall~\cite{WikerKlarbringBorrvall2007}. Channel-flow optimization and
extensions to steady Navier--Stokes flow were investigated in
\cite{GersborgHansenSigmundHaber2005,OlesenOkkelsBruus2006}, while Challis
and Guest~\cite{ChallisGuest2009} developed a level-set approach for Stokes
flow. Chen et al. ~\cite{ChenTanWangWang2026} proposed a prediction--correction
convolution-thresholding method for thermal-fluid topology optimization,
using Darcy penalization to represent solid regions. These studies establish the effectiveness of topology optimization
for fluid design. Nevertheless, a bulk resistance interpolation does not
by itself specify an independent BJS interface law; retaining that law
requires an explicit treatment of its geometric dependence.


A level-set function
is an implicit geometric representation~\cite{OsherSethian1988,OsherFedkiw2003}
and can support either shape updates or algorithms allowing changes of
connectivity. Early level-set approaches to optimal design include the work of Sethian
and Wiegmann~\cite{SethianWiegmann2000} and Osher and
Santosa~\cite{OsherSantosa2001}. Allaire, Jouve, and
Toader~\cite{AllaireJouveToader2002,AllaireJouveToader2004}, together with
Wang, Wang, and Guo~\cite{WangWangGuo2003}, established influential
sensitivity-based methods for structural optimization. The choice of
evolution equation remains important: transporting an existing interface
does not automatically provide a mechanism for nucleating new holes.
Combining shape and topological sensitivities addresses this issue
\cite{AllaireDeGournayJouveToader2005,AmstutzAndra2006}, and more recent
developments include quasi-Newton level-set algorithms by Blauth and
Sturm~\cite{BlauthSturm2023}. An alternative is the reaction--diffusion
method of Yamada et al.~\cite{Yamada2010}, in which a fictitious interface
energy regularizes direct updates of the level-set field. This latter
viewpoint is adopted in the present work. Level-set methods are also compatible with body-fitted discretizations.
Allaire, Dapogny, and Frey~\cite{AllaireDapognyFrey2014} combined implicit
interface evolution with explicit mesh representations. Feppon and
coauthors extended this strategy to coupled thermal fluid--structure
optimization~\cite{FepponAllaireBordeuCortialDapogny2019} and
three-dimensional computations using parallel solvers
\cite{FepponAllaireDapognyJolivet2020}; an overview is given by Dapogny and
Feppon~\cite{DapognyFeppon2023}.

Existence results for fluid topology optimization are available
in several settings. Wiker et al.~\cite{WikerKlarbringBorrvall2007}
proved existence for a regularized formulation involving Darcy
and Stokes flow. In the phase-field setting, Garcke and
Hecht~\cite{GarckeHecht2016Stokes} analyzed Stokes shape and topology
optimization, while Garcke et al.~\cite{GarckeHechtHinzeKahle2015}
established well-posedness and developed numerical approximations
for a diffuse-interface fluid optimization problem.
Surface-dependent objectives were subsequently considered
in~\cite{GarckeHechtHinzeKahleLam2016}.
Unlike these phase-field formulations, the present approach
uses a level-set design variable subject to an Eikonal constraint
motivated by signed-distance representations.
Existence results for fluid optimization in such constrained
level-set classes remain limited, particularly when the
fluid--porous coupling depends on the interface geometry.
The analytical challenge is to obtain compactness that preserves
the Eikonal constraint and controls the normal field and
tangential projector appearing in the diffuse BJS term,
in addition to the phase indicators in the bulk equations.

The present work instead uses a fixed-domain
level-set representation, allowing shape and topology changes during optimization.
It also addresses interface reconstruction from velocity observations.
Unlike Afshar and Hu~\cite{AfsharHu2022}, who estimated spatially
distributed permeability in a Navier--Stokes--Darcy model, we prescribe
the permeability and recover the fluid--porous interface through
a perimeter-regularized least-squares formulation without a volume
constraint. The main contributions are threefold:
\begin{itemize}
  \item We establish state well-posedness and existence of minimizers
  for volume-constrained viscous-energy minimization and
  perimeter-regularized interface reconstruction.

  \item We prove local continuous Fr\'echet differentiability of the
  control-to-state map and derive well-posed adjoint systems and
  level-set sensitivities including the geometric BJS contributions.

  \item We develop a reaction--diffusion optimization algorithm
  and validate it through Taylor tests and two- and three-dimensional
  numerical experiments.
\end{itemize}

The remainder of the paper is organized as follows. Section~2 introduces
the regularized model and establishes the existence results. Section~3
develops the differentiability and adjoint analysis and the continuous
level-set evolution. Section~4 describes the finite element
discretization and optimization algorithm. Numerical verification and
design experiments are presented in Section~5. \ref{app:existence}
contains the supporting lemmas and the proof of the optimal-design
existence theorem.


\section{Optimization model}
\label{sec:model}

\subsection{Sharp-interface model and interface conditions}

Let $\Omega\subset\mathbb R^d$, $d\in\{2,3\}$, be a bounded, connected
Lipschitz domain.  The free-flow and porous regions, denoted by $\Omega_f$
and $\Omega_p$, are disjoint Lipschitz subdomains satisfying $\overline\Omega=\overline{\Omega_f}\cup\overline{\Omega_p}$, $\Gamma:=\partial\Omega_f\cap\partial\Omega_p$,
where $\Gamma$ is their Lipschitz interface.  Let $\mathbf n_f$ and
$\mathbf n_p=-\mathbf n_f$ denote the outward unit normals of $\Omega_f$ and
$\Omega_p$, respectively, on $\Gamma$. 

\begin{figure}[!htbp]
  \centering
  \includegraphics[width=2.5in]{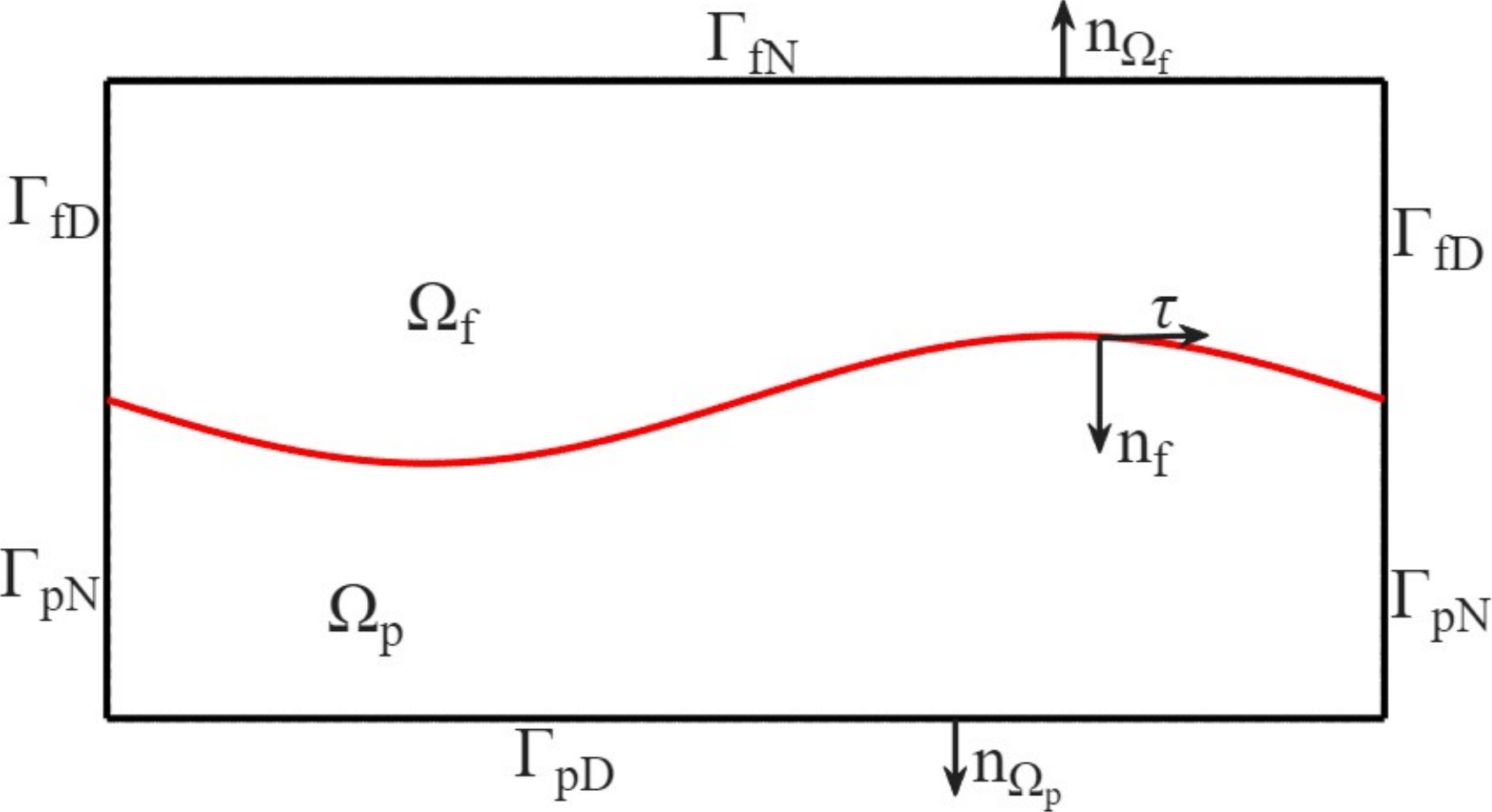}
  \caption{Decomposition of the design domain into the free-flow region
  $\Omega_f$, the porous region $\Omega_p$, and their interface $\Gamma$.}
  \label{fig:Domain}
\end{figure}

Let $\mu>0$ be the dynamic viscosity, let $K>0$ be a scalar permeability,
and let $\alpha\ge0$ be the Beavers--Joseph--Saffman coefficient.  The
free-flow stress is
\[
 \mathbf T_f(\mathbf u_f,p_f)
 :=-p_f\mathbf I+2\mu\mathcal E(\mathbf u_f),
 \qquad
 \mathcal E(\mathbf u_f)
 :=\tfrac12\bigl(D\mathbf u_f+D\mathbf u_f^{\top}\bigr).
\]
The strong equations in the two phases are
\[
 \begin{aligned}
 -\nabla\!\cdot\mathbf T_f(\mathbf u_f,p_f)&=\mathbf f_f,
 &\nabla\!\cdot\mathbf u_f&=0
 &&\text{in }\Omega_f,\\
 \frac{\mu}{K}\mathbf u_p+\nabla p_p&=\mathbf f_p,
 &\nabla\!\cdot\mathbf u_p&=0
 &&\text{in }\Omega_p.
 \end{aligned}
\]
When $\mathbf f_p=\mathbf0$, the second equation reduces to Darcy's law
\(\mathbf u_p=-(K/\mu)\nabla p_p\).

We define
\(
 \mathbf P_\tau:=\mathbf I-\mathbf n_f\otimes\mathbf n_f
\).
For sufficiently regular solutions, the interface conditions are
\begin{align}
 \mathbf u_f\cdot\mathbf n_f+\mathbf u_p\cdot\mathbf n_p&=0
 &&\text{on }\Gamma,\label{eq:interface_1}\\
 -\mathbf n_f\cdot\mathbf T_f(\mathbf u_f,p_f)\mathbf n_f&=p_p
 &&\text{on }\Gamma,\label{eq:interface_2}\\
 -\mathbf P_\tau\mathbf T_f(\mathbf u_f,p_f)\mathbf n_f
 &=\frac{\alpha\mu}{\sqrt K}\mathbf P_\tau\mathbf u_f
 &&\text{on }\Gamma.\label{eq:interface_3}
\end{align}

In the remainder of the paper, we work
with their level-set-based diffuse-interface representation.

\subsection{Level-set state equations}

Let $\phi\in W^{1,\infty}(\Omega)$ be a level-set function with
\begin{equation}\label{eq:LSdef}
 \phi<0\text{ in }\Omega_f,\qquad
 \phi=0\text{ on }\Gamma,\qquad
 \phi>0\text{ in }\Omega_p.
\end{equation}
Fix two independent regularization parameters: an interface width
$\varepsilon>0$ and an ersatz-material floor $\eta\in(0,1)$.  Define
\begin{equation}\label{eq:Heaviside}
 \mathcal H_\varepsilon(s):=
 \begin{cases}
 0,&s<-\varepsilon,\\
 \dfrac12\left[1+\dfrac{s}{\varepsilon}
 +\dfrac1\pi\sin\!\left(\dfrac{\pi s}{\varepsilon}\right)\right],
 &|s|\le\varepsilon,\\
 1,&s>\varepsilon,
 \end{cases}
 \qquad
 \delta_\varepsilon(s):=\mathcal H_\varepsilon'(s)
 =\begin{cases}
 \dfrac{1}{2\varepsilon}
 \left[1+\cos\!\left(\dfrac{\pi s}{\varepsilon}\right)\right],
 &|s|<\varepsilon,\\
 0,&|s|\ge\varepsilon.
 \end{cases}
\end{equation}
The regularized phase weights are
\begin{equation}\label{eq:Hfp}
 \mathcal H_p(\phi):=\eta+(1-\eta)\mathcal H_\varepsilon(\phi),
 \qquad
 \mathcal H_f(\phi):=\eta+(1-\eta)
 \bigl(1-\mathcal H_\varepsilon(\phi)\bigr).
\end{equation}
These regularized phase weights approximate the characteristic
functions of the two regions and satisfy
$\eta\le\mathcal H_i(\phi)\le1$ for $i\in\{f,p\}$. The strictly positive lower bound of
$\mathcal H_f$ supplies the uniform $H^1$ coercivity used below.

For a continuous integrand $g$, the diffuse-interface approximation is
\[
 \int_\Gamma g\ds
 \approx\int_\Omega g\,\delta_\varepsilon(\phi)
 |\nabla\phi|\dx.
\]
We assume that the level-set function is nondegenerate on the whole design
domain, namely,
\[
 |\nabla\phi|\ge c_0>0
 \qquad\text{a.e. in }\Omega.
\]
Consequently, the unit normal and the tangential projection are well defined
almost everywhere in $\Omega$ by
\begin{equation}\label{normalTanProj}
 \mathbf n(\phi):=\frac{\nabla\phi}{|\nabla\phi|},
 \qquad
 \mathbf P_\tau(\phi):=\mathbf I-\mathbf n(\phi)\otimes\mathbf n(\phi).
\end{equation}

We now introduce a single velocity field and a single pressure field on the
fixed design domain. First, we require the boundary datum to admit a fixed solenoidal lifting:
there exists $\mathbf u_D\in H^1(\Omega)^d$ such that
\begin{equation}\label{eq:velocity-lifting}
 \mathbf u_D=\mathbf u_d\text{ on }\partial\Omega,
 \qquad
 \nabla\!\cdot\mathbf u_D=0\text{ in }\Omega.
\end{equation}
The boundary datum $\mathbf u_d$ and the solenoidal lifting $\mathbf u_D$
are fixed independently of the design variable $\phi$.
Define
$$\mathcal U:=H^1_0(\Omega)^d,\
\mathcal P:=L^2_0(\Omega)
:=\left\{q\in L^2(\Omega):\int_\Omega q\dx=0\right\},\ \mathcal U_{\mathbf u_d}
:=\left\{\mathbf v\in H^1(\Omega)^d:  \mathbf v=\mathbf u_d\text{ on }\partial\Omega\right\}.$$ 
Let $\mathbf f\in L^2(\Omega)^d$ be a prescribed global force.  Given $\phi\in W^{1,\infty}(\Omega)$,
find
$(\mathbf u,p)\in\mathcal U_{\mathbf u_d}\times\mathcal P$ such that
\begin{equation}\label{eq:global-ls}
 \begin{aligned}
 a(\phi;\mathbf u,\mathbf v)+b(p,\mathbf v)
 &=\langle\mathbf f,\mathbf v\rangle
 &&\forall\mathbf v\in\mathcal U,\\
 b(q,\mathbf u)&=0
 &&\forall q\in\mathcal P,
 \end{aligned}
\end{equation}
where
\[
 \begin{aligned}
 a(\phi;\mathbf z,\mathbf v)
 &:={\int_\Omega}2\mu\mathcal E(\mathbf z):\mathcal E(\mathbf v)
       \mathcal H_f(\phi)\dx
   +{\int_\Omega}\frac{\mu}{K}\mathbf z\cdot\mathbf v
       \mathcal H_p(\phi)\dx\\
 &\quad
   +{\int_\Omega}\frac{\alpha\mu}{\sqrt K}
       \bigl(\mathbf P_\tau(\phi)\mathbf z\bigr)\cdot
       \bigl(\mathbf P_\tau(\phi)\mathbf v\bigr)
       \delta_\varepsilon(\phi)|\nabla\phi|\dx,\\
 b(q,\mathbf v)&:=-{\int_\Omega}q\nabla\!\cdot\mathbf v\dx,
 \qquad
 \langle\mathbf f,\mathbf v\rangle
 :={\int_\Omega}\mathbf f\cdot\mathbf v\dx.
 \end{aligned}
\]

 Setting
$\mathbf u=\mathbf w+\mathbf u_D$ gives the equivalent homogeneous problem
\begin{equation}\label{eq:global-ls-homogeneous}
 \begin{aligned}
 a(\phi;\mathbf w,\mathbf v)+b(p,\mathbf v)&=F_\phi(\mathbf v)
 &&\forall\mathbf v\in\mathcal U,\\
 b(q,\mathbf w)&=0
 &&\forall q\in\mathcal P,
 \end{aligned}
\end{equation}
where
\(
 F_\phi(\mathbf v):=\langle\mathbf f,\mathbf v\rangle
 -a(\phi;\mathbf u_D,\mathbf v)
\).

\begin{lemma}[Well-posedness of the state equations]
\label{lem:pde-wellposed}
Let $\mathbf f\in L^2(\Omega)^d$, let
$\mathbf u_D$ satisfy \eqref{eq:velocity-lifting}, and let
$\phi\in W^{1,\infty}(\Omega)$ satisfy
$|\nabla\phi|\ge c_0>0$ almost everywhere in $\Omega$.  Then
\eqref{eq:global-ls} has a unique solution
$(\mathbf u,p)\in\mathcal U_{\mathbf u_d}\times\mathcal P$.  Moreover,
\begin{equation}\label{eq:pde-apriori-estimate}
 \|\mathbf u\|_{H^1(\Omega)^d}+\|p\|_{L^2(\Omega)}
 \le C\left(
  \|\mathbf f\|_{L^2(\Omega)^d}
  +\|\mathbf u_D\|_{H^1(\Omega)^d}\right),
\end{equation}
where $C$ depends on $\Omega$, $\mu$, $K$, $\alpha$, $\eta$,
$\varepsilon$, and $\|\nabla\phi\|_{L^\infty(\Omega)}$.  Consequently,
$C$ can be chosen independently of $\phi$ on any set with a uniform
$W^{1,\infty}$ bound.
\end{lemma}

\begin{proof}
Throughout the proof, $C>0$ denotes a generic constant that depends only on $\Omega$, $\mu$, $K$, $\alpha$, $\eta$,
$\varepsilon$, and $\|\nabla\phi\|_{L^\infty(\Omega)}$. Since
$ 0\le\delta_\varepsilon(\phi)\le\varepsilon^{-1}$,
 $\eta\le\mathcal H_f(\phi),\mathcal H_p(\phi)\le 1$,
 $|\mathbf P_\tau(\phi)\mathbf\xi|\le|\mathbf\xi|$,
the interface contribution satisfies
\[
\begin{aligned}
\left|\int_\Omega\frac{\alpha\mu}{\sqrt K}
 \bigl(\mathbf P_\tau(\phi)\mathbf z\bigr)\cdot
 \bigl(\mathbf P_\tau(\phi)\mathbf v\bigr)
 \delta_\varepsilon(\phi)|\nabla\phi|\dx\right|\le
 \frac{\alpha\mu}{\sqrt K\,\varepsilon}
 \|\nabla\phi\|_{L^\infty(\Omega)}
 \|\mathbf z\|_{L^2(\Omega)^d}
 \|\mathbf v\|_{L^2(\Omega)^d}.
\end{aligned}
\]
The remaining two terms in $a$ are bounded by the Cauchy--Schwarz
inequality.  Hence, for all $\mathbf z,\mathbf v\in H^1(\Omega)^d$,
\begin{equation}\label{eq:a-continuity-proof}
 |a(\phi;\mathbf z,\mathbf v)|
 \le C\|\mathbf z\|_{H^1(\Omega)^d}
        \|\mathbf v\|_{H^1(\Omega)^d}.
\end{equation}
Moreover, for every $q\in\mathcal P$ and $\mathbf v\in\mathcal U$,
\begin{equation}\label{eq:b-continuity-proof}
 |b(q,\mathbf v)|
 \le \|q\|_{L^2(\Omega)}
       \|\nabla\!\cdot\mathbf v\|_{L^2(\Omega)}
 \le C\|q\|_{L^2(\Omega)}
       \|\mathbf v\|_{H^1(\Omega)^d}.
\end{equation}
The Darcy and interface contributions to $a(\phi;\mathbf v,\mathbf v)$ are
nonnegative.  Therefore, Korn's inequality on
$\mathcal U=H^1_0(\Omega)^d$ and the lower bound
$\mathcal H_f(\phi)\ge\eta$ yield
\[
 a(\phi;\mathbf v,\mathbf v)
 \ge2\mu\eta\|\mathcal E(\mathbf v)\|_{L^2(\Omega)}^2
 \ge C^{-1}\|\mathbf v\|_{H^1(\Omega)^d}^2,
 \qquad \forall\mathbf v\in\mathcal U.
\]
Using the Poincar\'e inequality and \eqref{eq:a-continuity-proof}, we obtain
\[
\begin{aligned}
 |F_\phi(\mathbf v)|
\le |\langle\mathbf f,\mathbf v\rangle|
      +|a(\phi;\mathbf u_D,\mathbf v)|\le C\left(\|\mathbf f\|_{L^2(\Omega)^d}
       +\|\mathbf u_D\|_{H^1(\Omega)^d}\right)
       \|\mathbf v\|_{H^1(\Omega)^d},
 \qquad \forall\mathbf v\in\mathcal U.
\end{aligned}
\]
Thus $F_\phi\in\mathcal U'$. Finally, the standard divergence inf--sup condition on a bounded Lipschitz
domain gives
\[
 \inf_{q\in\mathcal P\setminus\{0\}}
 \sup_{\mathbf v\in\mathcal U\setminus\{0\}}
 \frac{|b(q,\mathbf v)|}
 {\|q\|_{L^2(\Omega)}\|\mathbf v\|_{H^1(\Omega)^d}}
 \ge\beta_\Omega>0.
\]
The Babu\v{s}ka--Brezzi theorem applied to
\eqref{eq:global-ls-homogeneous} therefore yields a unique pair
$(\mathbf w,p)\in\mathcal U\times\mathcal P$ satisfying
\[
 \|\mathbf w\|_{H^1(\Omega)^d}+\|p\|_{L^2(\Omega)}
 \le C\|F_\phi\|_{\mathcal U'}
 \le C\left(\|\mathbf f\|_{L^2(\Omega)^d}
       +\|\mathbf u_D\|_{H^1(\Omega)^d}\right).
\]
Since $\mathbf u=\mathbf w+\mathbf u_D$, the triangle inequality gives
\eqref{eq:pde-apriori-estimate}.  Uniqueness of $(\mathbf u,p)$ follows from
the uniqueness of $(\mathbf w,p)$.
\end{proof}

\subsection{Level-set topology optimization}

We consider either the viscous-energy objective
\begin{equation}\label{eq:objective1}
 \mathcal J_1(\mathbf u,\phi)
 :=\int_\Omega2\mu|\mathcal E(\mathbf u)|^2
       \mathcal H_f(\phi)\dx,
\end{equation}
or the perimeter-regularized velocity-tracking objective
\begin{equation}\label{eq:objective2}
 \mathcal J_2(\mathbf u,\phi)
 :=\frac12\int_\Omega|\mathbf u-\mathbf u_{\rm obs}|^2\dx
   +\gamma_{\rm per}\mathcal P_\varepsilon(\phi),
 \qquad
 \mathcal P_\varepsilon(\phi)
 :=\int_\Omega\delta_\varepsilon(\phi)|\nabla\phi|\dx,
 \qquad
 \gamma_{\rm per}>0,
 \qquad
 \mathbf u_{\rm obs}\in L^2(\Omega)^d.
\end{equation}
For fixed $\varepsilon>0$, $\mathcal P_\varepsilon$ is a diffuse approximation
of the perimeter of the fluid--porous interface.

Let $V_0\in[\eta|\Omega|,|\Omega|]$ and choose
$C_{\mathcal A}>0$ sufficiently large so that
$\mathcal A_1\neq\varnothing$.
Since $\mathcal A_1\subset\mathcal A_2$, the inverse admissible
set is also nonempty.
Using the continuous representative of a $W^{1,\infty}$
function, define
\begin{equation}\label{eq:admissible-set}
 \mathcal A_1:=\left\{\phi\in W^{1,\infty}(\Omega):
 \begin{array}{l}
 \|\phi\|_{W^{1,\infty}(\Omega)}\le C_{\mathcal A},\\
 \phi\text{ is a viscosity solution of }|\nabla\phi|=1
       \text{ in }\Omega,\\
 \displaystyle\int_\Omega\mathcal H_f(\phi)\dx=V_0
 \end{array}\right\}.
\end{equation}
For the inverse problem, no volume constraint is imposed; its admissible
set is
\begin{equation}\label{eq:inverse-admissible-set}
 \mathcal A_{2}:=\left\{\phi\in W^{1,\infty}(\Omega):
 \begin{array}{l}
 \|\phi\|_{W^{1,\infty}(\Omega)}\le C_{\mathcal A},\\
 \phi\text{ is a viscosity solution of }|\nabla\phi|=1
       \text{ in }\Omega
 \end{array}\right\}.
\end{equation}
The viscosity formulation accommodates Lipschitz level-set
functions and preserves the Eikonal constraint under uniform
convergence~\cite{CrandallLions1983,CrandallIshiiLions1992}.
Together with the uniform $W^{1,\infty}$ bound, this stability
provides the compactness required for the existence analysis.
The same bound ensures that the state stability constant
is uniform over the admissible designs.
Signed-distance reinitialization is used separately
in the numerical algorithm.

For $j\in\{1,2\}$, let $S(\phi)=(\mathbf u(\phi),p(\phi))$ denote the
unique state solution pair supplied by Lemma~\ref{lem:pde-wellposed}.  The
reduced problem is
\begin{equation}\label{eq:optimization_problem}
 \min_{\phi\in\mathcal A_j}
 \mathcal J_j(\mathbf u(\phi),\phi).
\end{equation}

\begin{theorem}[Existence of an optimal design]
\label{thm:existence}
For each $j\in\{1,2\}$, assume $\mathcal A_j\ne\emptyset$.  Then problem
\eqref{eq:optimization_problem} has at least one solution
$\phi^*\in\mathcal A_j$.
\end{theorem}
The proof is given in Appendix~\ref{app:existence}, together with
the compactness and state-closedness results used in the argument.

\section{Sensitivity analysis for topology optimization}
\label{sec:sensitivity}

Throughout this section, the body force $\mathbf f$, the boundary datum
$\mathbf u_d$, and its solenoidal lifting $\mathbf u_D$ are independent of
the design variable.  For $j\in\{1,2\}$ and every
$\phi\in\mathcal A_j$, the Eikonal constraint
gives $|\nabla\phi|=1$ almost everywhere.  Hence
$\|\psi-\phi\|_{W^{1,\infty}(\Omega)}<1/2$ implies
$|\nabla\psi|\ge1/2$ almost everywhere, so the state equations and the
geometric coefficients are well defined in an open neighborhood of $\phi$.
All derivatives below are taken in this ambient neighborhood and may then
be restricted to feasible directions of $\mathcal A_j$.

For $\phi\in\mathcal A_j$, we define $$g_\phi:=|\nabla\phi|,\quad
 \mathbf n(\phi):=\frac{\nabla\phi}{g_\phi},\quad
 \mathbf P_\tau(\phi):=\mathbf I-\mathbf n(\phi)\otimes\mathbf n(\phi).$$
The standard differential identities for normalized level-set gradients
and tangential projections follow from the chain and product rules; see
\cite{DelfourZolesio2011,OsherFedkiw2003}.  Thus, for a direction
$h\in W^{1,\infty}(\Omega)$, one obtains
\begin{equation}\label{eq:geometric-directional-derivatives}
 \begin{aligned}
 g_\phi'[h]&=\mathbf n(\phi)\cdot\nabla h,\qquad
 \mathbf m_h:=D_\phi\mathbf n(\phi)[h]
   =\frac{\mathbf P_\tau(\phi)\nabla h}{g_\phi},\\
 D_\phi\mathbf P_\tau(\phi)[h]
   &=-\mathbf m_h\otimes\mathbf n(\phi)
     -\mathbf n(\phi)\otimes\mathbf m_h,\\
 \mathcal H_f'(\phi)[h]
   &=-(1-\eta)\delta_\varepsilon(\phi)h,\qquad
 \mathcal H_p'(\phi)[h]
   =(1-\eta)\delta_\varepsilon(\phi)h.
 \end{aligned}
\end{equation}
Moreover,
\begin{equation}\label{eq:surface-weight-derivative}
 D\!\left[\delta_\varepsilon(\phi)g_\phi\right][h]
 =\delta_\varepsilon'(\phi)h g_\phi
  +\delta_\varepsilon(\phi)\mathbf n(\phi)\cdot\nabla h.
\end{equation}
Consequently, the derivative of the state bilinear form with respect to
$\phi$ is
\begin{equation}\label{eq:bilinear-form-derivative}
 \begin{aligned}
 D_\phi a(\phi;\mathbf z,\mathbf v)[h]
 &=-\int_\Omega2\mu(1-\eta)
   \mathcal E(\mathbf z):\mathcal E(\mathbf v)
   \delta_\varepsilon(\phi)h\dx\\
 &\quad+\int_\Omega\frac{\mu}{K}(1-\eta)
   \mathbf z\cdot\mathbf v\,\delta_\varepsilon(\phi)h\dx\\
 &\quad+\int_\Omega\frac{\alpha\mu}{\sqrt K}
 \Big[\big(D_\phi\mathbf P_\tau(\phi)[h]\mathbf z\big)\cdot
          \big(\mathbf P_\tau(\phi)\mathbf v\big)
      +\big(\mathbf P_\tau(\phi)\mathbf z\big)\cdot
          \big(D_\phi\mathbf P_\tau(\phi)[h]\mathbf v\big)\Big]
   \delta_\varepsilon(\phi)g_\phi\dx\\
 &\quad+\int_\Omega\frac{\alpha\mu}{\sqrt K}
   \big(\mathbf P_\tau(\phi)\mathbf z\big)\cdot
   \big(\mathbf P_\tau(\phi)\mathbf v\big)
   \Big[\delta_\varepsilon'(\phi)h g_\phi
       +\delta_\varepsilon(\phi)\mathbf n(\phi)\cdot\nabla h\Big]\dx.
 \end{aligned}
\end{equation}

\subsection{Differentiability of the control-to-state map}

\begin{lemma}[Fr\'echet differentiability of the control-to-state map]
\label{lem:control-to-state-differentiability}
Let $S(\phi)=(\mathbf u(\phi),p(\phi))$ be the state map defined by
\eqref{eq:global-ls}. Set
$
\operatorname{dom}(S)
:= \bigl\{
\psi\in W^{1,\infty}(\Omega):
\operatorname*{ess\,inf}_{\Omega}|\nabla\psi|>0
\bigr\}.
$
Since every $\phi\in\mathcal A_1\cup\mathcal A_2$ satisfies
$|\nabla\phi|=1$ almost everywhere in $\Omega$, we have
$\mathcal A_1\cup\mathcal A_2\subset\operatorname{dom}(S)$.
The map
$S:\operatorname{dom}(S)\to H^1(\Omega)^d\times\mathcal P$
is continuously Fr\'echet differentiable, and
$S(\phi)\in\mathcal U_{\mathbf u_d}\times\mathcal P$ for every
$\phi\in\operatorname{dom}(S)$.  For such $\phi$ and
$h\in W^{1,\infty}(\Omega)$, its derivative $DS(\phi)[h]=(\mathbf z,\pi)
 \in\mathcal U\times\mathcal P$
is the unique solution of
\begin{equation}\label{eq:linearized-state}
 \begin{aligned}
 a(\phi;\mathbf z,\mathbf v)+b(\pi,\mathbf v)
 &=-D_\phi a(\phi;\mathbf u,\mathbf v)[h]
 &&\forall\mathbf v\in\mathcal U,\\
 b(q,\mathbf z)&=0
 &&\forall q\in\mathcal P,
 \end{aligned}
\end{equation}
where $(\mathbf u,p)=S(\phi)$.
\end{lemma}

\begin{proof}
Fix $\bar\phi\in\operatorname{dom}(S)$ and choose $0<\rho<
 \operatorname*{ess\,inf}_{\Omega}|\nabla\bar\phi|$.
Then $\mathcal V_{\bar\phi}
 :=\left\{\psi\in W^{1,\infty}(\Omega):
 \|\psi-\bar\phi\|_{W^{1,\infty}(\Omega)}<\rho\right\}
 \subset\operatorname{dom}(S)$.
Indeed, every $\psi\in\mathcal V_{\bar\phi}$ satisfies $\operatorname*{ess\,inf}_{\Omega}|\nabla\psi|
 \ge \operatorname*{ess\,inf}_{\Omega}|\nabla\bar\phi|
    -\|\nabla(\psi-\bar\phi)\|_{L^\infty(\Omega)^d}>0$.
This also proves that $\operatorname{dom}(S)$ is open in
$W^{1,\infty}(\Omega)$.  Set $X:=\mathcal U\times\mathcal P$ and write
$\mathbf u=\mathbf w+\mathbf u_D$.
For $(\phi,\mathbf w,p)\in\mathcal V_{\bar\phi}\times X$, define the
functional $\mathcal F(\phi,\mathbf w,p)\in X'$ by
\begin{equation}\label{eq:state-functional}
 \begin{aligned}
 \left\langle\mathcal F(\phi,\mathbf w,p),(\mathbf v,q)\right\rangle
 &:=
 a(\phi;\mathbf w+\mathbf u_D,\mathbf v)+b(p,\mathbf v)
 +b(q,\mathbf w+\mathbf u_D)-\langle\mathbf f,\mathbf v\rangle
 \end{aligned}
\end{equation}
for every $(\mathbf v,q)\in X$.
The functions $\mathcal H_\varepsilon$ and $\delta_\varepsilon$ are
$C^2$ and $C^1$, respectively.  On $\mathcal V_{\bar\phi}$, the maps
$\phi\mapsto g_\phi$, $\phi\mapsto\mathbf n(\phi)$, and
$\phi\mapsto\mathbf P_\tau(\phi)$ are continuously Fr\'echet differentiable
into $L^\infty(\Omega)$, $L^\infty(\Omega)^d$, and
$L^\infty(\Omega)^{d\times d}$, respectively.  Their derivatives are given
by \eqref{eq:geometric-directional-derivatives}.  The Nemytskii chain rule,
the product rule, and \eqref{eq:surface-weight-derivative} therefore show
that $\mathcal F:\mathcal V_{\bar\phi}\times X\to X'$ is continuously
Fr\'echet differentiable, with
\[
 \left\langle D_\phi\mathcal F(\phi,\mathbf w,p)[h],(\mathbf v,q)\right\rangle
 =D_\phi a(\phi;\mathbf w+\mathbf u_D,\mathbf v)[h].
\]
Here $D_\phi a$ is given by \eqref{eq:bilinear-form-derivative}; in
particular,
\begin{equation}\label{eq:bilinear-derivative-bound}
 |D_\phi a(\phi;\mathbf z,\mathbf v)[h]|
 \le C\|h\|_{W^{1,\infty}(\Omega)}
       \|\mathbf z\|_{H^1(\Omega)^d}
       \|\mathbf v\|_{H^1(\Omega)^d}
\end{equation}
locally uniformly with respect to $\phi$.

The partial derivative of $\mathcal F$ with respect to $(\mathbf w,p)$ is
\[
 \left\langle D_{(\mathbf w,p)}\mathcal F(\phi,\mathbf w,p)
 (\mathbf z,\pi),(\mathbf v,q)\right\rangle
 =a(\phi;\mathbf z,\mathbf v)+b(\pi,\mathbf v)+b(q,\mathbf z)
\]
for every $(\mathbf z,\pi),(\mathbf v,q)\in X$.
The coercivity and inf--sup estimates used in
Lemma~\ref{lem:pde-wellposed} show, by the Babu\v{s}ka--Brezzi theorem, that
this derivative is an isomorphism from $X$ onto $X'$.  The Banach
implicit function theorem now yields a locally unique continuously
Fr\'echet differentiable map $\phi\mapsto(\mathbf w(\phi),p(\phi))$.
Adding the fixed lifting $u_D$ to $w(\phi)$ yields the map
$\phi\mapsto(\mathbf w(\phi)+\mathbf u_D,p(\phi))$, which coincides with $S$
by uniqueness of the state solution.
Since $\bar\phi\in\operatorname{dom}(S)$ was arbitrary, $S$ is continuously
Fr\'echet differentiable throughout $\operatorname{dom}(S)$.

Differentiating
$\mathcal F(\phi,\mathbf w(\phi),p(\phi))=0$ in the direction $h$ gives
\eqref{eq:linearized-state}, because the fixed lifting $\mathbf u_D$ has
zero derivative with respect to $\phi$.
\end{proof}

\subsection{Adjoint equations and objective derivative}

For $j\in\{1,2\}$, define the bounded functional
$\ell_j\in\mathcal U'$ by
\begin{equation}\label{eq:objective-state-derivatives}
 \begin{aligned}
 \ell_1(\mathbf v)
 &:=D_{\mathbf u}\mathcal J_1(\mathbf u,\phi)[\mathbf v]
 =\int_\Omega4\mu\mathcal E(\mathbf u):\mathcal E(\mathbf v)
   \mathcal H_f(\phi)\dx,\\
 \ell_2(\mathbf v)
 &:=D_{\mathbf u}\mathcal J_2(\mathbf u,\phi)[\mathbf v]
 =\int_\Omega(\mathbf u-\mathbf u_{\rm obs})\cdot\mathbf v\dx.
 \end{aligned}
\end{equation}

The adjoint equations can equivalently be obtained from the Lagrangian
\begin{equation}\label{eq:Lagrangian}
 \begin{aligned}
 \mathcal L_j(\mathbf u,p,\phi;\mathbf y,r)
 :=\mathcal J_j(\mathbf u,\phi)
 -a(\phi;\mathbf u,\mathbf y)-b(p,\mathbf y)-b(r,\mathbf u)
 +\langle\mathbf f,\mathbf y\rangle.
 \end{aligned}
\end{equation}
Indeed, stationarity with respect to $(\mathbf y,r)$ gives the state system
\eqref{eq:global-ls}, while stationarity with respect to $(\mathbf u,p)$ in
directions from $\mathcal U\times\mathcal P$ gives
\eqref{eq:adjoint-weak}.
\begin{lemma}[Well-posedness of the adjoint system]
\label{lem:adjoint-wellposedness}
Let $j\in\{1,2\}$, $\phi\in\mathcal A_j$, and
$(\mathbf u,p)=S(\phi)$.  Then there exists a unique
adjoint pair $(\mathbf y_j,r_j)\in\mathcal U\times\mathcal P$ satisfying
\begin{equation}\label{eq:adjoint-weak}
 \begin{aligned}
 a(\phi;\mathbf y_j,\mathbf v)+b(r_j,\mathbf v)
 &=\ell_j(\mathbf v)
 &&\forall\mathbf v\in\mathcal U,\\
 b(q,\mathbf y_j)&=0
 &&\forall q\in\mathcal P.
 \end{aligned}
\end{equation}
\end{lemma}
\begin{proof}
Both functionals in \eqref{eq:objective-state-derivatives} belong to
$\mathcal U'$: the assertion for $\ell_1$ follows from
$\mathcal H_f(\phi)\in L^\infty(\Omega)$, and that for $\ell_2$ follows
from the continuous embedding $H^1(\Omega)^d\hookrightarrow L^2(\Omega)^d$.
The bilinear forms in \eqref{eq:adjoint-weak} are exactly those of the state
problem.  Their continuity, coercivity on the divergence-free kernel, and
inf--sup stability were established in Lemma~\ref{lem:pde-wellposed}.
The Babu\v{s}ka--Brezzi theorem therefore yields existence and uniqueness.
\end{proof}

\begin{theorem}[Directional derivative of the objective]
\label{thm:sensitivity}
Let $j\in\{1,2\}$, $\phi\in\mathcal A_j$, and
$h\in W^{1,\infty}(\Omega)$.  Let
$(\mathbf u,p)=S(\phi)$ and let $(\mathbf y_j,r_j)$ solve
\eqref{eq:adjoint-weak}.  Then the map
$\phi\mapsto\mathcal J_j(\mathbf u(\phi),\phi)$ is Fr\'echet differentiable,
and its derivative in the direction $h$ is
\begin{equation}\label{eq:objective-directional-derivative}
 \begin{aligned}
 D\!\left[\mathcal J_j(\mathbf u(\cdot),\cdot)\right](\phi)[h]
 &=\partial_\phi\mathcal J_j(\mathbf u,\phi)[h]\\
 &\quad+(1-\eta)\int_\Omega
 \left[2\mu\mathcal E(\mathbf u):\mathcal E(\mathbf y_j)
       -\frac{\mu}{K}\mathbf u\cdot\mathbf y_j\right]
 \delta_\varepsilon(\phi)h\dx\\
 &\quad+\int_\Omega\frac{\alpha\mu}{\sqrt K}
 \delta_\varepsilon(\phi)\mathbf Q_j\cdot\nabla h\dx
 -\int_\Omega\frac{\alpha\mu}{\sqrt K}
 (\mathbf u_\tau\cdot\mathbf y_{j,\tau})
 \delta_\varepsilon'(\phi)g_\phi h\dx,
 \end{aligned}
\end{equation}
where
\begin{equation}\label{eq:bjs-direction-vector}
 \begin{aligned}
 \mathbf u_\tau&:=\mathbf P_\tau(\phi)\mathbf u,\qquad
 \mathbf y_{j,\tau}:=\mathbf P_\tau(\phi)\mathbf y_j,\\
 \mathbf Q_j&:=(\mathbf u\cdot\mathbf n(\phi))\mathbf y_{j,\tau}
 +(\mathbf y_j\cdot\mathbf n(\phi))\mathbf u_\tau
 -(\mathbf u_\tau\cdot\mathbf y_{j,\tau})\mathbf n(\phi),
 \end{aligned}
\end{equation}
and the explicit derivatives of the two objectives are
\begin{equation}\label{eq:explicit-objective-design-derivatives}
 \begin{aligned}
 \partial_\phi\mathcal J_1(\mathbf u,\phi)[h]
 &=-\int_\Omega2\mu(1-\eta)|\mathcal E(\mathbf u)|^2
   \delta_\varepsilon(\phi)h\dx,\\
 \partial_\phi\mathcal J_2(\mathbf u,\phi)[h]
 &=\gamma_{\rm per}\int_\Omega
 \left[
 \delta_\varepsilon'(\phi)h\,g_\phi
 +\delta_\varepsilon(\phi)\mathbf n(\phi)\cdot\nabla h
 \right]\dx.
 \end{aligned}
\end{equation}
\end{theorem}

\begin{proof}
Let $(\mathbf z,\pi)=S'(\phi)h$ be given by
Lemma~\ref{lem:control-to-state-differentiability}.  The chain rule and
\eqref{eq:objective-state-derivatives} give
\begin{equation}\label{eq:chain-rule-objective}
 D\!\left[\mathcal J_j(\mathbf u(\cdot),\cdot)\right](\phi)[h]
 =\partial_\phi\mathcal J_j(\mathbf u,\phi)[h]+\ell_j(\mathbf z).
\end{equation}
Testing the adjoint equations with $(\mathbf v,q)=(\mathbf z,\pi)$ and the
linearized state equations with $(\mathbf v,q)=(\mathbf y_j,r_j)$, and using
the symmetry of $a$, yields
\begin{equation}\label{eq:adjoint-elimination}
 \begin{aligned}
 \ell_j(\mathbf z)
 &=a(\phi;\mathbf y_j,\mathbf z)+b(r_j,\mathbf z)\\
 &=a(\phi;\mathbf z,\mathbf y_j)
   +b(\pi,\mathbf y_j)+b(r_j,\mathbf z)\\
 &=-D_\phi a(\phi;\mathbf u,\mathbf y_j)[h],
 \end{aligned}
\end{equation}
where $b(\pi,\mathbf y_j)=0$ follows from the second adjoint equation.
Thus
\[
 D\!\left[\mathcal J_j(\mathbf u(\cdot),\cdot)\right](\phi)[h]
 =\partial_\phi\mathcal J_j(\mathbf u,\phi)[h]
  -D_\phi a(\phi;\mathbf u,\mathbf y_j)[h].
\]
It remains to expand the last term.  From
\eqref{eq:geometric-directional-derivatives},
\[
 \begin{aligned}
 \Big[\big(D_\phi\mathbf P_\tau(\phi)[h]\mathbf u\big)\cdot\mathbf y_{j,\tau}
 +\mathbf u_\tau\cdot
       \big(D_\phi\mathbf P_\tau(\phi)[h]\mathbf y_j\big)\Big]g_\phi=-\Big[(\mathbf u\cdot\mathbf n(\phi))\mathbf y_{j,\tau}
 +(\mathbf y_j\cdot\mathbf n(\phi))\mathbf u_\tau\Big]\cdot\nabla h.
 \end{aligned}
\]
Substitution of this identity and
\eqref{eq:surface-weight-derivative} into
\eqref{eq:bilinear-form-derivative} gives the last three terms in
\eqref{eq:objective-directional-derivative}.  Finally, differentiating
\eqref{eq:objective1} and \eqref{eq:objective2} with the state held fixed
gives \eqref{eq:explicit-objective-design-derivatives}.  This proves the
result.
\end{proof}

\subsection{Reaction--diffusion level-set evolution}
\label{subsec:reaction-diffusion}

Following the reaction--diffusion level-set method of Yamada et
al.~\cite{Yamada2010}, we add a fictitious interface energy to regularize
the design update and control geometrical complexity.  The volume equality,
used only for the viscous-energy problem, is represented by
\begin{equation}\label{eq:volume-functional}
 \mathcal C(\phi)
 :=\int_\Omega\mathcal H_f(\phi)\dx-V_0,
 \qquad
 D\mathcal C(\phi)[h]
 =-(1-\eta)\int_\Omega\delta_\varepsilon(\phi)h\dx.
\end{equation}
For $j=1$, let $\lambda\in\mathbb R$ be the volume multiplier and let
$\rho>0$ be the augmented-Lagrangian penalty parameter.  For the inverse
problem $j=2$, set $\lambda=\rho=0$, so that no volume term is present.
Let $\tau>0$ be the fictitious-interface regularization parameter.  We define
\begin{equation}\label{eq:yamada-energy}
 \mathcal E_{j,\lambda,\rho,\tau}(\phi)
 :=\mathcal J_j(\mathbf u(\phi),\phi)
   +\lambda\mathcal C(\phi)
   +\frac{\rho}{2}\mathcal C(\phi)^2
   +\frac{\tau}{2}\int_\Omega|\nabla\phi|^2\dx.
\end{equation}
On $\mathcal A_j$, the Eikonal constraint implies
$\frac{\tau}{2}\int_\Omega|\nabla\phi|^2\dx
=\frac{\tau}{2}|\Omega|$.
The diffusion term smooths the numerical updates, allowing
flexibility in the intermediate iterates beyond $\mathcal A_j$.
Theorem~\ref{thm:existence} establishes the existence of optimal
designs in this class, while convergence of the numerical
evolution to such designs requires further analysis.

Thus the derivative of the volume terms is
\begin{equation}\label{eq:augmented-volume-derivative}
 D\!\left[\lambda\mathcal C+\frac{\rho}{2}\mathcal C^2\right](\phi)[h]
 =\Lambda(\phi,\lambda)D\mathcal C(\phi)[h],
 \qquad
 \Lambda(\phi,\lambda):=\lambda+\rho\mathcal C(\phi).
\end{equation}
For $j=1$, the penalty term reacts immediately to the current volume error,
while the multiplier removes the residual error through outer updates.
During each inner reaction--diffusion evolution, $\lambda$ and $\rho$ are
held fixed.  These terms vanish identically for $j=2$.

The parameter $s\ge0$ used below is an artificial optimization parameter,
not physical time. Assuming sufficient regularity for the reduced objective
derivative to extend continuously to $H^1(\Omega)$,
we formally consider the $L^2(\Omega)$-gradient flow
of \eqref{eq:yamada-energy}:
Find $\phi(s)$ such
that
\begin{equation}\label{eq:yamada-flow-weak}
 \frac{1}{\kappa}\int_\Omega\partial_s\phi\,\psi\dx
 +D\!\left[\mathcal J_j(\mathbf u(\cdot),\cdot)\right](\phi)[\psi]
 +\Lambda(\phi,\lambda)D\mathcal C(\phi)[\psi]
 +\tau\int_\Omega\nabla\phi\cdot\nabla\psi\dx=0
\end{equation}
for every $\psi\in H^1(\Omega)$, where $\kappa>0$ is the mobility. The derivative is given by Theorem~\ref{thm:sensitivity}.
If there exists $G_{j,\lambda,\rho}(\phi)\in L^2(\Omega)$ such that
\[
D\!\left[
\mathcal J_j(\mathbf u(\cdot),\cdot)+\lambda C+\frac{\rho}{2}C^2
\right](\phi)[\psi]
=
\int_\Omega G_{j,\lambda,\rho}(\phi)\psi\,{\rm d}x
\qquad\forall\psi\in H^1(\Omega),
\]
then, for sufficiently regular $\phi$, \eqref{eq:yamada-flow-weak} yields
the formal strong form
\begin{equation}\label{eq:yamada-flow-strong}
 \partial_s\phi
 =-\kappa\bigl(G_{j,\lambda,\rho}(\phi)-\tau\Delta\phi\bigr)
 \quad\text{in }\Omega,
 \qquad
 \nabla\phi\cdot\mathbf n_{\partial\Omega}=0
 \quad\text{on }\partial\Omega.
\end{equation}
To prevent a progressive loss of the
signed-distance property, the level-set function is periodically
reinitialized by the fast marching method~\cite{Sethian1996}.  More precisely, after a fixed
number of descent steps, $\phi$ is replaced by the signed distance to its
current zero level set,
\begin{equation}\label{eq:fast-marching-reinitialization}
 \phi^{\rm re}(x)
 :=\begin{cases}
  \operatorname{dist}(x,\Gamma(\phi)),&\phi(x)>0,\\
  -\operatorname{dist}(x,\Gamma(\phi)),&\phi(x)<0,
 \end{cases}
 \qquad
 \Gamma(\phi):=\{x\in\Omega:\phi(x)=0\},
\end{equation}
with zero value on $\Gamma(\phi)$. 

\section{Numerical method}
\label{sec:numerical-realization}

\subsection{Finite element spaces and discrete problems}
\label{subsec:finite-element-spaces}

Let $\mathcal T_h$ be a shape-regular simplicial triangulation of $\Omega$,
and let $h$ denote its mesh size.  The level-set space is
\[
 \mathcal W_h
 :=\left\{\psi_h\in C^0(\overline\Omega):
 \psi_h|_T\in\mathbb P_1(T)\quad\forall T\in\mathcal T_h\right\}
 \subset H^1(\Omega).
\]
For each $T\in\mathcal T_h$, let $b_T$ be the element bubble and set
$\mathbb P_1^{\rm b}(T):=\mathbb P_1(T)\oplus\operatorname{span}\{b_T\}$.
The MINI velocity and pressure spaces introduced by Arnold, Brezzi, and
Fortin~\cite{ArnoldBrezziFortin1984} are
\begin{align*}
 \mathcal U_h
 &:=
 \left\{\mathbf v_h\in H^1_0(\Omega)^d:
 \mathbf v_h|_T\in\mathbb P_1^{\rm b}(T)^d
 \quad\forall T\in\mathcal T_h\right\},\\
 \mathcal P_h
 &:=
 \left\{q_h\in C^0(\overline\Omega)\cap L^2_0(\Omega):
 q_h|_T\in\mathbb P_1(T)
 \quad\forall T\in\mathcal T_h\right\}.
\end{align*}
Let $\mathbf u_{D,h}$ be a finite element lifting of the prescribed
boundary datum, and define
$\mathcal U_{h,\mathbf u_d}:=\mathbf u_{D,h}+\mathcal U_h$.
Thus the discrete state velocity belongs to
$\mathcal U_{h,\mathbf u_d}$, whereas adjoint
velocities belong to $\mathcal U_h$.

For a fixed $\zeta>0$, define
\begin{equation}\label{eq:discrete-geometric-quantities}
 g_h(\phi_h):=\sqrt{|\nabla\phi_h|^2+\zeta^2},\qquad
 \mathbf n_h(\phi_h):=\frac{\nabla\phi_h}{g_h(\phi_h)},\qquad
 \mathbf P_{\tau,h}(\phi_h)
 :=\mathbf I-\mathbf n_h(\phi_h)\otimes\mathbf n_h(\phi_h).
\end{equation}
These quantities are evaluated elementwise. The diffuse interface measure is approximated by
$\delta_\varepsilon(\phi_h)g_h(\phi_h)\dx$.

For $\phi_h\in\mathcal W_h$, define
\begin{equation}\label{eq:discrete-state-form}
\begin{aligned}
a_h(\phi_h;\mathbf z_h,\mathbf v_h)
&:=
\int_\Omega 2\mu\mathcal E(\mathbf z_h):\mathcal E(\mathbf v_h)
      \mathcal H_f(\phi_h)\dx
+\int_\Omega\frac{\mu}{K}\mathbf z_h\cdot\mathbf v_h
      \mathcal H_p(\phi_h)\dx\\
&\quad+
\int_\Omega\frac{\alpha\mu}{\sqrt K}
\bigl(\mathbf P_{\tau,h}(\phi_h)\mathbf z_h\bigr)\cdot
\bigl(\mathbf P_{\tau,h}(\phi_h)\mathbf v_h\bigr)
\delta_\varepsilon(\phi_h)g_h(\phi_h)\dx.
\end{aligned}
\end{equation}
 Introduce the mixed bilinear form
\begin{equation}\label{eq:discrete-mixed-form}
\begin{aligned}
B_h\bigl(\phi_h;(\mathbf z_h,\pi_h),(\mathbf v_h,q_h)\bigr)
:=
a_h(\phi_h;\mathbf z_h,\mathbf v_h)
+b(\pi_h,\mathbf v_h)+b(q_h,\mathbf z_h)
\end{aligned}
\end{equation}
and the corresponding load functional
$L_h(\mathbf v_h)
:=
\int_\Omega\mathbf f\cdot\mathbf v_h\dx$.

The discrete state problem is: Find
$(\mathbf u_h,p_h)\in\mathcal U_{h,\mathbf u_d}\times\mathcal P_h$
such that
\begin{equation}\label{eq:discrete-state}
B_h\bigl(\phi_h;(\mathbf u_h,p_h),(\mathbf v_h,q_h)\bigr)
=
L_h(\mathbf v_h)
\qquad
\forall(\mathbf v_h,q_h)\in\mathcal U_h\times\mathcal P_h.
\end{equation}

For each fixed $\phi_h$, the form $a_h$ is continuous.
Since $\mathcal H_f(\phi_h)\geq\eta>0$ and the Darcy and BJS
terms are nonnegative, Korn's inequality yields 
\[
a_h(\phi_h;\mathbf v_h,\mathbf v_h)
\geq
2\mu\eta\|\mathcal E(\mathbf v_h)\|_{L^2(\Omega)}^2
\geq
c\mu\eta\|\mathbf v_h\|_{H^1(\Omega)^d}^2
\qquad
\forall\mathbf v_h\in\mathcal U_h.
\]
Together with the discrete inf--sup condition of the MINI pair
on shape-regular meshes, this coercivity gives existence and
uniqueness of the discrete state solution after subtracting
the prescribed finite element lifting.

The discrete objectives are
\begin{equation}\label{eq:discrete-objectives}
\begin{aligned}
\mathcal J_{1,h}(\mathbf u_h,\phi_h)
&:=
\int_\Omega 2\mu|\mathcal E(\mathbf u_h)|^2
      \mathcal H_f(\phi_h)\dx,\\
\mathcal J_{2,h}(\mathbf u_h,\phi_h)
&:=
\frac12\int_\Omega|\mathbf u_h-\mathbf u_{\rm obs}|^2\dx
+\gamma_{\rm per}\int_\Omega
  \delta_\varepsilon(\phi_h)g_h(\phi_h)\dx.
\end{aligned}
\end{equation}
For $j\in\{1,2\}$, define
\begin{equation}\label{eq:discrete-objective-state-derivatives}
\begin{aligned}
\ell_{1,h}(\mathbf v_h)
:=
\int_\Omega 4\mu\mathcal E(\mathbf u_h):\mathcal E(\mathbf v_h)
      \mathcal H_f(\phi_h)\dx,\quad
\ell_{2,h}(\mathbf v_h)
:=
\int_\Omega(\mathbf u_h-\mathbf u_{\rm obs})\cdot\mathbf v_h\dx.
\end{aligned}
\end{equation}
The discrete adjoint is the transpose of the
discrete state operator: Find
$(\mathbf y_{j,h},r_{j,h})\in\mathcal U_h\times\mathcal P_h$
such that
\begin{equation}\label{eq:discrete-adjoint}
B_h\bigl(\phi_h;(\mathbf v_h,q_h),(\mathbf y_{j,h},r_{j,h})\bigr)
=
\ell_{j,h}(\mathbf v_h)
\qquad
\forall(\mathbf v_h,q_h)\in\mathcal U_h\times\mathcal P_h.
\end{equation} 
The discrete adjoint problem is uniquely solvable because its
matrix is the transpose of the invertible discrete state matrix.
In particular, \eqref{eq:discrete-adjoint} is valid for both objectives;
only its right-hand side changes with $j$.

\subsection{Discrete design derivative and sensitivity field}
\label{subsec:discrete-sensitivity}

For $\psi_h\in\mathcal W_h$, differentiation of
\eqref{eq:discrete-geometric-quantities} gives
\begin{equation}\label{eq:discrete-geometric-derivatives}
 \begin{aligned}
 &Dg_h(\phi_h)[\psi_h]
 =\mathbf n_h(\phi_h)\cdot\nabla\psi_h,\quad
 \mathbf m_h[\psi_h]
 :=D\mathbf n_h(\phi_h)[\psi_h]
 =\frac{\mathbf P_{\tau,h}(\phi_h)\nabla\psi_h}
 {g_h(\phi_h)},\\
 &D\mathbf P_{\tau,h}(\phi_h)[\psi_h]
 =-\mathbf m_h[\psi_h]\otimes\mathbf n_h(\phi_h)
  -\mathbf n_h(\phi_h)\otimes\mathbf m_h[\psi_h].
 \end{aligned}
\end{equation}
With the state and adjoint velocities held fixed,
the design derivative of $a_h$ is
\begin{equation}\label{eq:discrete-a-derivative}
 \begin{aligned}
 &D_{\phi_h}a_h(\phi_h;\mathbf u_h,\mathbf y_{j,h})[\psi_h]\\
 ={}&
 -\int_\Omega2\mu(1-\eta)
 \mathcal E(\mathbf u_h):\mathcal E(\mathbf y_{j,h})
 \delta_\varepsilon(\phi_h)\psi_h\dx
 +\int_\Omega\frac{\mu}{K}(1-\eta)
 \mathbf u_h\cdot\mathbf y_{j,h}
 \delta_\varepsilon(\phi_h)\psi_h\dx\\
 &+\int_\Omega\frac{\alpha\mu}{\sqrt K}
 \Big[
 \big(D\mathbf P_{\tau,h}(\phi_h)[\psi_h]\mathbf u_h\big)
 \cdot
 \big(\mathbf P_{\tau,h}(\phi_h)\mathbf y_{j,h}\big)
 +\big(\mathbf P_{\tau,h}(\phi_h)\mathbf u_h\big)
 \cdot
 \big(D\mathbf P_{\tau,h}(\phi_h)[\psi_h]\mathbf y_{j,h}\big)
 \Big]
 \delta_\varepsilon(\phi_h)g_h(\phi_h)\dx\\
 &+\int_\Omega\frac{\alpha\mu}{\sqrt K}
 \big(\mathbf P_{\tau,h}(\phi_h)\mathbf u_h\big)\cdot
 \big(\mathbf P_{\tau,h}(\phi_h)\mathbf y_{j,h}\big)
 \Big[
 \delta_\varepsilon'(\phi_h)\psi_h g_h(\phi_h)
 +\delta_\varepsilon(\phi_h)
 \mathbf n_h(\phi_h)\cdot\nabla\psi_h
 \Big]\dx.
 \end{aligned}
\end{equation}
The explicit design derivatives of the objectives are
\begin{equation}\label{eq:discrete-objective-design-derivatives}
 \begin{aligned}
 \partial_{\phi_h}\mathcal J_{1,h}
 (\mathbf u_h,\phi_h)[\psi_h]
 &=-\int_\Omega2\mu(1-\eta)|\mathcal E(\mathbf u_h)|^2
 \delta_\varepsilon(\phi_h)\psi_h\dx,\\
 \partial_{\phi_h}\mathcal J_{2,h}
 (\mathbf u_h,\phi_h)[\psi_h]
 &=\gamma_{\rm per}\int_\Omega
 \Big[
 \delta_\varepsilon'(\phi_h)\psi_h g_h(\phi_h)
 +\delta_\varepsilon(\phi_h)
 \mathbf n_h(\phi_h)\cdot\nabla\psi_h
 \Big]\dx.
 \end{aligned}
\end{equation}

The discrete volume equality and its derivative, used only for
$j=1$, are
\begin{equation}\label{eq:discrete-volume}
 \mathcal C_h(\phi_h)
 :=\int_\Omega\mathcal H_f(\phi_h)\dx-V_0,
 \qquad
 D\mathcal C_h(\phi_h)[\psi_h]
 =-(1-\eta)\int_\Omega
 \delta_\varepsilon(\phi_h)\psi_h\dx.
\end{equation}
For $j=1$, let $\lambda\in\mathbb R$ and $\rho>0$, and introduce
the discrete effective multiplier
\begin{equation}\label{eq:discrete-effective-multiplier}
 \Lambda_h(\phi_h,\lambda)
 :=\lambda+\rho\mathcal C_h(\phi_h).
\end{equation}
For $j=2$, set $\Lambda_h(\phi_h,\lambda):=0$.

Differentiating the discrete state equation
\eqref{eq:discrete-state} and using the discrete adjoint
\eqref{eq:discrete-adjoint} gives the derivative of the
discrete reduced objective. Adding the augmented-volume
derivative for $j=1$ defines the functional
$\mathfrak D_{j,h}(\phi_h,\lambda)\in\mathcal W_h'$ by
\begin{equation}\label{eq:discrete-sensitivity}
 \begin{aligned}
 \left\langle
 \mathfrak D_{j,h}(\phi_h,\lambda),\psi_h
 \right\rangle
 &:=
 \partial_{\phi_h}\mathcal J_{j,h}
 (\mathbf u_h,\phi_h)[\psi_h]
 -D_{\phi_h}a_h
 (\phi_h;\mathbf u_h,\mathbf y_{j,h})[\psi_h]\\
 &\quad-\Lambda_h(\phi_h,\lambda)(1-\eta)
 \int_\Omega\delta_\varepsilon(\phi_h)\psi_h\dx
 \qquad\forall\psi_h\in\mathcal W_h.
 \end{aligned}
\end{equation}
For $j=2$, the last term in
\eqref{eq:discrete-sensitivity} vanishes.

The raw sensitivity may contain element-scale oscillations,
particularly in the diffuse BJS contribution.
We therefore compute an $H^1$ Riesz representative
$\widetilde G_{j,h}\in\mathcal W_h$: for a prescribed
filter length $\ell_s>0$,
\begin{equation}\label{eq:filtered-sensitivity}
 \int_\Omega\widetilde G_{j,h}\psi_h\dx
 +\ell_s^2\int_\Omega
 \nabla\widetilde G_{j,h}\cdot\nabla\psi_h\dx
 =
 \left\langle
 \mathfrak D_{j,h}(\phi_h,\lambda),\psi_h
 \right\rangle
 \quad\forall\psi_h\in\mathcal W_h.
\end{equation}
This step smooths the sensitivity associated with the objective
and augmented-volume terms, while the diffusion term is treated
implicitly in the level-set update
\eqref{eq:phi-update-discrete}.
In the computations, $\ell_s$ is chosen as a small multiple
of the mesh size.

\subsection{Optimization algorithm}
\label{subsec:discrete-reaction-diffusion}

Let $\Delta s>0$ be the artificial step size, let $\kappa>0$ be the
mobility, and let $\tau>0$ be the fictitious-interface regularization
parameter from \eqref{eq:yamada-energy}.
The following describes the baseline algorithm.
Both $\kappa$ and $\tau$ are kept fixed throughout the optimization.
Given $\phi_h^n$ and the filtered sensitivity
$\widetilde G_{j,h}^n$, we use the following semi-implicit
reaction--diffusion update with filtered sensitivity:
find $\phi_h^{n+1}\in W_h$ such that
\begin{equation}\label{eq:phi-update-discrete}
 \begin{aligned}
 \frac{1}{\kappa\Delta s}
 \int_\Omega(\phi_h^{n+1}-\phi_h^n)\psi_h\dx
 +\tau\int_\Omega\nabla\phi_h^{n+1}\cdot\nabla\psi_h\dx
 =-\int_\Omega\widetilde G_{j,h}^n\psi_h\dx
 \quad\forall\psi_h\in\mathcal W_h.
 \end{aligned}
\end{equation}
The homogeneous Neumann condition
$\nabla\phi_h^{n+1}\cdot\mathbf n_{\partial\Omega}=0$ is natural in
\eqref{eq:phi-update-discrete}.  In matrix notation,
\begin{equation}\label{eq:phi-update-matrix}
 \bigl(M_h+\kappa\tau\Delta s\,K_h\bigr)\Phi^{n+1}
 =M_h\Phi^n-\kappa\Delta s\,M_h\widetilde G^n,
\end{equation}
where $M_h$ and $K_h$ are the scalar mass and stiffness matrices.

For $j=1$, the augmented-Lagrangian multiplier is held fixed during
$N_\lambda\ge1$ inner design steps.  After the level-set update and any
scheduled reinitialization, it is updated by the relaxed dual step
\begin{equation}\label{eq:multiplier-update}
 \lambda^{n+1}
 =
 \begin{cases}
 \lambda^n+\sigma\mathcal C_h(\phi_h^{n+1}),
 &(n+1)\equiv0\pmod{N_\lambda}
   \ \text{and}\ 
   |\mathcal C_h(\phi_h^{n+1})|>\mathrm{tol}_{\mathcal C},\\
 \lambda^n,&\text{otherwise},
 \end{cases}
 \qquad \sigma>0.
\end{equation}
The penalty parameter $\rho$ remains fixed, whereas $\sigma$ controls the
multiplier relaxation.  Updating the multiplier only after several inner
steps avoids changing the augmented functional at every design update.
For $j=2$, no multiplier update is performed and $\lambda=\rho=0$.

Every $N_{\rm re}$ level-set updates, the fast marching method
\eqref{eq:fast-marching-reinitialization} is applied.
For convergence
monitoring, we use the fixed-penalty merit function
\begin{equation}\label{eq:fixed-penalty-merit}
 \mathcal M_{j,\rho,h}(\phi_h)
 :=\mathcal J_{j,h}(\mathbf u_h(\phi_h),\phi_h)
   +\frac{\rho}{2}\mathcal C_h(\phi_h)^2.
\end{equation}
For $j=2$, this reduces to
$\mathcal M_{2,0,h}=\mathcal J_{2,h}$.
The iteration is terminated when
\begin{equation}\label{eq:merit-stopping-criterion}
 \frac{\left|\mathcal M_{j,\rho,h}(\phi_h^n)
 -\mathcal M_{j,\rho,h}(\phi_h^{n-1})\right|}
 {\max\left\{1,\left|\mathcal M_{j,\rho,h}(\phi_h^{n-1})\right|\right\}}
 \le \mathrm{tol}_{\mathcal M}.
\end{equation}

\begin{algorithm}[htbp]
\caption{Reaction--diffusion level-set optimization}
\label{alg:optimization}
\SetAlgoLined
\KwData{$j$, $\phi_h^0$, $\tau$, $\kappa$, $\Delta s$, $\ell_s$,
$N_{\rm re}$, $N_{\max}$, and $\mathrm{tol}_{\mathcal M}$; additionally
$\lambda^0$, $V_0$, $\rho$, $\sigma$, $N_\lambda$,
 and $\mathrm{tol}_{\mathcal C}$ for $j=1$.}

\For{$n=0,1,\ldots,N_{\max}-1$}{
  Solve the state problem \eqref{eq:discrete-state} for
  $(\mathbf u_h^n,p_h^n)$ and evaluate
  $\mathcal M_{j,\rho,h}(\phi_h^n)$\;

  \If{$n\ge1$ and \eqref{eq:merit-stopping-criterion} holds}{
    \textbf{break}\;
  }

  Solve the adjoint problem \eqref{eq:discrete-adjoint} for
  $(\mathbf y_{j,h}^n,r_{j,h}^n)$\;

  Compute the design sensitivity from
  \eqref{eq:discrete-sensitivity} and
  \eqref{eq:filtered-sensitivity}\;

  Update $\phi_h^{n+1}$ using \eqref{eq:phi-update-discrete}\;

  Reinitialize $\phi_h^{n+1}$ by the fast marching method every
  $N_{\rm re}$ iterations\;

  \If{$j=1$}{
    Update $\lambda^{n+1}$ using \eqref{eq:multiplier-update}\;
  }
}
\end{algorithm}









\section{Numerical examples}

\subsection{Numerical assessment of the state solver and objective derivatives}

Before applying the state solver to topology optimization, 
we assess its spatial accuracy using a manufactured solution
and examine the response of a fixed free-flow--porous
configuration to the permeability.
Hereafter, $N$
denotes the nominal mesh resolution, with target spacing $h=1/N$; for a
Cartesian mesh, it is exactly the number of uniform subdivisions in each
coordinate direction. Convergence rates are computed using the measured
maximum element diameter $h_{\max}$.

\subsubsection{Numerical accuracy of the regularized Stokes--Darcy solver}

\begin{figure}[!htbp]
 \begin{subfigure}[h]{0.32\textwidth}
  \centering
  \includegraphics[width=1.6in]{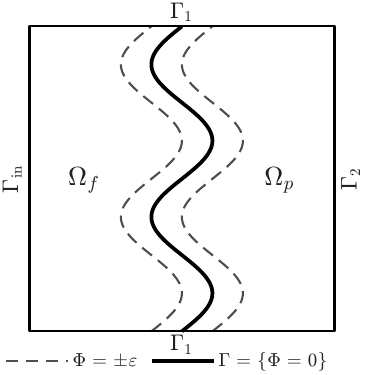}
  \caption{}
 \end{subfigure}
 \begin{subfigure}[h]{0.32\textwidth}
  \centering
  \includegraphics[width=1.8in]{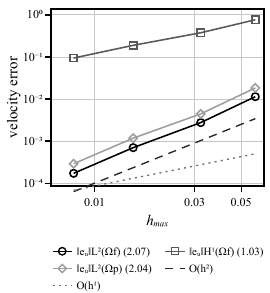}
  \caption{}
 \end{subfigure}
 \begin{subfigure}[h]{0.32\textwidth}
  \centering
  \includegraphics[width=1.8in]{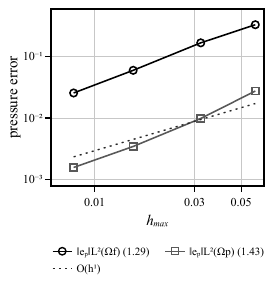}
  \caption{}
 \end{subfigure}
 \caption{Manufactured-solution verification: (a) computational geometry,
 (b) velocity errors, and (c) pressure errors. In (a), the solid and dashed
 curves denote $\Phi=0$ and $\Phi=\pm\varepsilon$, respectively.
 Parenthetical values in (b)--(c) are the fitted rates based on the four
 measured values of $h_{\max}$.}
 \label{fig:state-solver-accuracy}
\end{figure}

We consider the square domain $\Omega=(0,1)^2$ shown in
Fig.~\ref{fig:state-solver-accuracy}(a), where
$\Gamma_{\rm in}=\{0\}\times(0,1)$,
$\Gamma_2=\{1\}\times(0,1)$, and
$\Gamma_1=(0,1)\times\{0,1\}$. The two regions and their S-shaped interface
are defined by
\begin{equation}
 \Phi(x,y)=x-\frac12-0.1\sin(4\pi y),\qquad
 \Omega_f=\{\Phi<0\},\quad \Omega_p=\{\Phi>0\},\quad
 \Gamma=\{\Phi=0\}.
 \label{eq:state-verification-level-set}
\end{equation}
Homogeneous velocity data are imposed on $\partial\Omega$. Let
$\Psi=\sin^2(\pi x)\sin^2(\pi y)$ and prescribe
\begin{equation}
 \mathbf U=(\partial_y\Psi,-\partial_x\Psi)^{\top}
 =\begin{pmatrix}
 \pi\sin^2(\pi x)\sin(2\pi y)\\[-1mm]
 -\pi\sin(2\pi x)\sin^2(\pi y)
 \end{pmatrix},\qquad
 P=\sin(2\pi x)\sin(2\pi y).
 \label{eq:state-verification-exact}
\end{equation}
These fields satisfy $\nabla\!\cdot\mathbf U=0$,
$\mathbf U=\mathbf0$ on $\partial\Omega$, and
$\int_\Omega P\,{\rm d}x=0$.
The manufactured load $\mathbf f$ is obtained by
substituting $(\mathbf U,P)$ and the analytical level-set function $\Phi$
into the regularized momentum operator. The geometric factors in the load
use the same $\zeta$-regularized definitions as the discrete operator.

The global MINI pair $[P_1\oplus b_T]^2/P_1$ and the interpolated level-set
function $\phi_h=I_h\Phi$ are used. Throughout the refinement
sequence, $\mu=K=\alpha=1$, $\varepsilon=0.1$, $\eta=10^{-3}$, and
$\zeta=10^{-8}$. All volume integrals are evaluated with quadrature order
10. The four nominal resolutions are $N=24,48,96,$ and $192$, and the rates
are fitted against the corresponding measured values of $h_{\max}$.
The fitted rates of the velocity errors are $2.07$ in
$L^2(\Omega_f)$, $1.03$ in $H^1(\Omega_f)$, and $2.04$ in
$L^2(\Omega_p)$, consistent with the expected MINI velocity behavior; see
Fig.~\ref{fig:state-solver-accuracy}(b). The pressure errors decrease under
refinement, with fitted $L^2$ rates $1.29$ in $\Omega_f$ and $1.43$ in
$\Omega_p$; see Fig.~\ref{fig:state-solver-accuracy}(c). These pressure rates
are empirical results for the present smooth manufactured solution and are
not asserted as general a priori estimates.

\subsubsection{Effect of permeability on pressure-driven seepage}

We retain the geometry in Eq.~\eqref{eq:state-verification-level-set},
set the body force to zero, and impose the normal tractions
$\mathbf T_\varepsilon(\mathbf u,p;\phi)\mathbf n_{\partial\Omega}
=-p_{\rm in}\mathbf n_{\partial\Omega}$
on $\Gamma_{\rm in}$ and
$\mathbf T_\varepsilon(\mathbf u,p;\phi)\mathbf n_{\partial\Omega}
=-p_{\rm out}\mathbf n_{\partial\Omega}$
on $\Gamma_2$, with $p_{\rm in}=1$ and $p_{\rm out}=0$;
no-slip conditions are imposed on $\Gamma_1$.
Here
$
\mathbf T_\varepsilon(\mathbf u,p;\phi)
:=-p\mathbf I+2\mu\mathcal H_f(\phi)\mathcal E(\mathbf u)
$
is the regularized stress, and $\mathbf n_{\partial\Omega}$
is the outward unit normal to $\partial\Omega$.
The velocity test space is
$\{\mathbf v\in H^1(\Omega)^d:
\mathbf v=0\text{ on }\Gamma_1\}$,
and the pressure space is $L^2(\Omega)$ without a zero-mean constraint.
The prescribed tractions enter the weak momentum equation through the
boundary load
$
-\int_{\Gamma_{\rm in}}
p_{\rm in}\mathbf n_{\partial\Omega}\cdot\mathbf v\ds
-\int_{\Gamma_2}
p_{\rm out}\mathbf n_{\partial\Omega}\cdot\mathbf v\ds.
$
The finite element spaces are adjusted accordingly.

We use a Cartesian mesh with $N=160$ subdivisions in each
coordinate direction and set $h=1/N$. The parameters are
$\mu=\alpha=1$, $\eta=10^{-8}$, $\varepsilon=0.75h$, and
$\zeta=10^{-8}$, and all volume and boundary integrals are evaluated with
quadrature order 8. Only the permeability is varied:
$K\in\{1,10^{-2},10^{-4}\}$.
We use the outlet flux
$Q_{\rm out}=\int_{\Gamma_2}\mathbf u_h\!\cdot
\mathbf n_{\partial\Omega}\ds$
to quantify the transmitted flow.

\begin{table}[!htbp]
\centering
\small
\caption{Outlet flux for the pressure-driven problem on the $N=160$ mesh.}
\label{tab:state-permeability-flux}
\begin{tabular}{ccc}
\toprule
$K$ & $Q_{\rm out}$ & $Q_{\rm out}/Q_{\rm out}(K=1)$\\
\midrule
$1$ &
$1.7570\!\times\!10^{-1}$ &
$1$\\
$10^{-2}$ &
$1.6579\!\times\!10^{-2}$ &
$9.4357\!\times\!10^{-2}$\\
$10^{-4}$ &
$2.1612\!\times\!10^{-4}$ &
$1.2300\!\times\!10^{-3}$\\
\bottomrule
\end{tabular}
\end{table}

\begin{figure}[!htbp]
 \begin{subfigure}[h]{0.33\textwidth}
  \centering
  \includegraphics[width=1.75in]{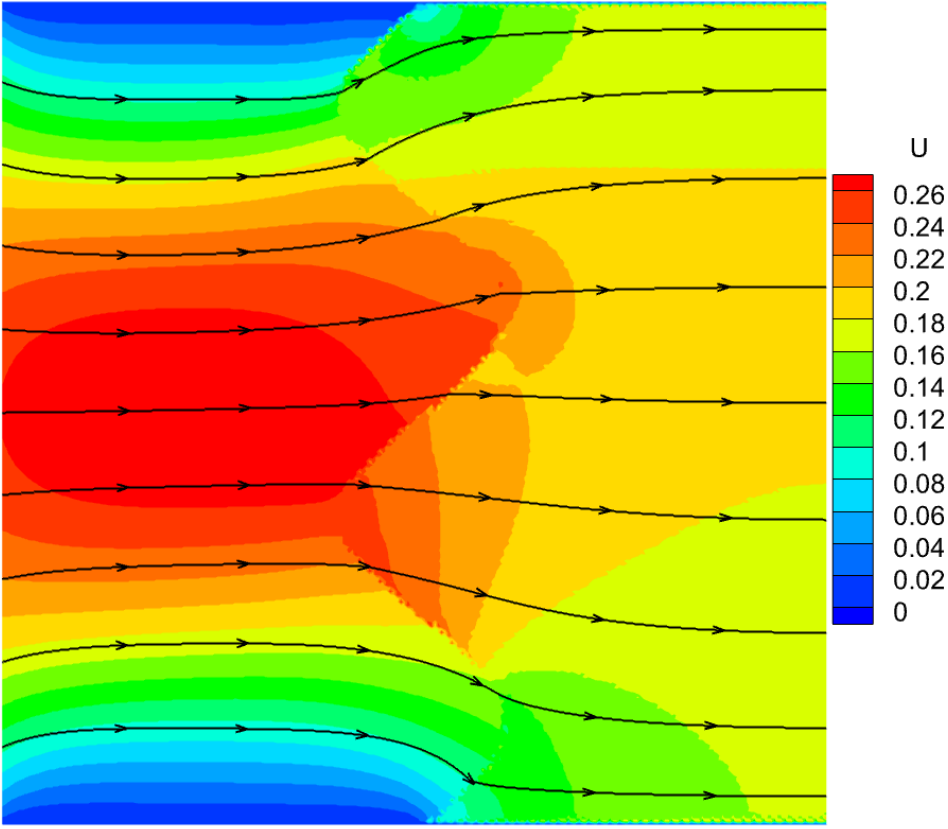}
  \caption{$K=1$}
 \end{subfigure}
 \begin{subfigure}[h]{0.33\textwidth}
  \centering
  \includegraphics[width=1.8in]{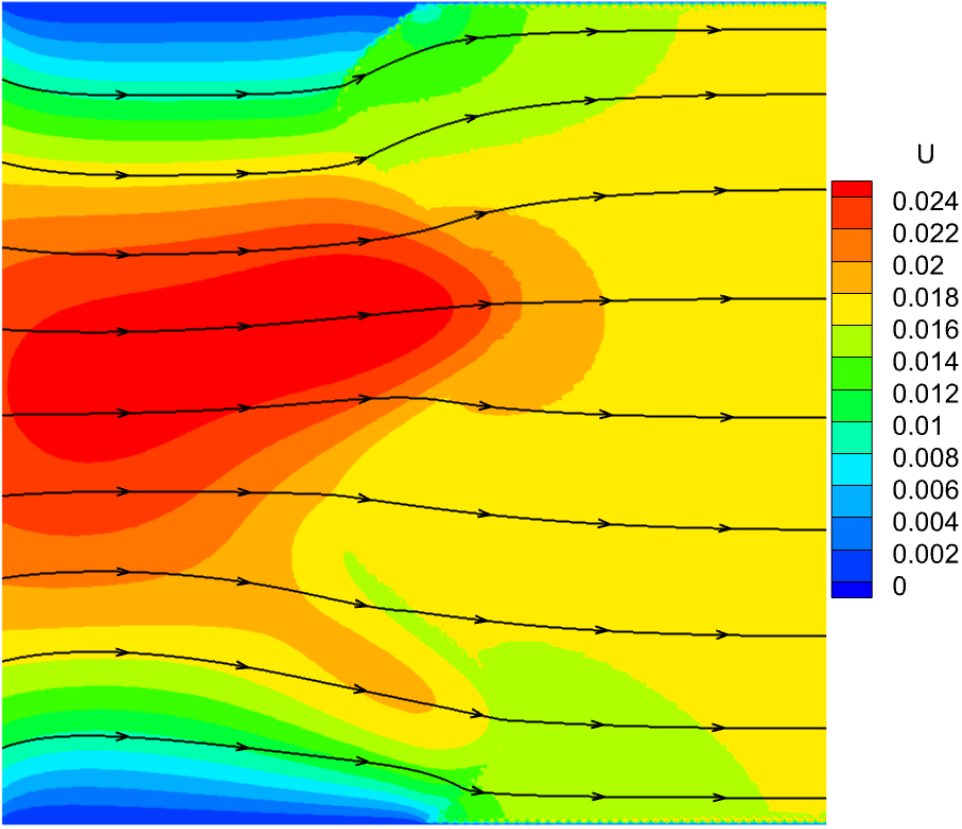}
  \caption{$K=10^{-2}$}
 \end{subfigure}
 \begin{subfigure}[h]{0.33\textwidth}
  \centering
  \includegraphics[width=1.85in]{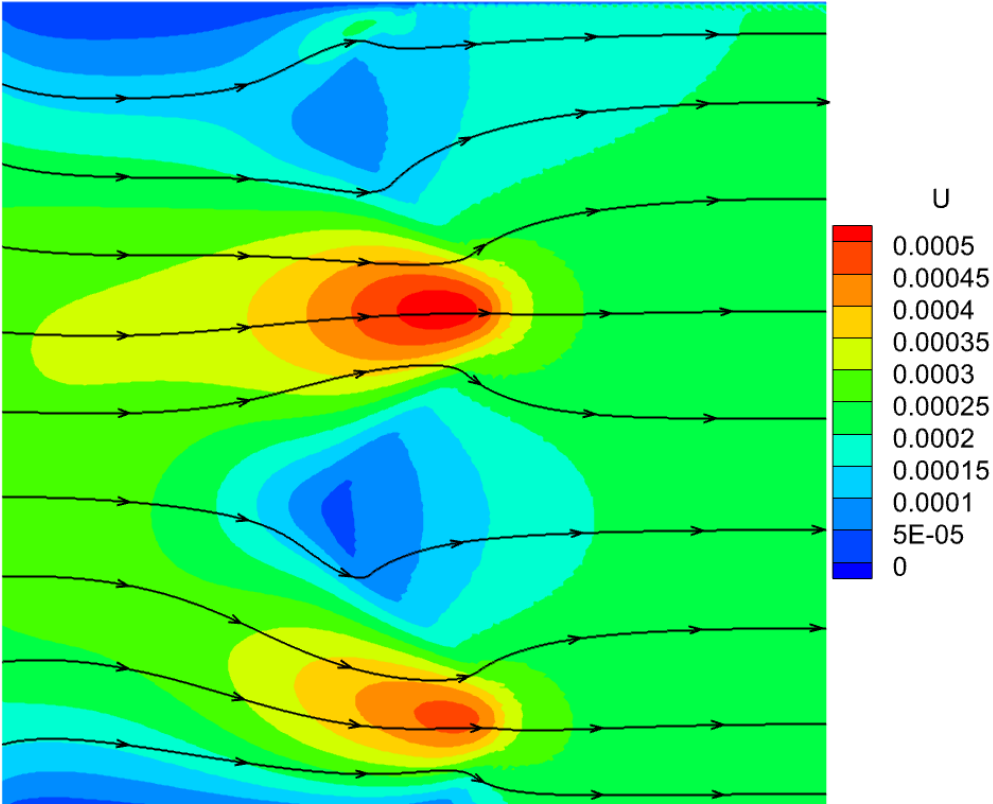}
  \caption{$K=10^{-4}$}
 \end{subfigure}
 \caption{Pressure-driven flow for three permeabilities:
 (a) $K=1$, (b) $K=10^{-2}$, and (c) $K=10^{-4}$.
 Colors indicate velocity magnitude, and the overlaid streamlines indicate
 the flow direction.}
 \label{fig:state-permeability-fields}
\end{figure}

Table~\ref{tab:state-permeability-flux} and
Fig.~\ref{fig:state-permeability-fields} show that the transmitted flow
decreases as $K$ is reduced in the fixed regularized model. The outlet flux
at $K=10^{-4}$ is approximately $1.23\times10^{-3}$ of that at $K=1$, a
reduction of nearly three orders of magnitude. The prescribed
traction difference provides a common driving condition for comparing the
response of the fixed configuration to changes in $K$. Since varying $K$
changes both the Darcy coefficient $\mu/K$ and the diffuse BJS coefficient
$\alpha\mu/\sqrt{K}$, the comparison represents their combined effect in the
regularized model.

\subsubsection{Taylor verification of the objective derivatives}
\label{subsec:taylor-verification}

We test the derivatives of \eqref{eq:discrete-objectives} on fixed meshes,
using MINI velocity--pressure elements and piecewise linear level sets.
Reinitialization, sensitivity filtering, augmented-volume
terms, and level-set updates are not used; all data, quadrature rules, and
regularization parameters remain fixed.

For $\phi_h^s=\phi_h+s\psi_h$, recompute the state from
\eqref{eq:discrete-state} and write
$\mathcal J_{j,h}^s:=\mathcal J_{j,h}(\mathbf u_h(\phi_h^s),\phi_h^s)$.
We evaluate the directional derivative
$d_{j,h}:=\left.\frac{\mathrm d}{\mathrm ds}
\mathcal J_{j,h}^{s}\right|_{s=0}$
using the adjoint formula without the augmented-volume term.
We measure $R_{1,j}(s):=|\mathcal J_{j,h}^s-\mathcal J_{j,h}^0-sd_{j,h}|$
and compute $p_{1,j}(s):=\log_2(R_{1,j}(2s)/R_{1,j}(s))$.
For sufficiently small $s$, a consistent implementation is
expected to give $R_{1,j}(s)=O(s^2)$.
The perturbations need not preserve the volume or Eikonal constraints.

The two-dimensional test uses $\Omega=(0,1)^2$ with 2048 triangles
and 1089 vertices. With $I_h$ denoting nodal interpolation, take
$\phi_h=I_h[x-\tfrac12-0.04\sin(4\pi y)]$ and the
$L^\infty$-normalized direction $I_h[\sin(\pi x)\sin(4\pi y)]$.
Set $\mu=K=\alpha=1$, $\varepsilon=0.1$, $\eta=10^{-3}$,
$\zeta=10^{-8}$, and $\gamma_{\rm per}=0.02$.
The velocity boundary data are homogeneous, and
$\mathbf u_{\rm obs}=0.8(\partial_y\Psi,-\partial_x\Psi)$ with
$\Psi=\sin^2(\pi x)\sin^2(\pi y)$.
The manufactured force is fixed at the base geometry.
Both objectives are tested; the explicit perimeter contribution is
approximately $-5.65644\times10^{-2}$ and is therefore nonzero.

The three-dimensional test checks $\mathcal J_{1,h}$ on $(0,1)^3$
with 3072 tetrahedra and 729 vertices. The base level-set
field is constructed from a prescribed union of seven spheres, followed by
one signed-distance reconstruction. Boundary velocities have balanced
fluxes and remain fixed. We set $\mu=\alpha=1$,
$K=10^{-4}$, $\varepsilon=0.15$, $\eta=10^{-3}$, and
$\zeta=10^{-4}$.
The direction is
$\psi_h=I_h[\sin(\pi x)\sin(\pi y)\sin(\pi z)
(1+0.3\cos(2\pi y))]$.
All geometric derivatives are retained, including both tangential
factors in the BJS term.

\begin{table}[htbp]
 \centering
 \captionsetup{font=footnotesize}
 \caption{Taylor remainders $R_{1,j}$ and observed orders
 $p_{1,j}$ in two and three dimensions.}
 \label{tab:taylor-2d-3d}
 \setlength{\tabcolsep}{3pt}
 \begin{tabular}{ccccc ccc}
 \toprule
 \multicolumn{5}{c}{Two dimensions}
 &\multicolumn{3}{c}{Three dimensions}\\
 \cmidrule(lr){1-5}\cmidrule(lr){6-8}
 $s_{\rm 2D}$ & $R_{1,1}$ & $p_{1,1}$ & $R_{1,2}$ & $p_{1,2}$
 &$s_{\rm 3D}$ & $R_{1,1}$ & $p_{1,1}$\\
 \midrule
$3.125\times10^{-4}$ & $3.541\times10^{-5}$ & -- & $2.266\times10^{-6}$ & -- & $5.000\times10^{-4}$ & $6.559\times10^{-6}$ & -- \\
$1.563\times10^{-4}$ & $8.863\times10^{-6}$ & 1.998 & $5.652\times10^{-7}$ & 2.004 & $2.500\times10^{-4}$ & $1.643\times10^{-6}$ & 1.997 \\
$7.813\times10^{-5}$ & $2.217\times10^{-6}$ & 1.999 & $1.411\times10^{-7}$ & 2.002 & $1.250\times10^{-4}$ & $4.111\times10^{-7}$ & 1.999 \\
$3.906\times10^{-5}$ & $5.544\times10^{-7}$ & 2.000 & $3.526\times10^{-8}$ & 2.001 & $6.250\times10^{-5}$ & $1.028\times10^{-7}$ & 1.999 \\
$1.953\times10^{-5}$ & $1.386\times10^{-7}$ & 2.000 & $8.812\times10^{-9}$ & 2.000 & $3.125\times10^{-5}$ & $2.571\times10^{-8}$ & 2.000 \\
 \bottomrule
 \end{tabular}
\end{table}

Table~\ref{tab:taylor-2d-3d} reports the five smallest
perturbation sizes retained from each halving sequence. The
remainders in all three tests exhibit second-order decay over the reported
perturbation range, consistent with the implemented objective derivatives.
These are fixed-mesh derivative-consistency tests, rather
than spatial discretization or optimization convergence tests.

\subsection{Viscous-energy minimization}

\subsubsection{Multi-inlet and multi-outlet channel design}
\begin{figure}[!htbp]
  \begin{subfigure}[hb]{0.5\textwidth}
    \centering
    \includegraphics[width=1.75in]
    {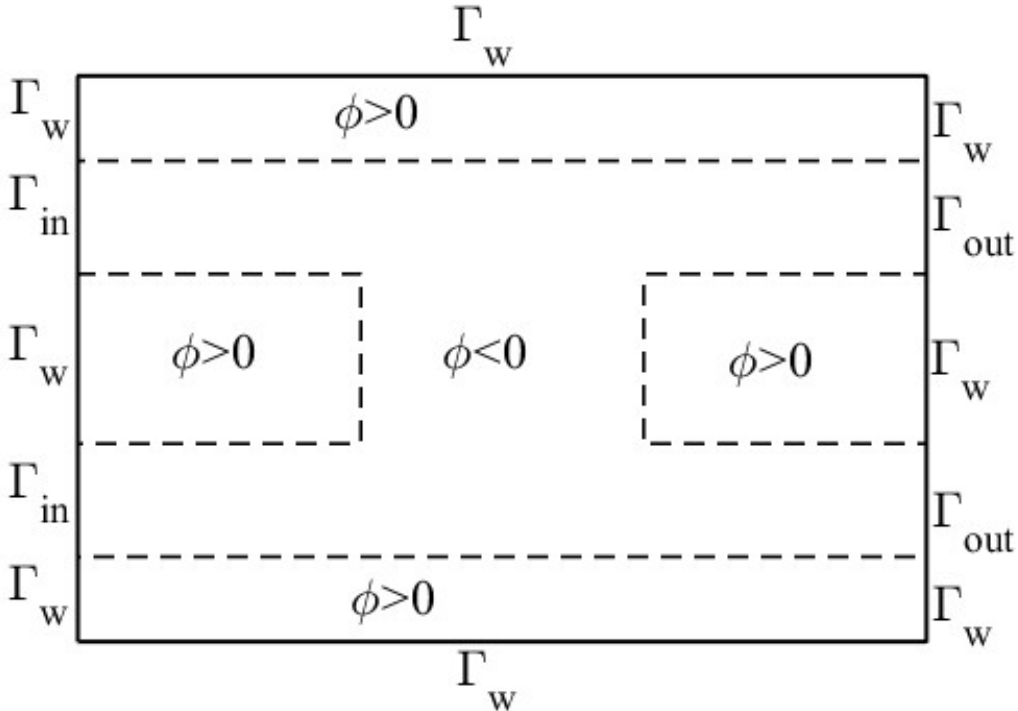}
    \caption{Domain and boundary conditions}
  \end{subfigure}
  \begin{subfigure}[hb]{0.5\textwidth}
    \centering
    \includegraphics[width=1.75in]
    {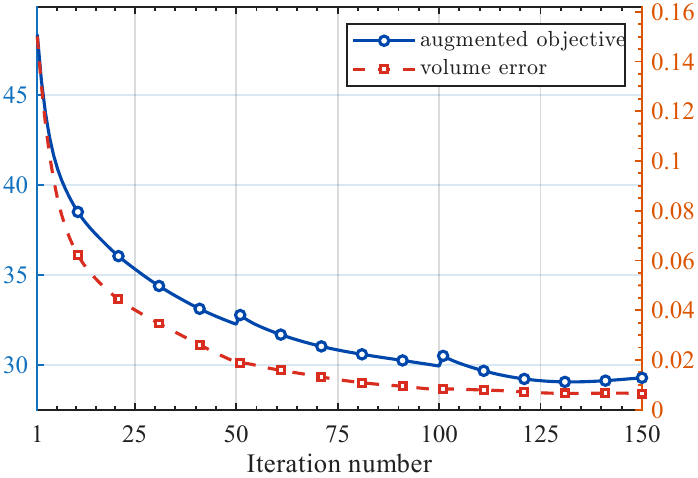}
    \caption{Optimization histories}
  \end{subfigure}
  \begin{subfigure}[h]{0.32\textwidth}
    \centering
    \includegraphics[width=1.6in]
    {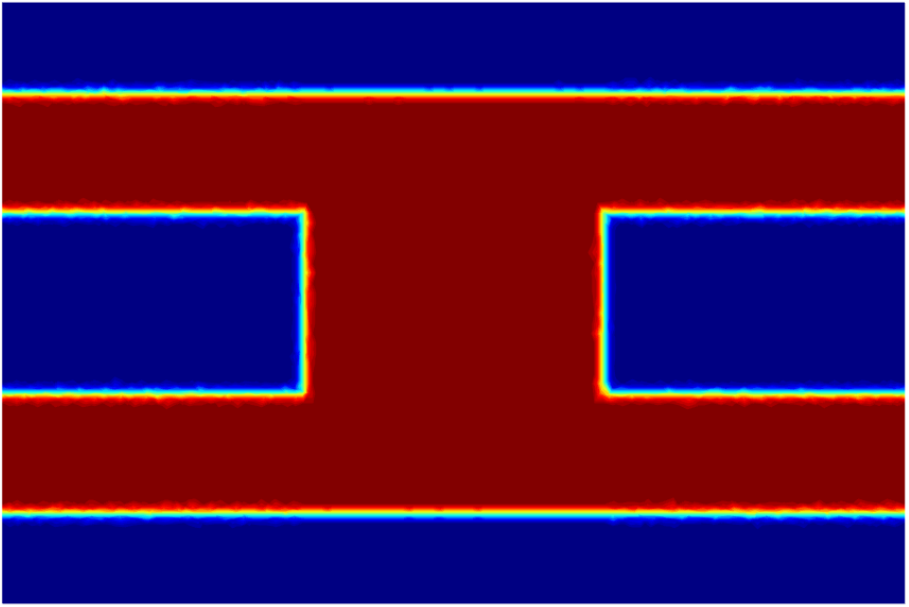}
    \caption{Initial Heaviside function}
  \end{subfigure}
  \begin{subfigure}[h]{0.32\textwidth}
    \centering
    \includegraphics[width=1.8in]
    {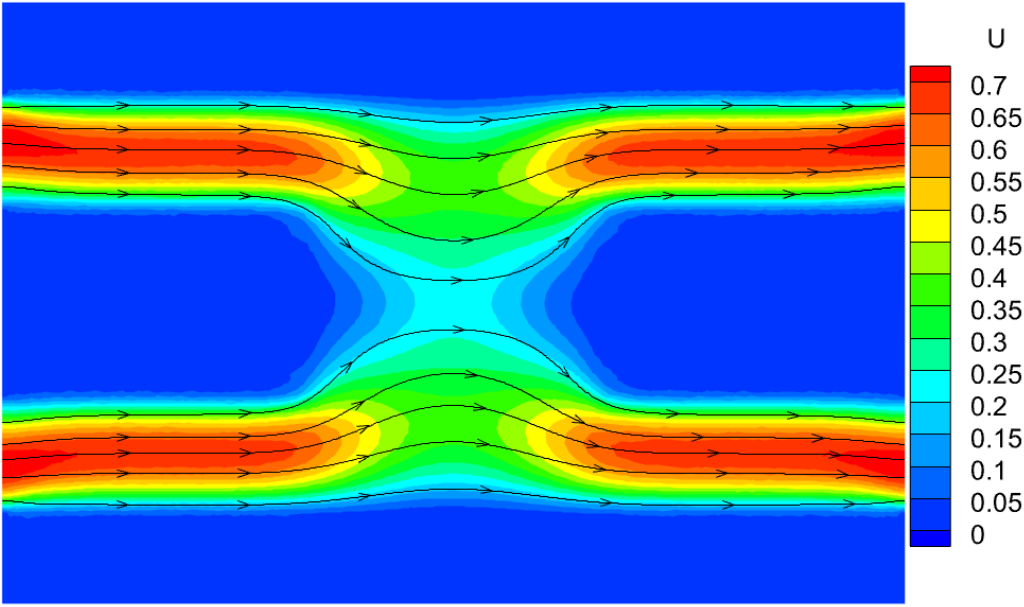}
    \caption{Velocity magnitude of the initial design}
  \end{subfigure}
  \begin{subfigure}[h]{0.32\textwidth}
    \centering
    \includegraphics[width=1.8in]
    {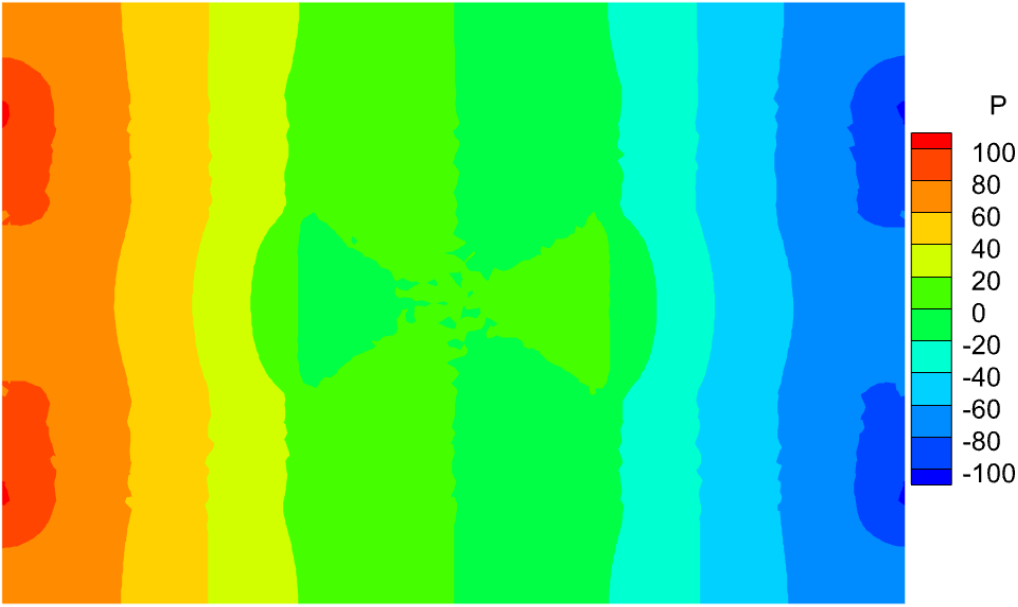}
    \caption{Pressure of the initial design}
  \end{subfigure}

  \begin{subfigure}[h]{0.32\textwidth}
    \centering
    \includegraphics[width=1.6in]
    {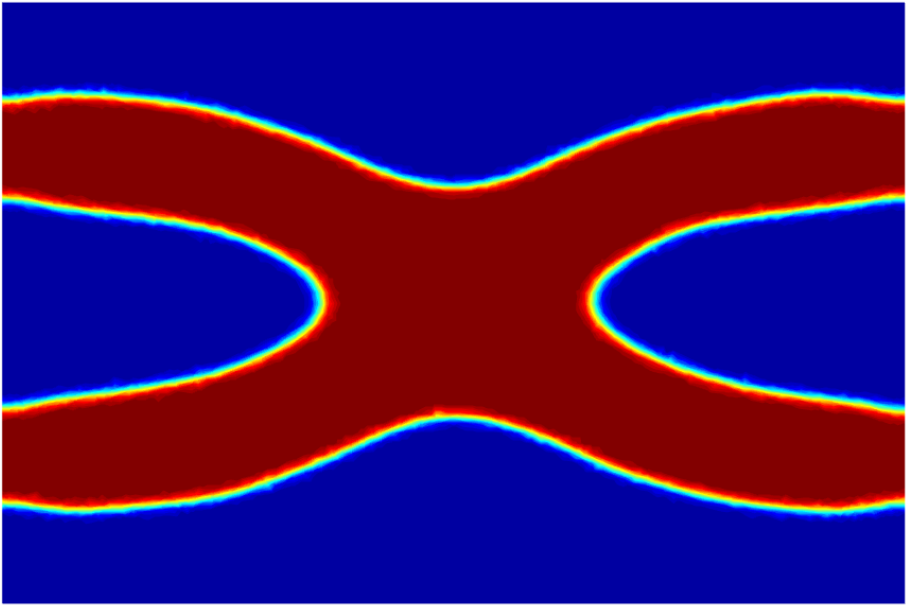}
    \caption{Optimized Heaviside function}
  \end{subfigure}
  \begin{subfigure}[h]{0.32\textwidth}
    \centering
    \includegraphics[width=1.8in]
    {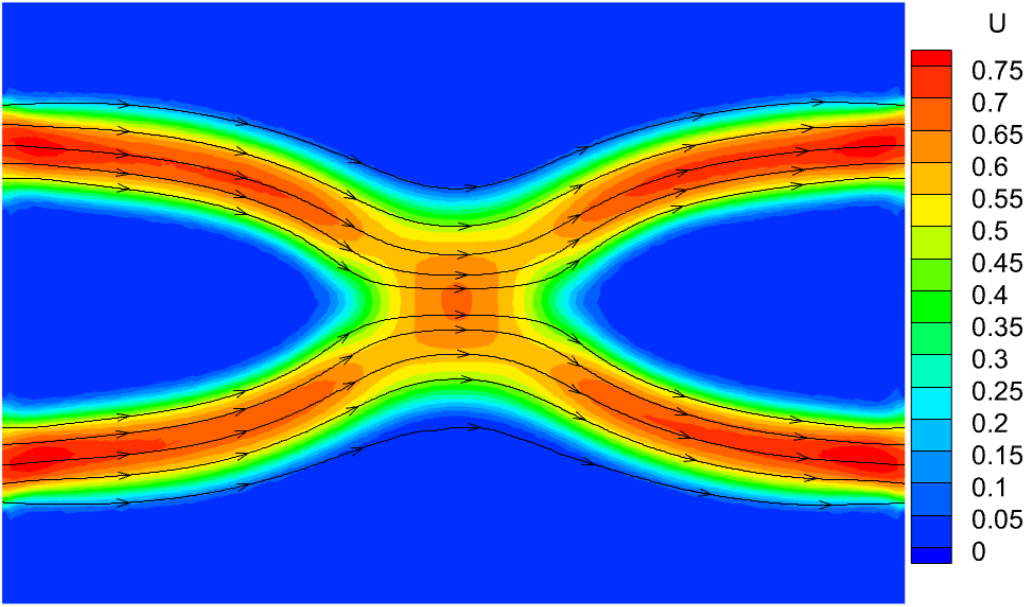}
    \caption{Velocity magnitude of the optimized design}
  \end{subfigure}
  \begin{subfigure}[h]{0.32\textwidth}
    \centering
    \includegraphics[width=1.8in]
    {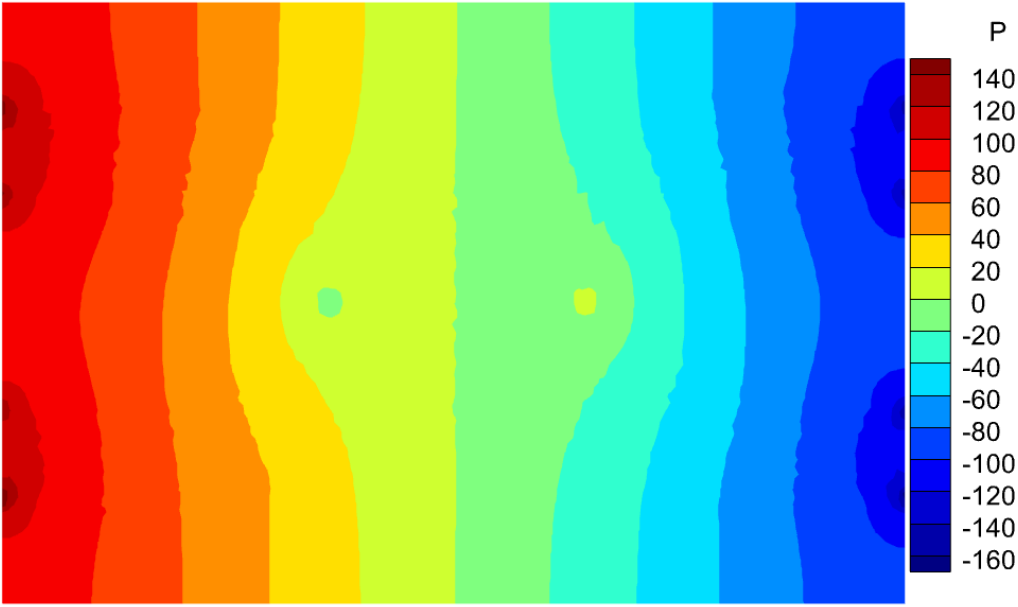}
    \caption{Pressure of the optimized design}
  \end{subfigure}

  \caption{Shape optimization of the two-channel configuration:
  (a) computational domain and boundary conditions;
  (b) histories of the augmented objective and the absolute
  diffuse-volume residual;
  (c)--(e) initial Heaviside function, velocity magnitude, and
  pressure, respectively; and
  (f)--(h) corresponding optimized Heaviside function, velocity
  magnitude, and pressure, respectively.}
  \label{TwoChannels}
\end{figure}

	{\r Perfusable channel networks can be fabricated within hydrogels for
	tissue engineering \cite{miller2012rapid}. Motivated by such systems,
	we consider an idealized channel design problem with two inlets and
	two outlets.} We minimize the
	viscous-energy objective \eqref{eq:objective1} for the coupled
	Stokes--Darcy system in
	$\Omega=(0,1.5)\times(-0.5,0.5)$, subject to a prescribed diffuse
	fluid volume. As shown in Fig.~\ref{TwoChannels}(a), the boundary is
	decomposed as
	$\partial\Omega=\Gamma_{\rm in}\cup\Gamma_w\cup\Gamma_{\rm out}$,
	where
	$
	\Gamma_{\rm in}
	=\{(0,y):0.15\leq |y|\leq0.35\}$,
	$
	\Gamma_{\rm out}
	=\{(1.5,y):0.15\leq |y|\leq0.35\}.
	$
	The prescribed velocity on the inlet and outlet segments is
	$
	\mathbf u_d(x,y)
	=
	\left(
	300\left|
	(y-0.35)(y-0.15)(y+0.35)(y+0.15)
	\right|,0
	\right)^{\top},
	$
	and the no-slip condition $\mathbf u=\mathbf0$ is imposed on
	$\Gamma_w$.

The physical parameters are
$\mu=1$, $K=10^{-4}$, and $\alpha=1$. The diffuse-interface parameters
are $\eta=10^{-3}$, $\zeta=10^{-8}$, and
$\varepsilon=h_{\max}$. For the reaction--diffusion update, we use
$\Delta s=0.01$, $\kappa=5\times10^{-4}$, $\tau=5$, and the
sensitivity-filter length $\ell_s=3h_{\max}$. The augmented-Lagrangian
parameters are
$\lambda^0=120$, $\rho=600$, $\sigma=0.02$, and $N_\lambda=5$.
The prescribed diffuse fluid volume is
$
V_0
=
0.8\int_\Omega
\mathcal H_f(\phi_h^0)\dx$.
The reaction displacement in one update is limited to
$0.2\varepsilon$. During the final $30\%$ of the prescribed iterations,
the reaction step is multiplied by a cosine relaxation factor that
decreases from one to zero. The level-set function is reinitialized
by the fast marching method every $50$ iterations. We use
$N_{\max}=150$ iterations. The mesh generated with the refinement
parameter $n_h=3$ contains $12{,}308$ triangles and $6{,}305$ vertices.
The initial level-set function is
$
\phi_0(x,y)
=
\min\!\left\{
\min\bigl(|y-0.25|-0.1,\ |y+0.25|-0.1\bigr),
\max\bigl(|x-0.75|,|y|\bigr)-0.25
\right\}.
$
This initialization contains two horizontal flow passages and a
central connecting region.

Figures~\ref{TwoChannels}(b)--(h) present the numerical results.
Figure~\ref{TwoChannels}(b) reports the augmented objective and the
absolute diffuse-volume residual over the prescribed $150$
iterations. The augmented objective decreases from $48.3547$ to
$29.2801$, while the absolute volume residual decreases from
$0.150055$ to $6.49681\times10^{-3}$. The computation terminates upon
reaching $N_{\max}$. Figures~\ref{TwoChannels}(c)--(h) compare the
diffuse fluid indicator, velocity magnitude, and pressure for the
initial and optimized designs. The optimized free-flow region forms
connected paths between the inlet and outlet segments and narrows in
the central part of the domain. The larger velocity magnitudes occur
primarily within these paths, while the velocity remains smaller in
the porous background.

\subsubsection{L-shaped channel design}

\begin{figure*}[htbp]
    \centering

    \begin{tabular}{c c c c}
        \multirow{2}{*}[2.0cm]{
            \begin{subfigure}[c]{0.25\textwidth}
                \centering
                \includegraphics[width=\linewidth]
                {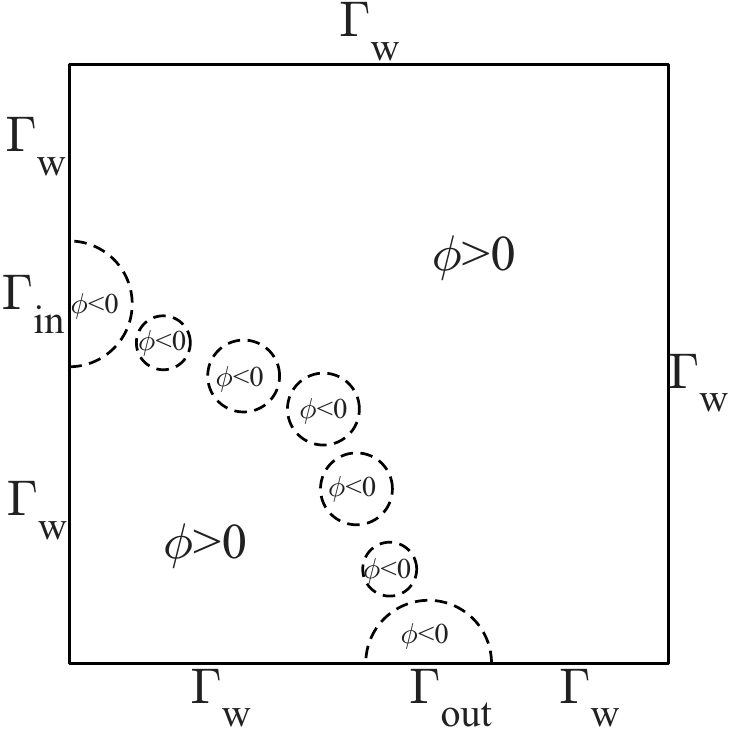}
                \caption{Domain and boundary conditions}
                \label{fig:case1}
            \end{subfigure}
        }
        &
        \begin{subfigure}[b]{0.215\textwidth}
            \centering
            \includegraphics[width=\linewidth]
            {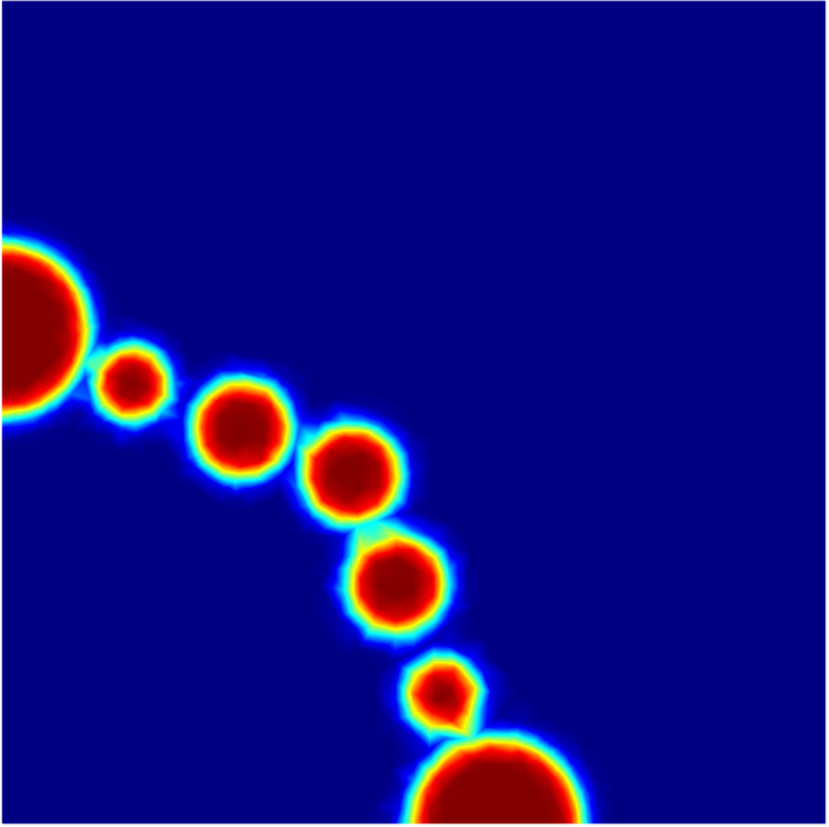}
            \caption{Initial Heaviside function:
            augmented objective $2275.34$;
            volume residual $-0.0953$}
            \label{fig:case2}
        \end{subfigure}
        &
        \begin{subfigure}[b]{0.215\textwidth}
            \centering
            \includegraphics[width=\linewidth]
            {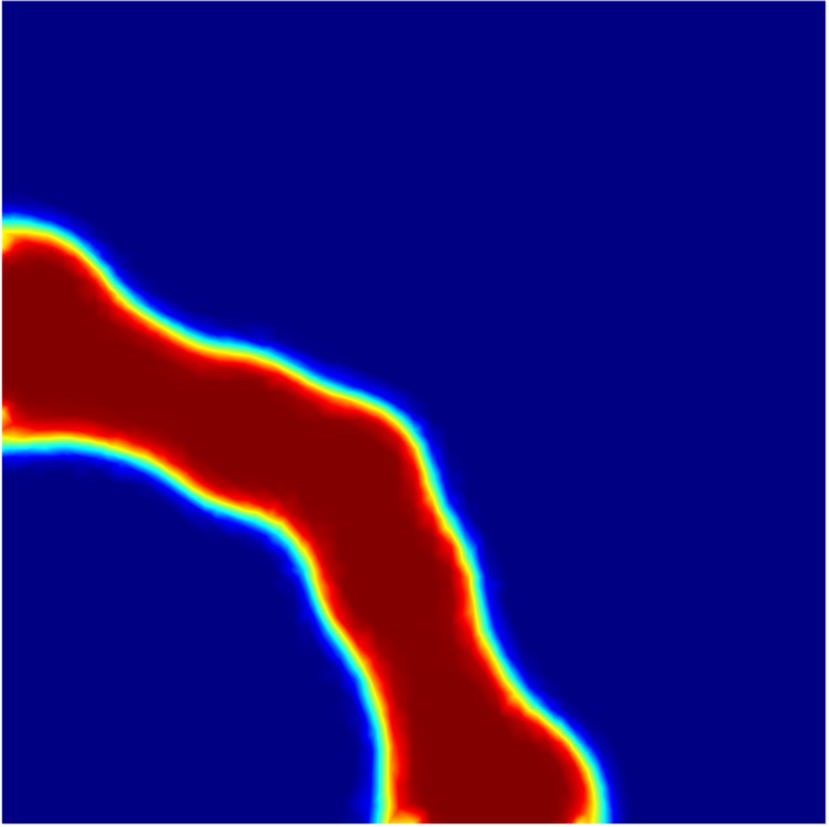}
            \caption{Heaviside function at iteration 50:
            augmented objective $1634.88$;
            volume residual $-4.90\times10^{-3}$}
            \label{fig:case3}
        \end{subfigure}
        &
        \begin{subfigure}[b]{0.215\textwidth}
            \centering
            \includegraphics[width=\linewidth]
            {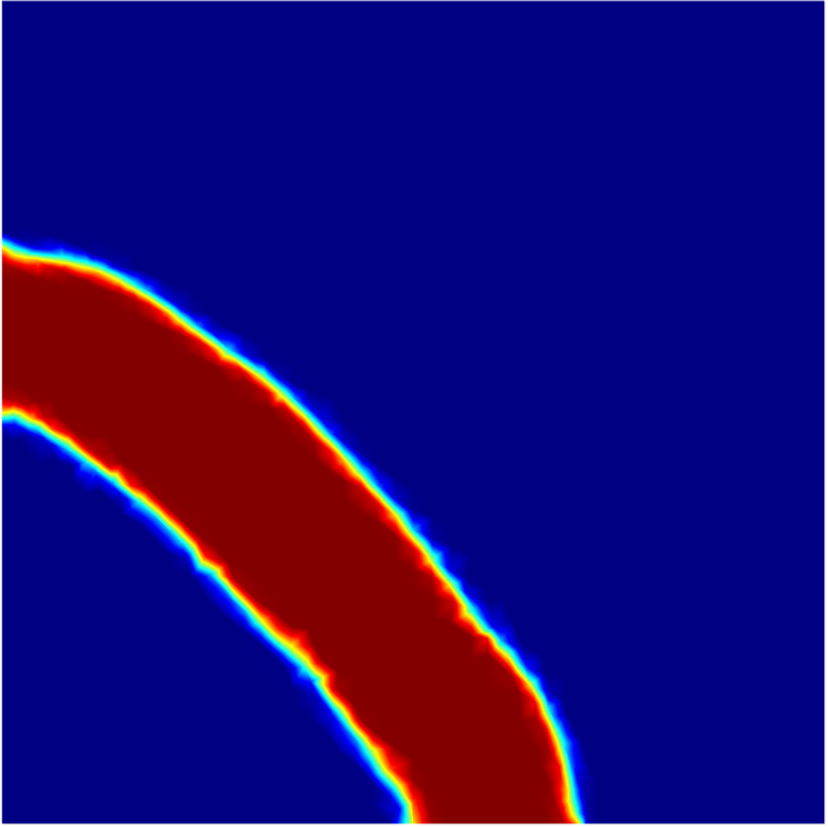}
            \caption{Optimized Heaviside function:
            augmented objective $1079.81$;
            volume residual $5.44\times10^{-3}$}
            \label{fig:case4}
        \end{subfigure}
        \\[0.1cm]

        &
        \begin{subfigure}[t]{0.215\textwidth}
            \centering
            \includegraphics[width=\linewidth]
            {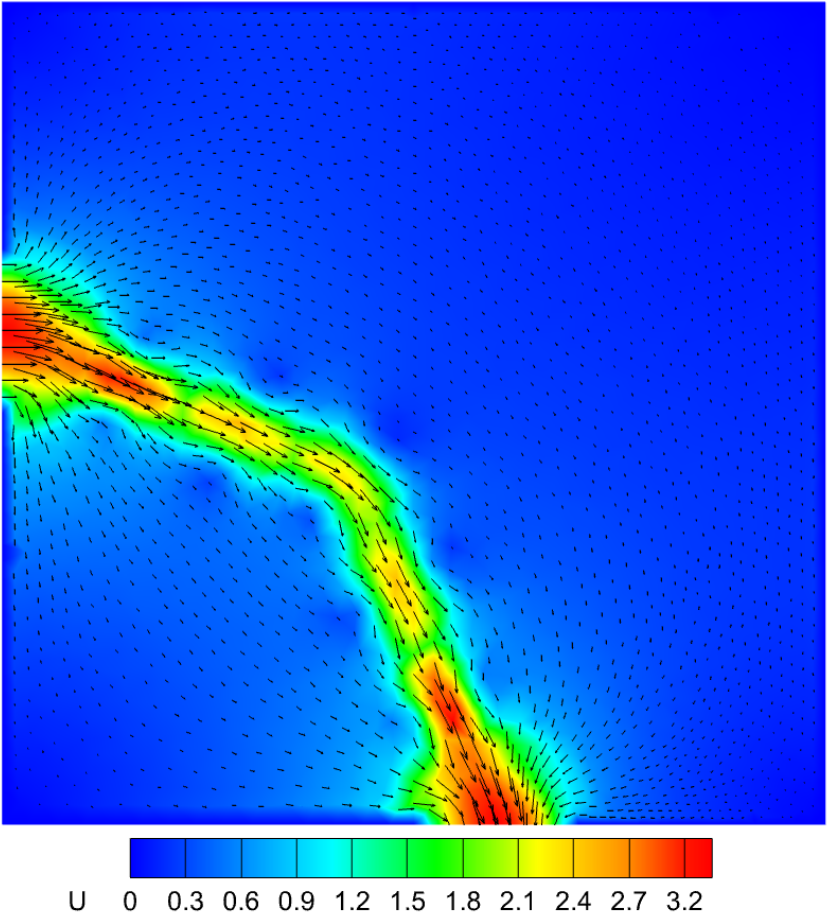}
            \caption{Velocity of the initial shape}
            \label{fig:case5}
        \end{subfigure}
        &
        \begin{subfigure}[t]{0.215\textwidth}
            \centering
            \includegraphics[width=\linewidth]
            {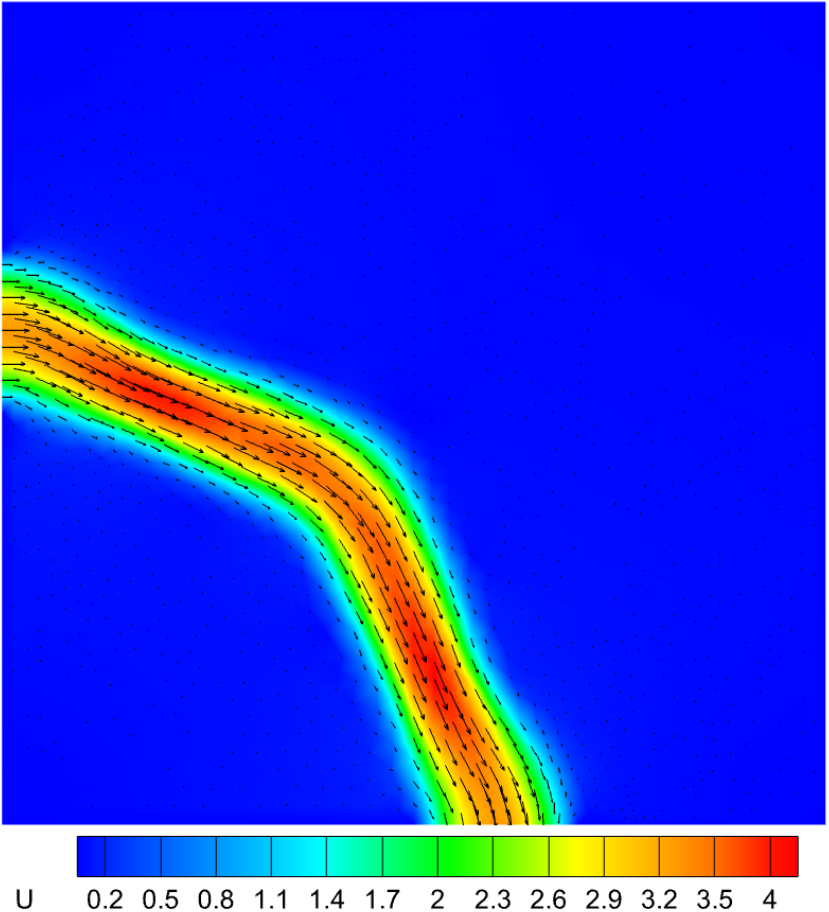}
            \caption{Velocity at iteration 50}
            \label{fig:case6}
        \end{subfigure}
        &
        \begin{subfigure}[t]{0.215\textwidth}
            \centering
            \includegraphics[width=\linewidth]
            {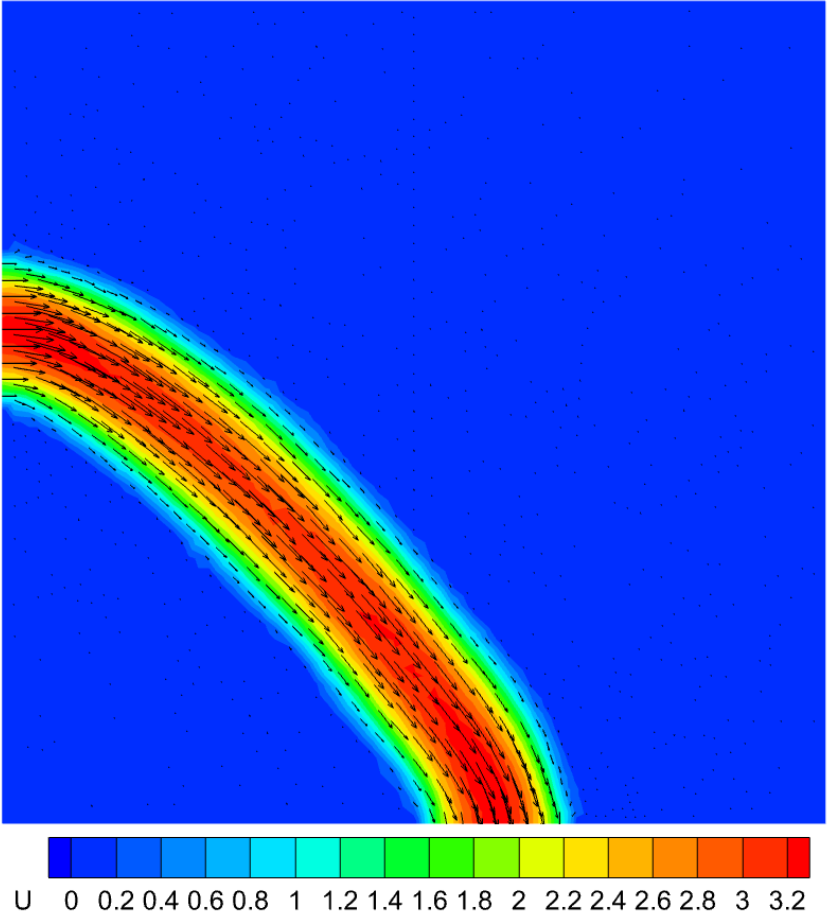}
            \caption{Velocity of optimized shape}
            \label{fig:case7}
        \end{subfigure}
    \end{tabular}

    \caption{Topology optimization of the L-shaped channel configuration:
    (a) computational domain and boundary conditions;
    (b)--(d) initial, iteration-50, and optimized
    Heaviside functions, respectively; and
    (e)--(g) corresponding velocity magnitudes.}
    \label{L_shape}
\end{figure*}

We consider the minimization of the viscous-energy objective for a coupled
Stokes--Darcy system in the unit square $\Omega=(0,1)^2$, subject to a
prescribed diffuse fluid volume.
As shown in Fig.~\ref{L_shape}(a), the inlet is located on the
left boundary and the outlet on the bottom boundary, defined respectively by
$\Gamma_{\rm in}=\{(0,y)\mid 0.5\leq y\leq0.7\}$
and
$\Gamma_{\rm out}=\{(x,0)\mid 0.5\leq x\leq0.7\}$.
The remaining parts of the boundary are denoted by $\Gamma_w$.
Non-homogeneous Dirichlet boundary conditions are prescribed as
$\left.\mathbf{u}\right|_{\Gamma_{\rm in}}
=\bigl(330(y-0.5)(0.7-y),0\bigr)^\top$
and
$\left.\mathbf{u}\right|_{\Gamma_{\rm out}}
=\bigl(0,-330(x-0.5)(0.7-x)\bigr)^\top$,
and the no-slip condition $\mathbf{u}=\mathbf{0}$ is imposed on $\Gamma_w$.
The physical parameters are $\mu=0.5$, $K=5\times10^{-5}$,
and $\alpha=1$.
We set $\eta=10^{-3}$, $\varepsilon=h_{\max}$, $\zeta=10^{-8}$,
$\Delta s=0.003$, $\kappa=0.002$, and $\ell_s=6h_{\max}$.
The prescribed diffuse fluid volume is $V_0=0.18|\Omega|$.
The augmented-Lagrangian parameters are $\lambda^0=0$ and
$\rho=5\times10^5$. The multiplier is updated every
$N_\lambda=5$ iterations according to
$
\lambda\leftarrow\lambda+0.1\rho C(\phi_h).
$

We use a two-stage continuation strategy in which the
viscous-energy term $\mathcal J_1$ is replaced by
$c_E\mathcal J_1$ in the augmented functional and the
corresponding design sensitivity.
During the first 30 updates, we set $c_E=0.1$ and $\tau=0.5$;
thereafter, we use $c_E=10$ and $\tau=5$.
All other algorithmic settings remain unchanged.
No periodic reinitialization is applied after the initial
signed-distance construction ($N_{\rm re}=0$).
We use $N_{\max}=450$ and $\mathrm{tol}_{\mathcal M}=10^{-4}$.
The interface displacement in one update is limited to $0.025\varepsilon$.
The unstructured mesh contains $5{,}298$ triangles and $2{,}750$ vertices.
The initial fluid region consists of seven disconnected circular seeds
arranged along a bent path between the inlet and outlet, as shown in
Fig.~\ref{L_shape}(b). Before optimization, the corresponding discrete
level-set function is converted to a signed-distance function by the
fast marching method.

Figures~\ref{L_shape}(b)--(d) show the initial, intermediate (iteration 50),
and optimized phase fields, and Figs.~\ref{L_shape}(e)--(g) show the
corresponding velocity magnitudes.
The seven initial fluid components merge into a single passage connecting
the inlet and outlet before iteration 20, and the passage remains connected
during the subsequent evolution. At iteration 50, the diffuse-volume
residual is $-4.90\times10^{-3}$. The calculation terminates at iteration
210 when the relative change of the fixed-penalty merit function falls below
$10^{-4}$. The final augmented-Lagrangian value is
$1.07981\times10^3$, and the final diffuse-volume residual is
$5.44\times10^{-3}$. The velocity field is concentrated mainly in the
resulting free-flow passage.

\subsubsection{Diffuser design}

\begin{figure}[!htbp]
 \begin{subfigure}[hb]{0.33\textwidth}
        \centering
        \includegraphics[width=1.8 in]
        {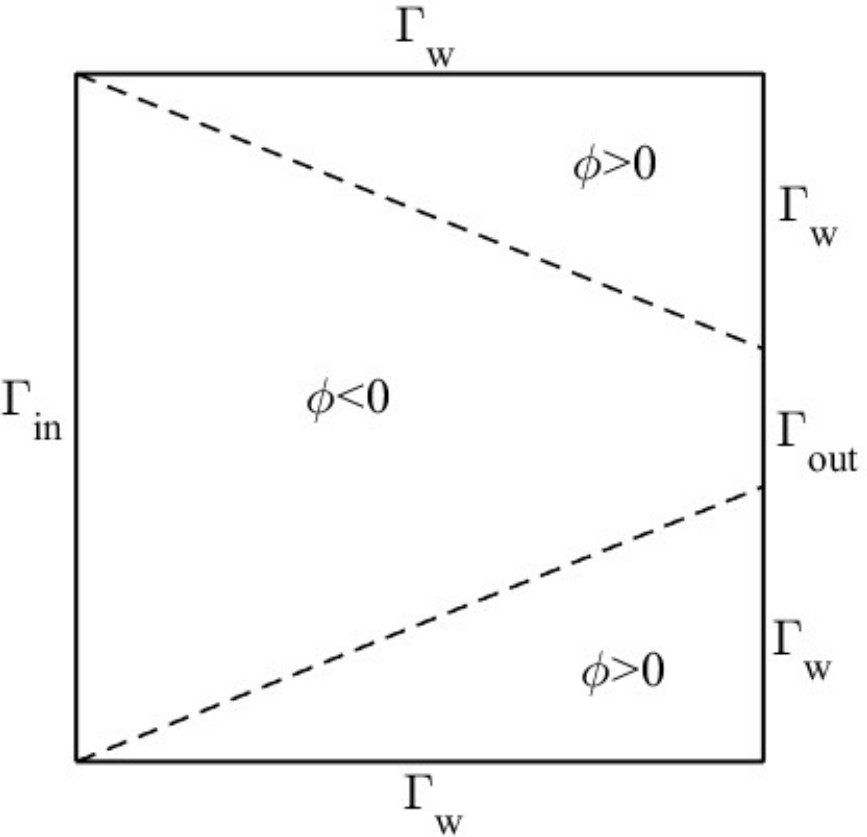}
    \caption{Domain and boundary conditions}
 \end{subfigure}
 \begin{subfigure}[hb]{0.33\textwidth}
        \centering
        \includegraphics[width=1.7 in]
        {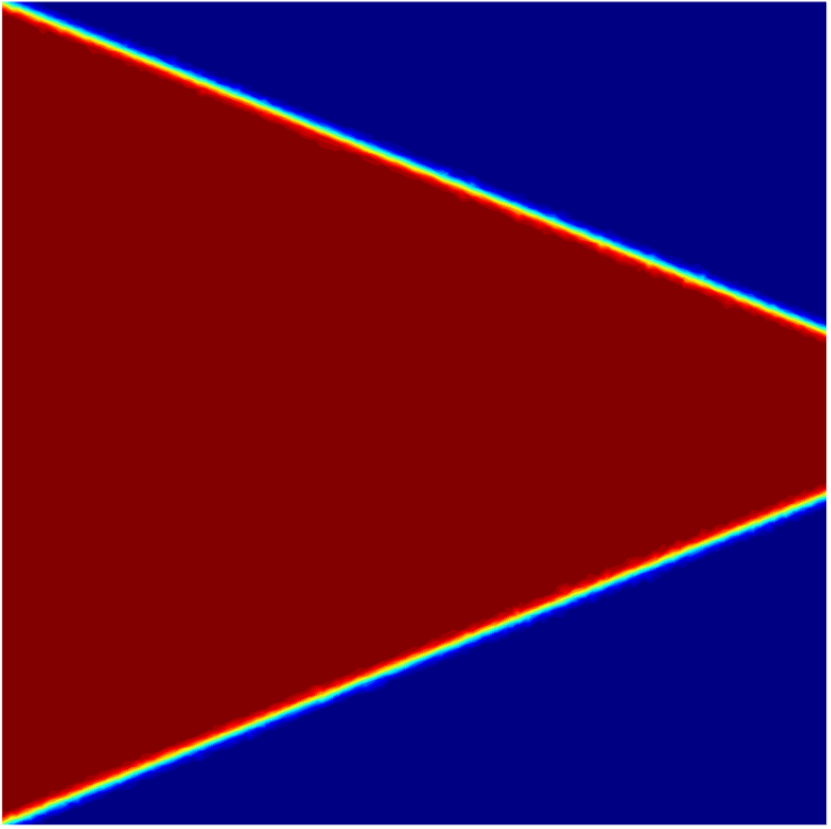}
    \caption{Initial Heaviside function}
 \end{subfigure}
 \begin{subfigure}[hb]{0.33\textwidth}
    \centering
    \includegraphics[width=1.9 in]
    {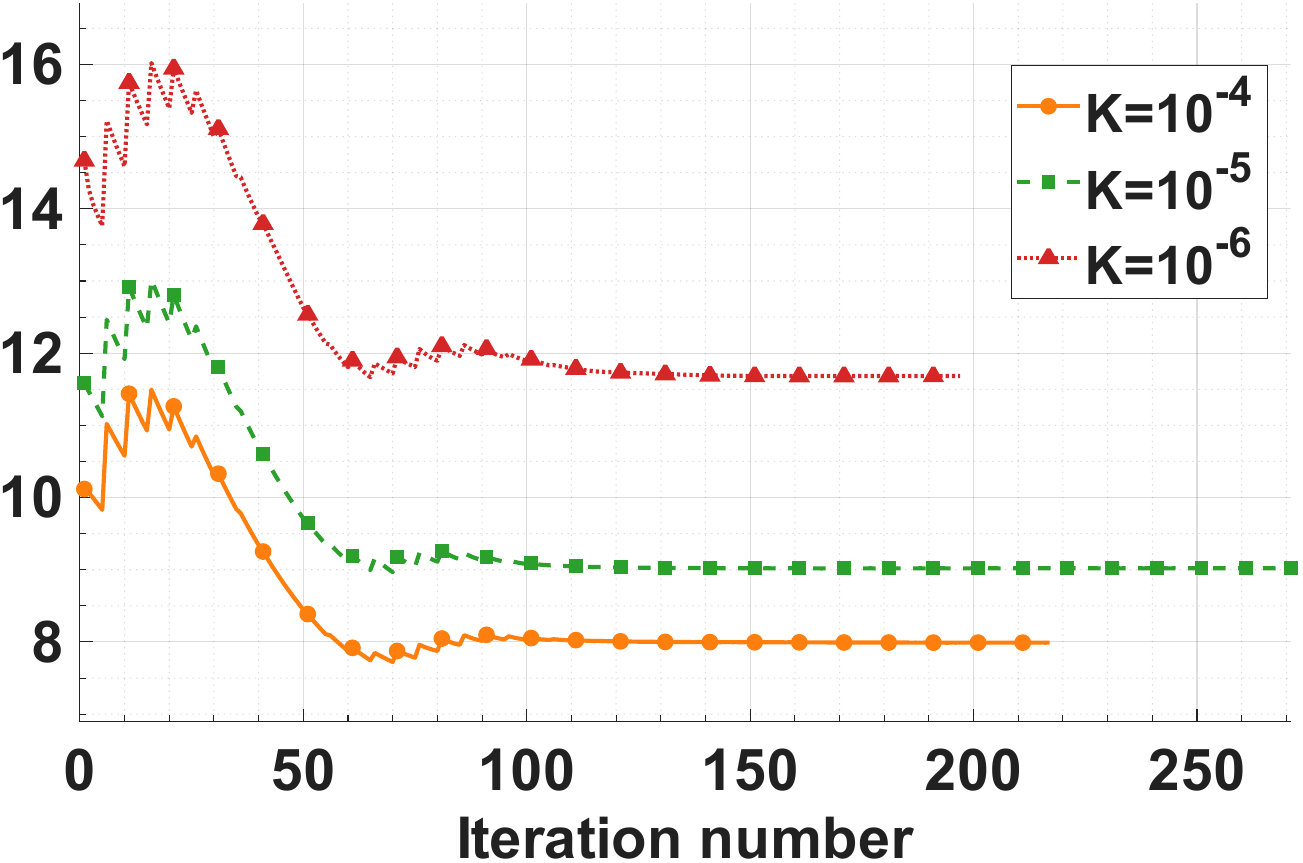}
    \caption{Augmented-Lagrangian histories}
 \end{subfigure}
\\
 \begin{subfigure}[h]{0.33\textwidth}
    \centering
    \includegraphics[width=1.7 in]
    {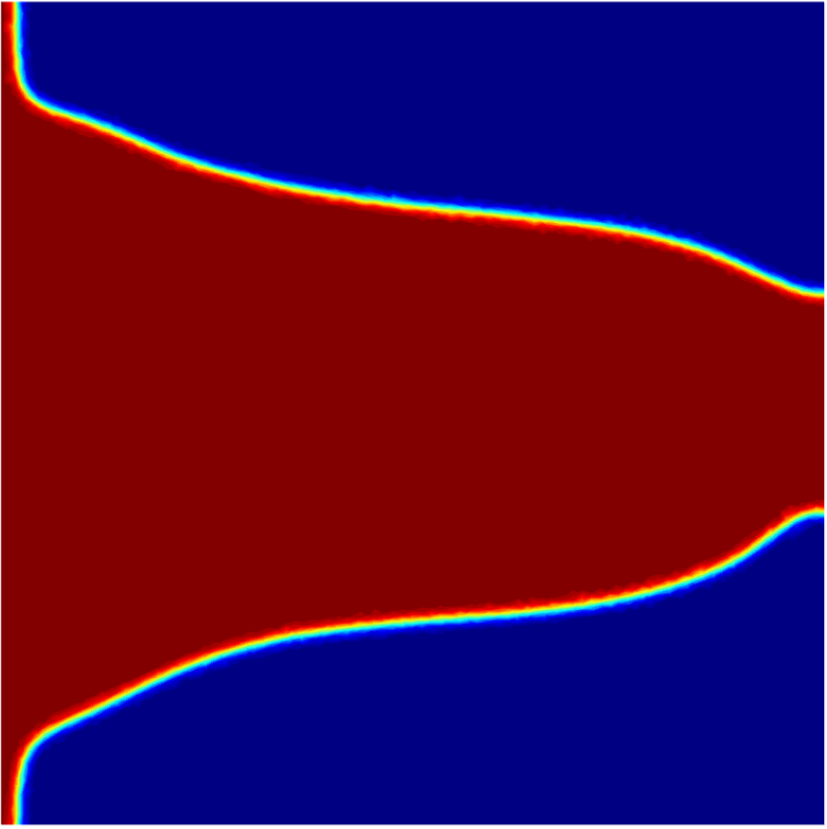}
    \caption{Optimized Heaviside function, $K=10^{-4}$:
    augmented objective $7.99054$;
    volume residual $-9.61\times10^{-5}$}
 \end{subfigure}
 \begin{subfigure}[h]{0.33\textwidth}
        \centering
        \includegraphics[width=1.7 in]
        {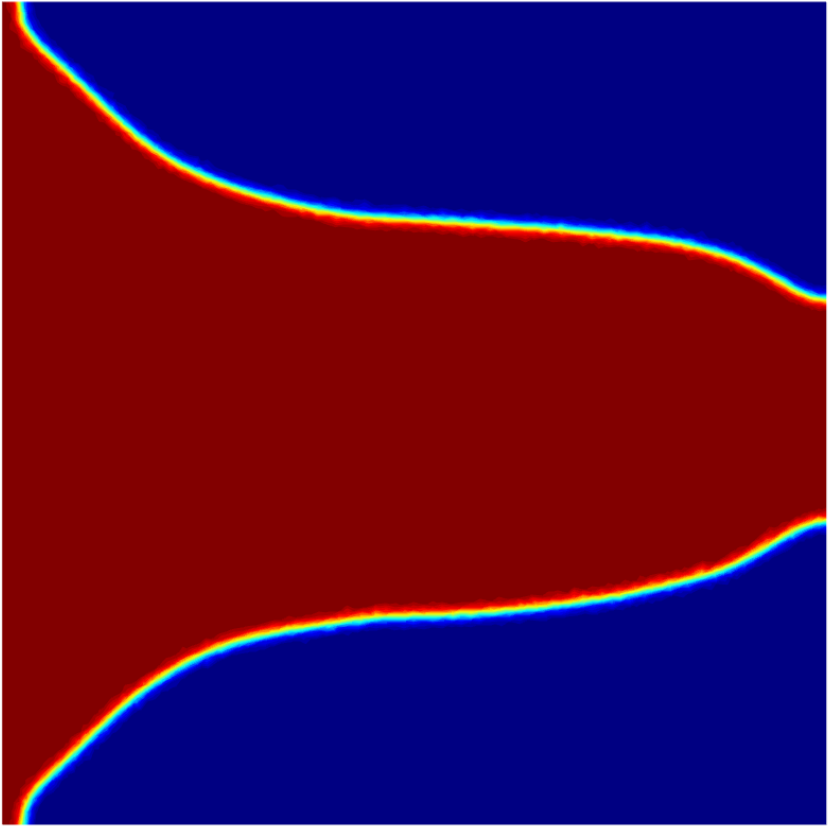}
    \caption{Optimized Heaviside function, $K=10^{-5}$:
    augmented objective $9.02048$;
    volume residual $4.04\times10^{-6}$}
 \end{subfigure}
 \begin{subfigure}[h]{0.33\textwidth}
    \centering
    \includegraphics[width=1.7 in]
    {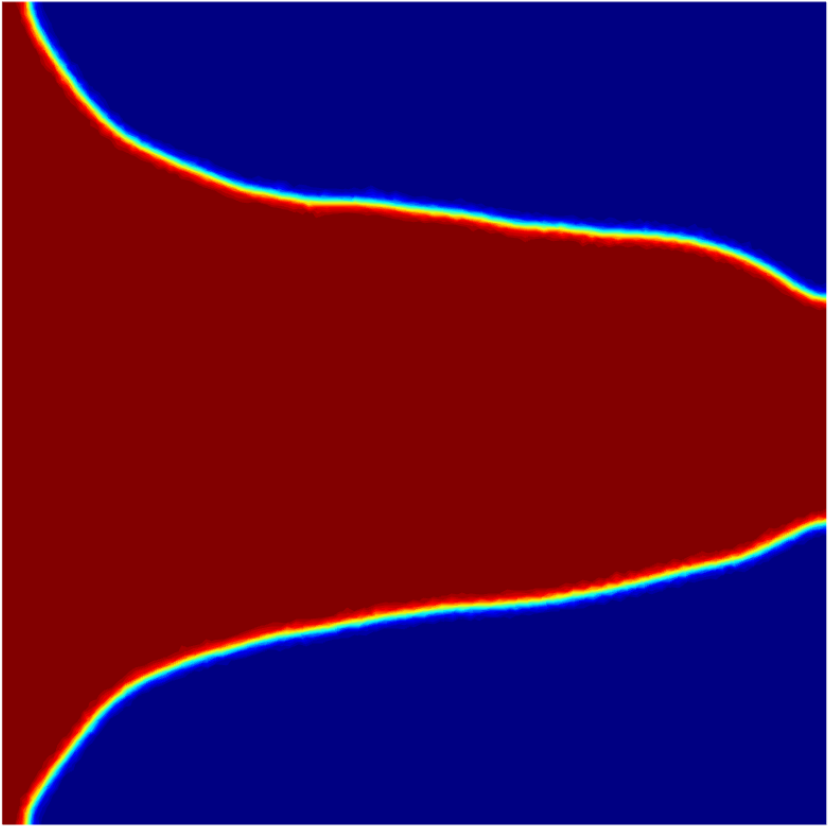}
    \caption{Optimized Heaviside function, $K=10^{-6}$:
    augmented objective $11.6843$;
    volume residual $2.05\times10^{-4}$}
 \end{subfigure}
\\
 \begin{subfigure}[h]{0.33\textwidth}
    \centering
    \includegraphics[width=1.7 in]
    {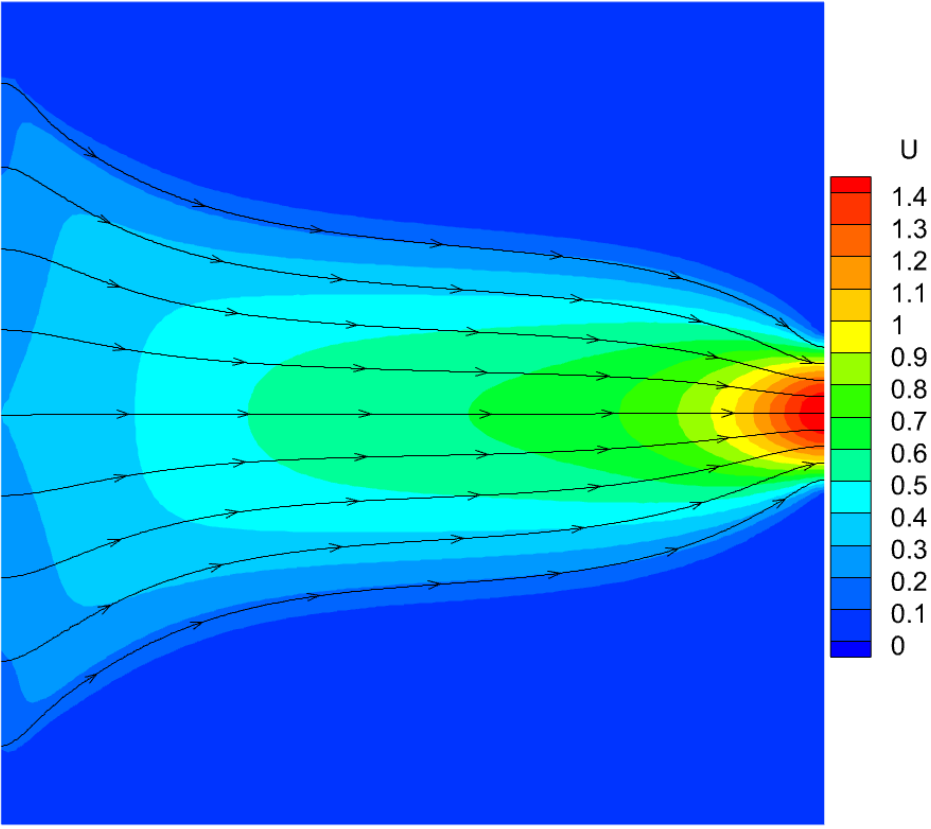}
    \caption{Velocity of optimized shape, $K=10^{-4}$}
 \end{subfigure}
 \begin{subfigure}[h]{0.33\textwidth}
        \centering
        \includegraphics[width=1.7 in]
        {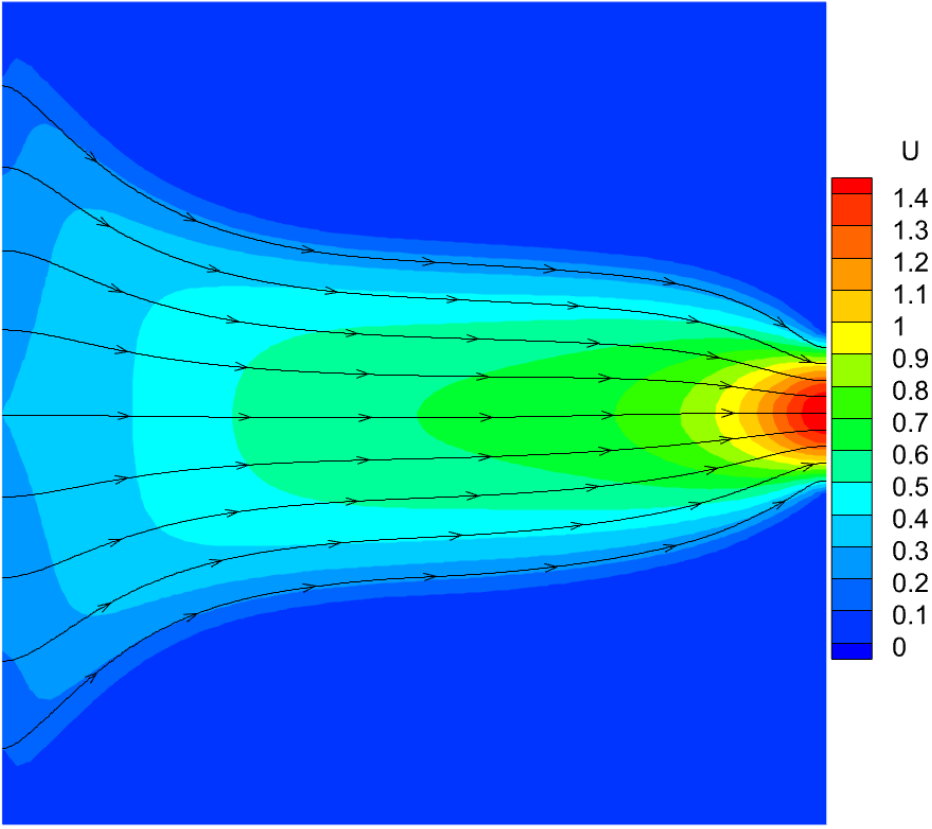}
    \caption{Velocity of optimized shape, $K=10^{-5}$}
 \end{subfigure}
 \begin{subfigure}[h]{0.33\textwidth}
    \centering
    \includegraphics[width=1.7 in]
    {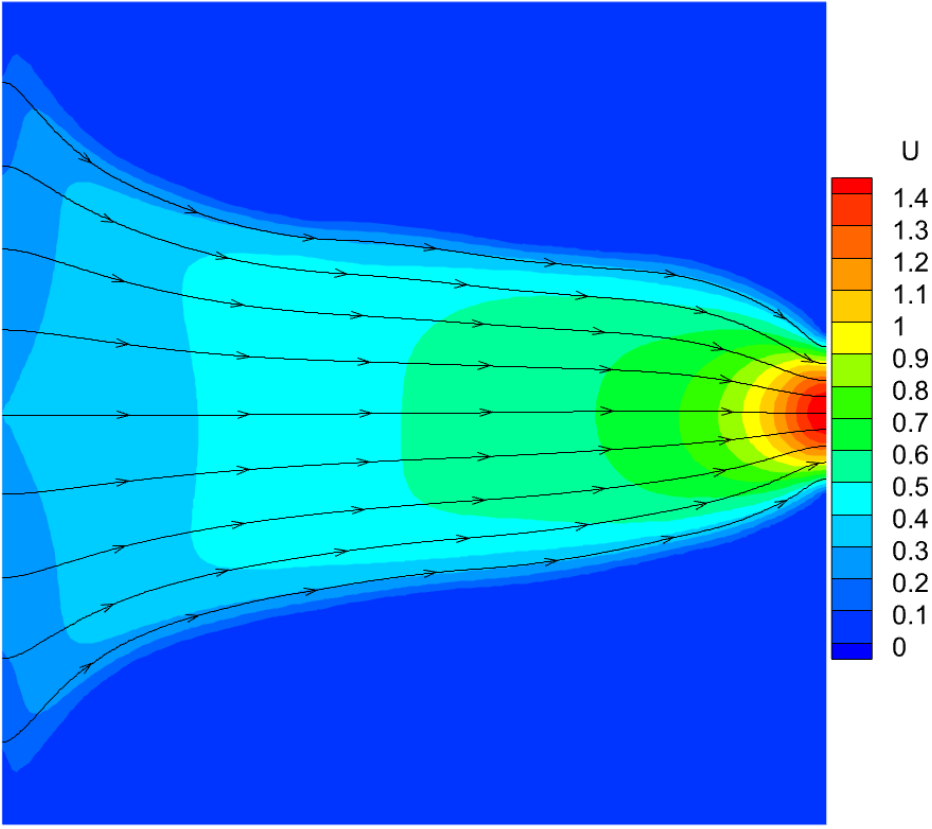}
    \caption{Velocity of optimized shape, $K=10^{-6}$}
 \end{subfigure}
 \caption{Topology optimization of the diffuser configuration for three
 porous-medium permeabilities: (a) computational domain and boundary
 conditions; (b) initial Heaviside function;
 (c) augmented-Lagrangian histories for
 $K=10^{-4}$, $10^{-5}$, and $10^{-6}$;
 (d)--(f) optimized Heaviside functions for
 $K=10^{-4}$, $10^{-5}$, and $10^{-6}$, respectively;
 (g)--(i) corresponding velocity magnitudes, respectively.}
 \label{Diffusion}
\end{figure}

We consider the minimization of the viscous-energy objective
\eqref{eq:objective1} for a Stokes--Darcy diffuser with porous walls in the
unit square $\Omega=(0,1)^2$. As shown in Fig.~\ref{Diffusion}(a), the inlet
and outlet are
$\Gamma_{\mathrm{in}}=\{(0,y)\mid 0\leq y\leq1\}$ and
$\Gamma_{\mathrm{out}}=\{(1,y)\mid 0.4\leq y\leq0.6\}$, respectively, and
the remaining boundary is denoted by $\Gamma_w$. The prescribed velocities are
$
\left.\mathbf{u}\right|_{\Gamma_{\mathrm{in}}}
=\bigl(1.2y(1-y),0\bigr)^{\top}
$
and
$
\left.\mathbf{u}\right|_{\Gamma_{\mathrm{out}}}
=\bigl(150(y-0.4)(0.6-y),0\bigr)^{\top},
$
and the no-slip condition $\mathbf{u}=\mathbf{0}$ is imposed on
$\Gamma_w$. The initial level-set function is
$\phi_0(x,y)=-\min\{y-0.4x,\,1-0.4x-y\}$, as displayed in
Fig.~\ref{Diffusion}(b).

The unstructured mesh contains $15{,}082$ triangles and
$7{,}702$ vertices.
We use
$\mu=1$, $\alpha=1$,
$\eta=10^{-3}$, $\zeta=10^{-8}$, $\varepsilon=h_{\max}$,
$\Delta s=0.005$, $\kappa=0.01$, $\tau=1$, and
$\ell_s=4h_{\max}$. The prescribed diffuse fluid volume is
$V_0=0.88\int_\Omega\mathcal H_f(\phi_h^0)\,\mathrm{d}x$.
The augmented-Lagrangian parameters are $\lambda^0=0$, $\rho=900$,
$\sigma=300$, and $N_\lambda=5$, with the multiplier
updated according to
$\lambda\leftarrow\lambda+\sigma C(\phi_h)$ every five iterations.
The explicit reaction displacement in one update is limited to
$0.08\varepsilon$, and fast-marching reinitialization is applied after every
update ($N_{\rm re}=1$). The fluid phase is prescribed only on the physical
inlet and outlet boundaries by setting $\phi_h=-2\varepsilon$.
We set $N_{\max}=300$ and use the stopping tolerance
$\mathrm{tol}_{\mathcal M}=10^{-4}$.
To examine the influence of permeability, we consider
$K=10^{-4}$, $10^{-5}$, and $10^{-6}$ while keeping all other parameters
fixed.

Figure~\ref{Diffusion}(c) shows the augmented-Lagrangian histories.
The runs for $K=10^{-4}$, $10^{-5}$, and $10^{-6}$ terminate after
217, 300, and 197 iterations, respectively; the $K=10^{-5}$ run reaches
$N_{\max}$, while the other two satisfy the prescribed merit-change
criterion. Figures~\ref{Diffusion}(d)--(f) present the corresponding final diffuse
fluid indicators.
In each case, the free-flow region forms a single connected
passage between the inlet and outlet.
The corresponding velocity magnitudes are shown in
Figs.~\ref{Diffusion}(g)--(i) using the same color scale. The velocity is
concentrated within the free-flow passage and reaches its largest values near
the outlet, while it remains small in the porous region.

\subsubsection{Three-dimensional channel design with different permeabilities}
\label{sec:3d-permeability-design}

We next investigate the influence of permeability on three-dimensional
channel design in the unit cube $\Omega=(0,1)^3$. Two circular inlets on
$x=0$ are centered at $(0,0.25,0.35)$ and $(0,0.75,0.65)$, and two
outlets on $x=1$ have the same respective $y$- and $z$-coordinates.
All four ports have radius $0.12$. Axial parabolic velocity profiles
are prescribed at the ports and normalized so that each port carries
the same discrete flux, approximately $0.0223114$. The velocity vanishes
on the remaining boundary. The four port centers lie in the plane
$z=0.6y+0.2$, although the state equations are solved throughout the
three-dimensional domain.

The objective is the viscous dissipation $\mathcal J_1$ defined in
\eqref{eq:objective1}; Darcy and interface dissipation enter the state
equations but are not added to the objective. We compare
$K=10^{-7},10^{-6},10^{-5}$ with $\mu=1$, $\alpha=1$,
$\eta=10^{-3}$, $\varepsilon=0.05$, and $\zeta=10^{-4}$.
The same connected four-port initial design and the same mesh of
$162{,}000$ tetrahedra and $29{,}791$ vertices are used in all cases.
The velocity--pressure pair is discretized with the MINI element,
and the level-set function with continuous piecewise linear elements.
The prescribed diffuse fluid volume is
$V_0=0.10$, where $V=\int_\Omega\mathcal H_f(\phi_h)\dx$.
Each computation uses 30 optimization iterations.

For numerical scaling, the implementation divides the objective by its
initial value and applies the augmented Lagrangian to the relative
volume residual $V/V_0-1$, with penalty parameter $200$ and initial
multiplier zero.  
With $c_h:=V/V_0-1$, trial updates are accepted using the
augmented energy
\[
\frac{\mathcal J_{1,h}}{\mathcal J_{1,h}^{0}}
+\lambda c_h+\frac{\rho}{2}c_h^2
+\frac{\tau}{2}\int_\Omega|\nabla\phi_h|^2\,\mathrm dx,
\]
where $\mathcal J_{1,h}^{0}$ is the initial objective value
for the corresponding permeability, and $\lambda$ is held
fixed during each line search.
 The multiplier is updated after each accepted step as
$\lambda\leftarrow\lambda+0.25\rho c_h$. The diffusion weight and
sensitivity smoothing length are both $0.30$. These numerical settings are fixed across the three runs.
Consequently, this experiment compares the outputs of the same
normalized optimization procedure, rather than holding the relative
weight of the unnormalized viscous objective and diffusion term fixed.

Figure~\ref{fig:3d-permeability-histories} shows
the unnormalized viscous objective
$\mathcal J_{1,h}$ and the diffuse fluid volume $V$
from iterations 0 to 30.
The viscous objective decreases at every
recorded step, from $109.2986$ to $38.9397$, from $51.7098$ to $18.5932$,
and from $21.3056$ to $12.0065$, respectively. These changes correspond
to reductions of $64.37\%$, $64.04\%$, and $43.65\%$.
The volume remains close to the prescribed value: the maximum relative
deviation over all three histories is less than $0.25\%$.
At iteration 30, the relative volume errors are approximately
$-0.06772\%$, $0.00209\%$, and $0.00071\%$, respectively.
The histories show objective reduction over the prescribed 30
optimization iterations.

Figure~\ref{fig:3d-permeability-fields} compares the resulting interfaces
and velocity fields. All three computed fluid regions remain connected.
The overall four-port topology is retained, while the central junction
and branch profiles vary with permeability. In the displayed view,
the $K=10^{-5}$ design has a broader central junction than the
$K=10^{-7}$ design. The common velocity scale also shows a more pronounced
central high-speed region for $K=10^{-7}$. Thus, permeability affects
both the local channel geometry and the velocity distribution without
changing the observed connectivity in these three runs. Since changing
$K$ modifies both the Darcy coefficient $\mu/K$ and the interface
coefficient $\alpha\mu/\sqrt K$, the differences reflect both effects;
the smaller objective values at larger $K$ should not be interpreted
as a ranking of geometries under an identical physical model.

\begingroup
\begin{figure}[!htbp]
 \centering
 \begin{subfigure}[t]{0.32\textwidth}
  \centering
  \includegraphics[width=\linewidth]{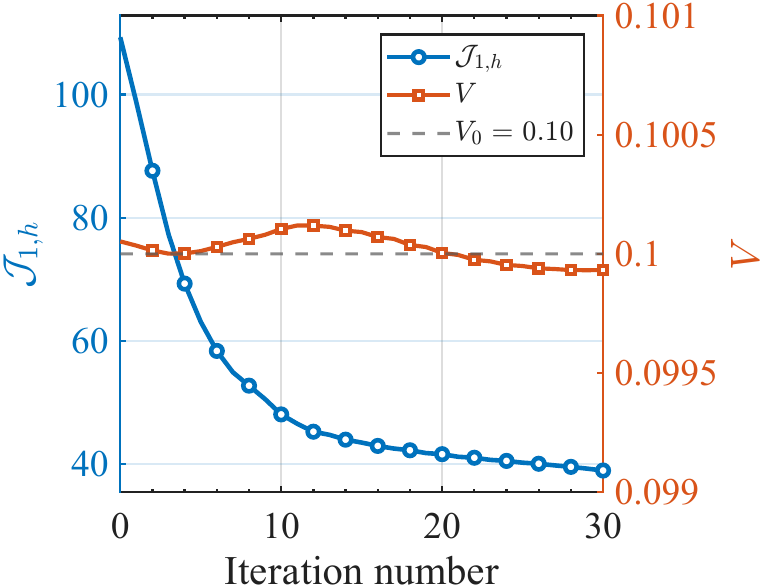}
  \caption{$K=10^{-7}$}
 \end{subfigure}\hfill
 \begin{subfigure}[t]{0.32\textwidth}
  \centering
  \includegraphics[width=\linewidth]{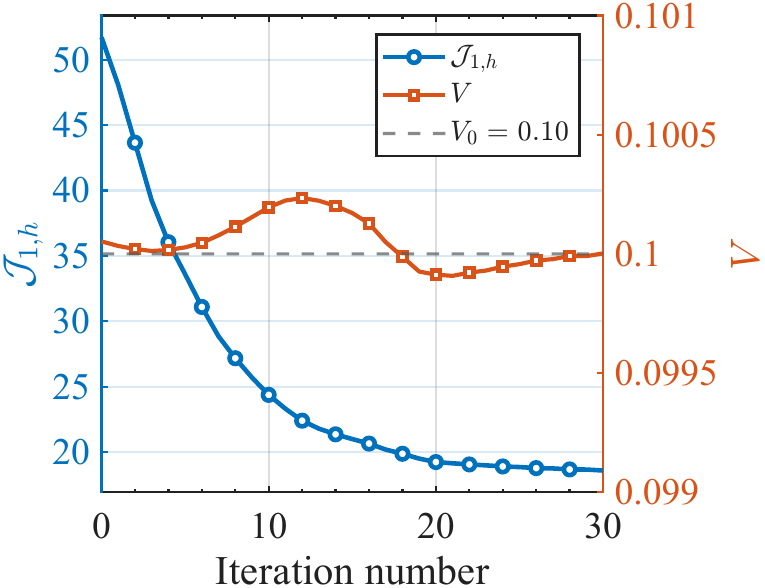}
  \caption{$K=10^{-6}$}
 \end{subfigure}\hfill
 \begin{subfigure}[t]{0.32\textwidth}
  \centering
  \includegraphics[width=\linewidth]{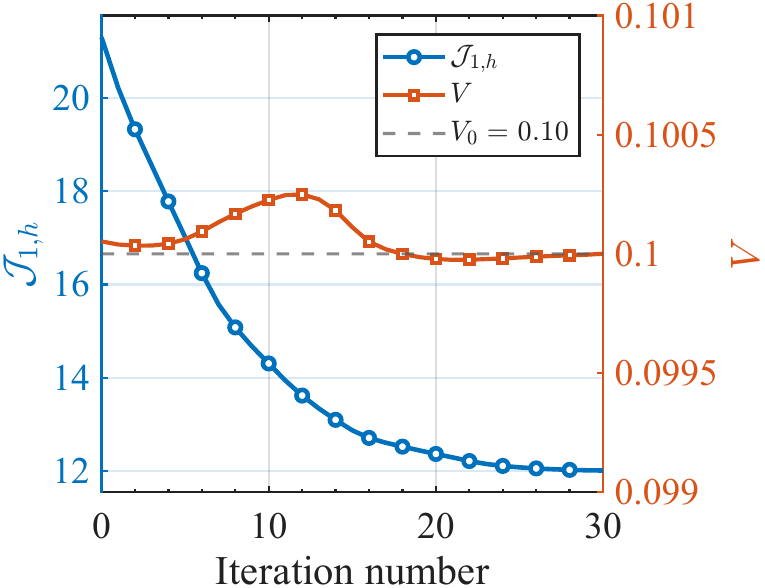}
  \caption{$K=10^{-5}$}
 \end{subfigure}
 \caption{Three-dimensional optimization histories. 
 }
 \label{fig:3d-permeability-histories}
\end{figure}

\begin{figure}[!htbp]
 \centering
 \begin{subfigure}[t]{0.32\textwidth}
  \includegraphics[width=\linewidth]{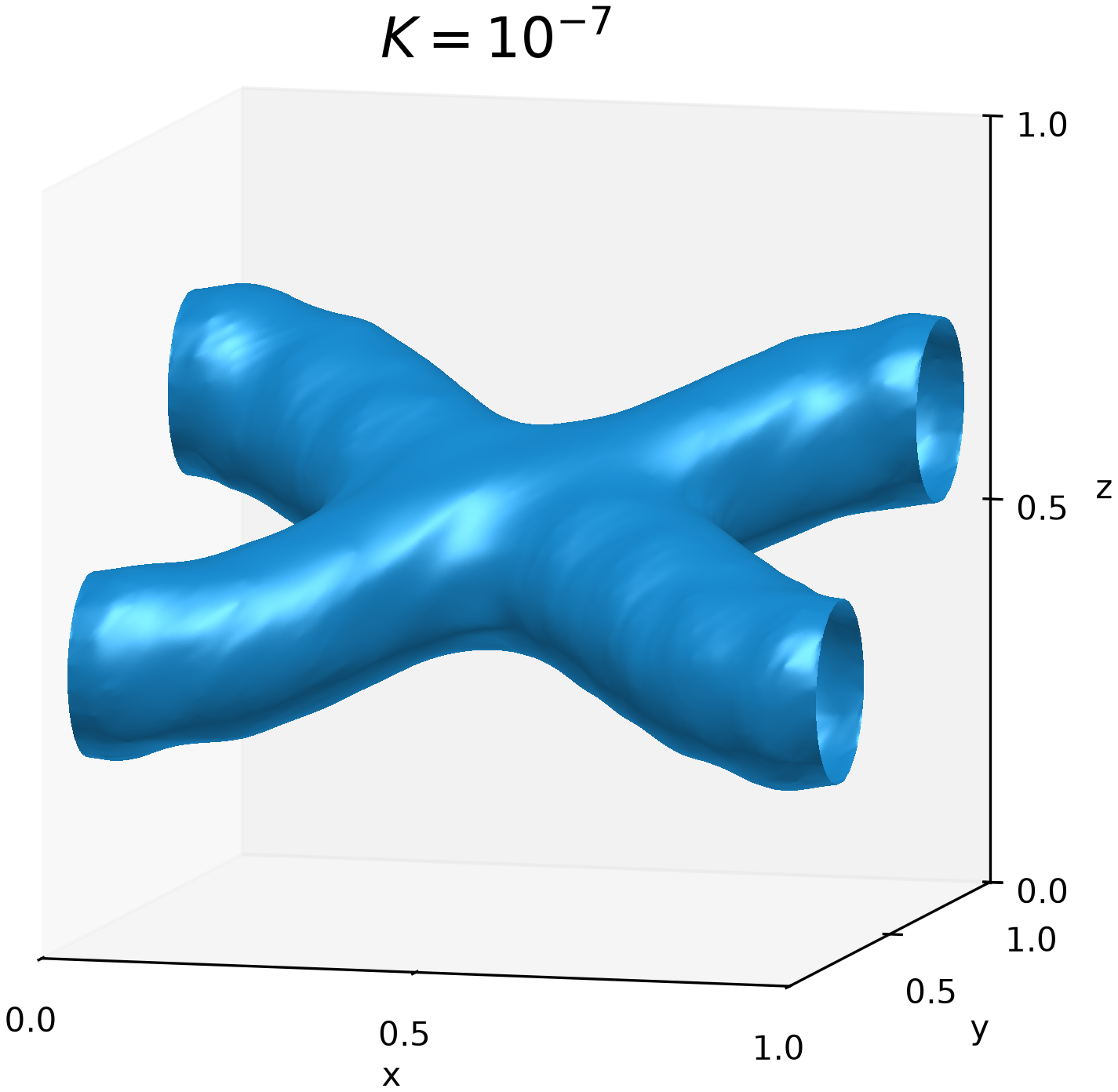}
  \caption{$K=10^{-7}$: interface.}
 \end{subfigure}\hfill
 \begin{subfigure}[t]{0.32\textwidth}
  \includegraphics[width=\linewidth]{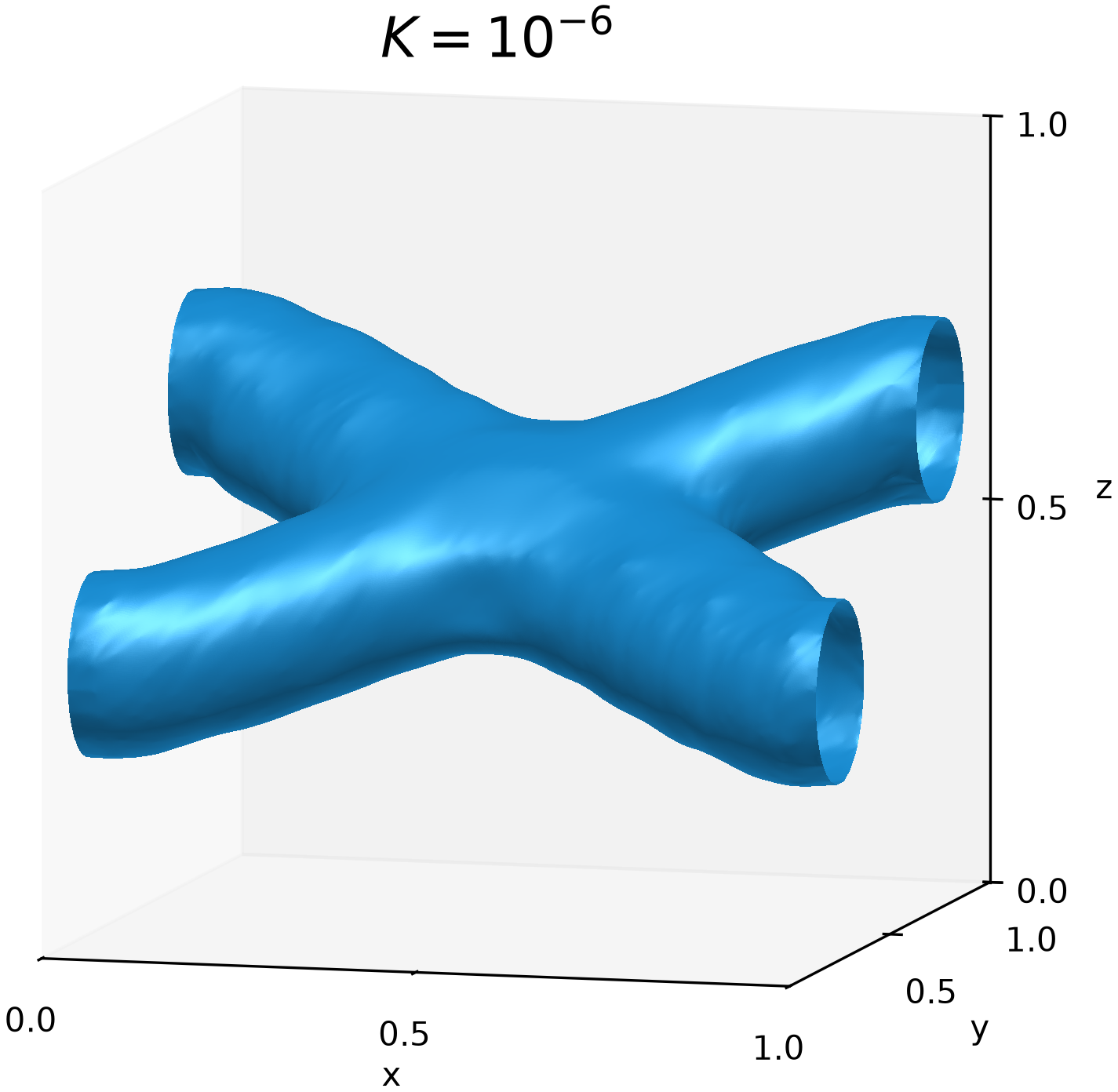}
  \caption{$K=10^{-6}$: interface.}
 \end{subfigure}\hfill
 \begin{subfigure}[t]{0.32\textwidth}
  \includegraphics[width=\linewidth]{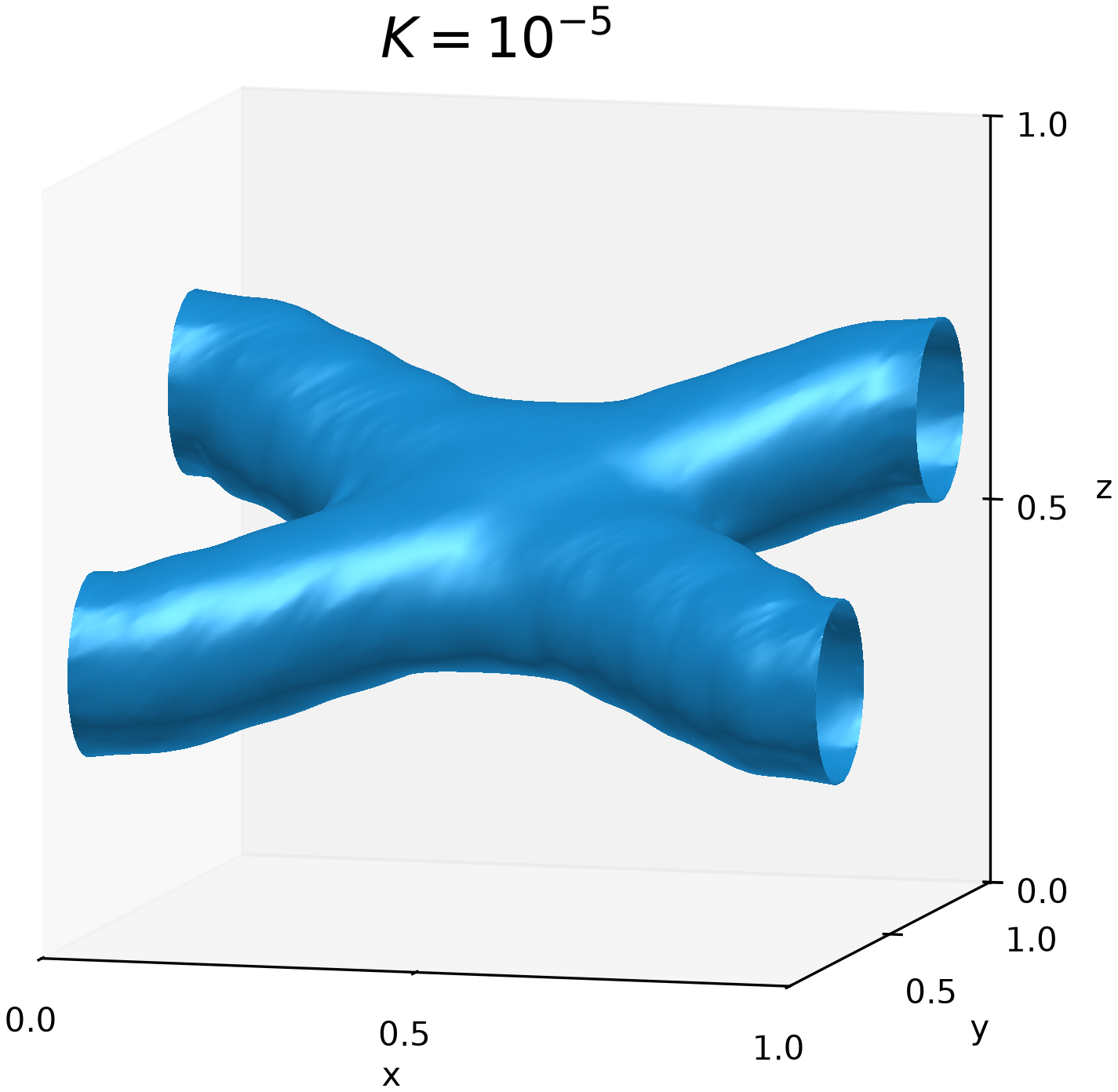}
  \caption{$K=10^{-5}$: interface.}
 \end{subfigure}
 \par\medskip
 \begin{subfigure}[t]{0.32\textwidth}
  \includegraphics[width=\linewidth]{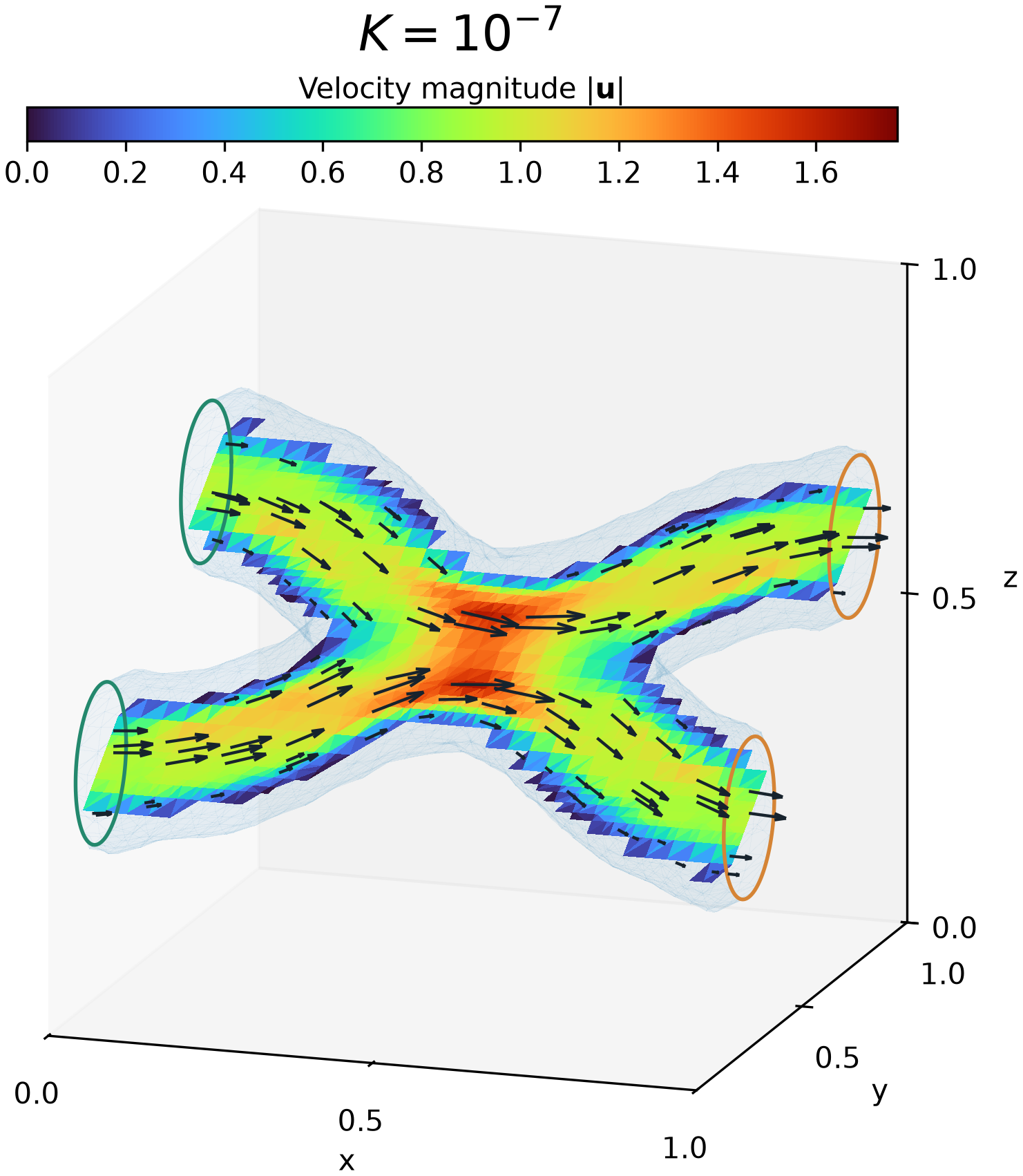}
  \caption{$K=10^{-7}$: velocity.}
 \end{subfigure}\hfill
 \begin{subfigure}[t]{0.32\textwidth}
  \includegraphics[width=\linewidth]{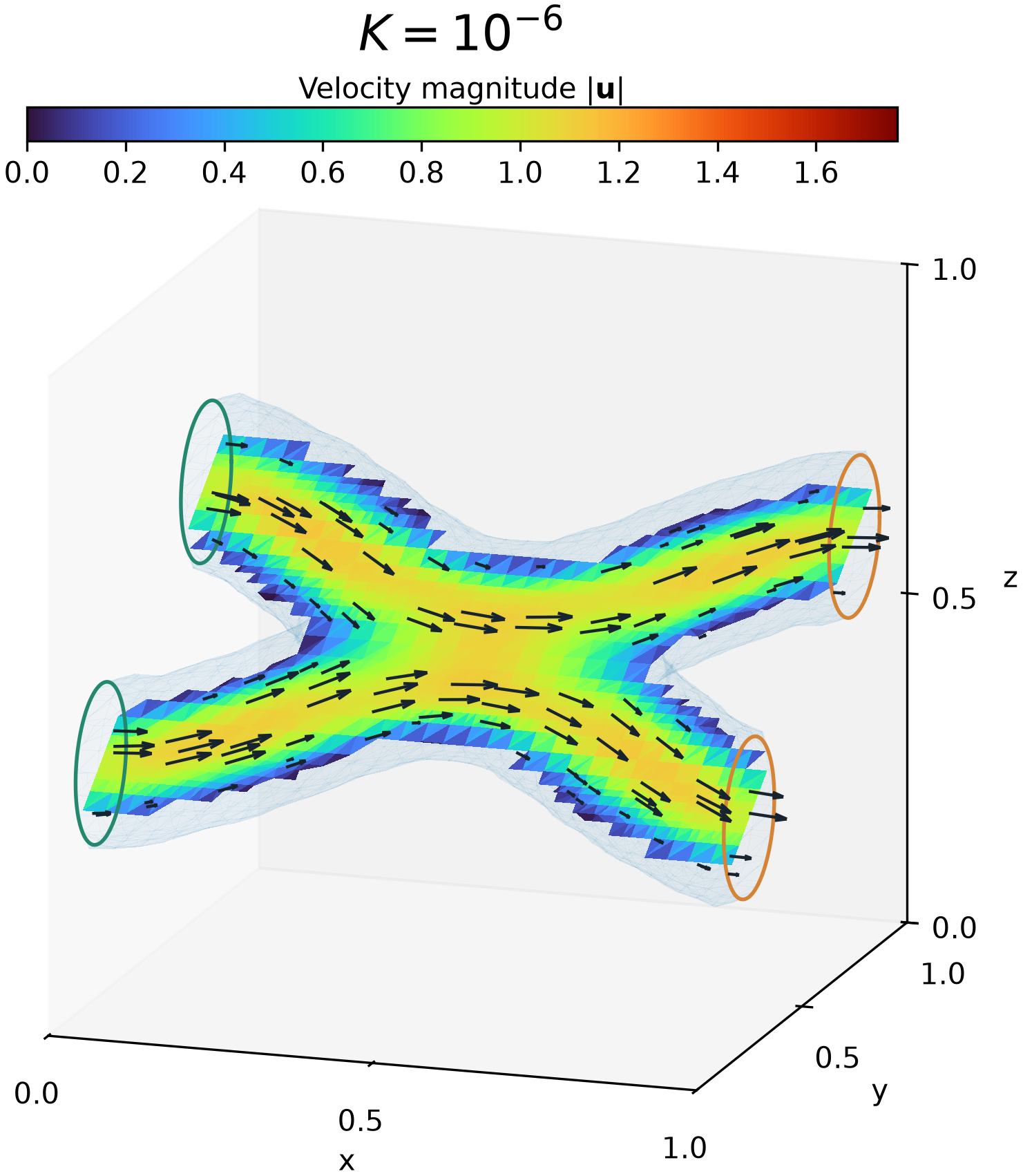}
  \caption{$K=10^{-6}$: velocity.}
 \end{subfigure}\hfill
 \begin{subfigure}[t]{0.32\textwidth}
  \includegraphics[width=\linewidth]{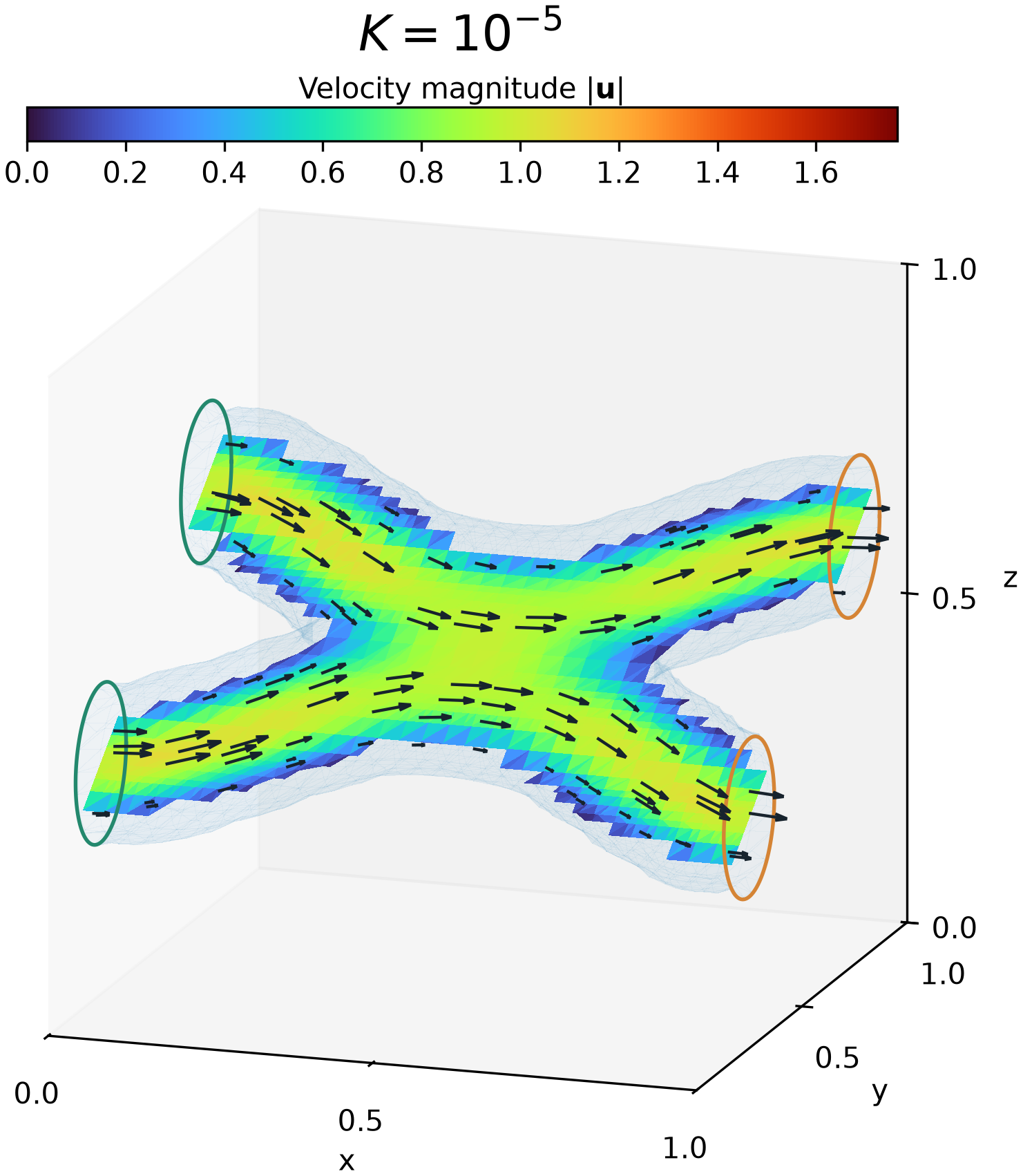}
  \caption{$K=10^{-5}$: velocity.}
 \end{subfigure}
 \caption{Three-dimensional designs and velocity fields: (a)–(c) optimized surfaces; (d)–(f) sampled three-dimensional velocity vectors.}
 \label{fig:3d-permeability-fields}
\end{figure}
\endgroup

\subsection{Inverse design via observation fitting}
\subsubsection{Complex riverbed sediment profile inversion based on flow-velocity tracking}
\begin{figure}[!htbp]
	\begin{subfigure}[b]{0.33\textwidth}
		\centering
		\includegraphics[width=1.8 in]{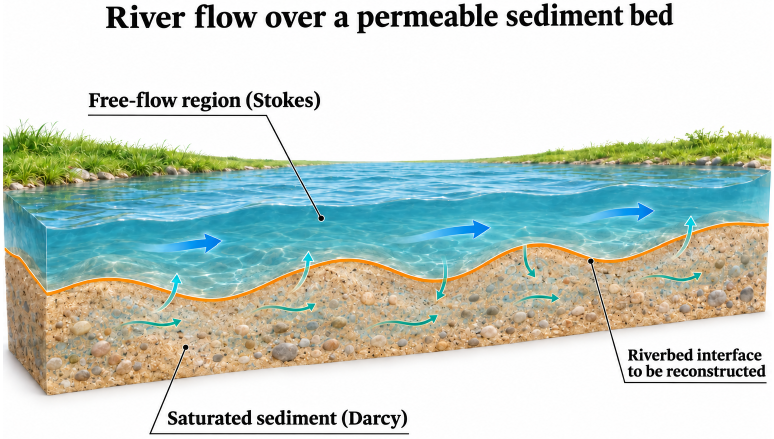}
		\caption{{\r Physical motivation}}
	\end{subfigure}\hfill
	\begin{subfigure}[b]{0.33\textwidth}
		\centering
		\includegraphics[width=1.8 in]{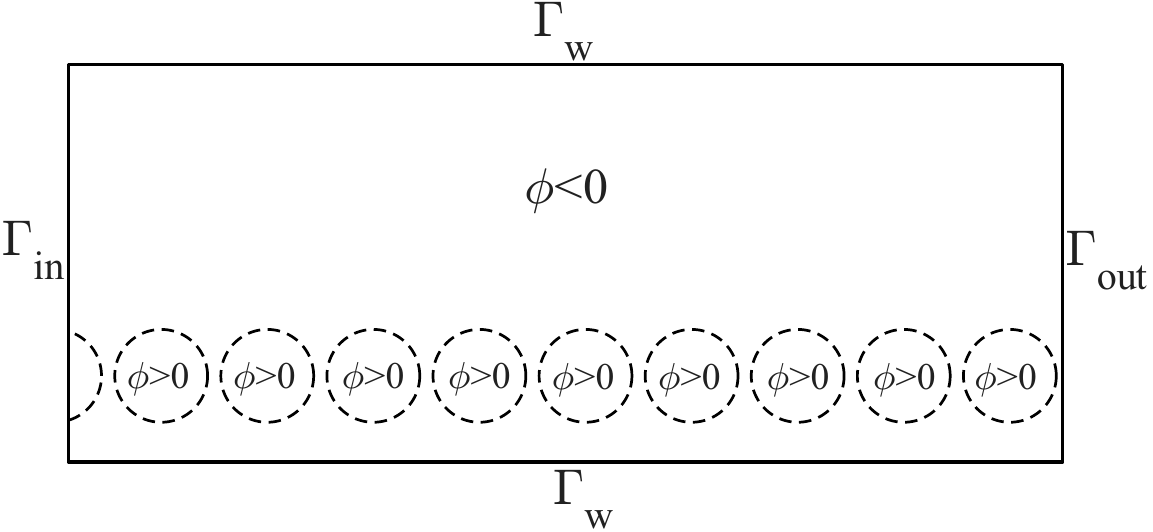}
		\caption{Domain and boundary conditions}
	\end{subfigure}\hfill
	\begin{subfigure}[b]{0.33\textwidth}
		\centering
		\includegraphics[width=1.8 in]{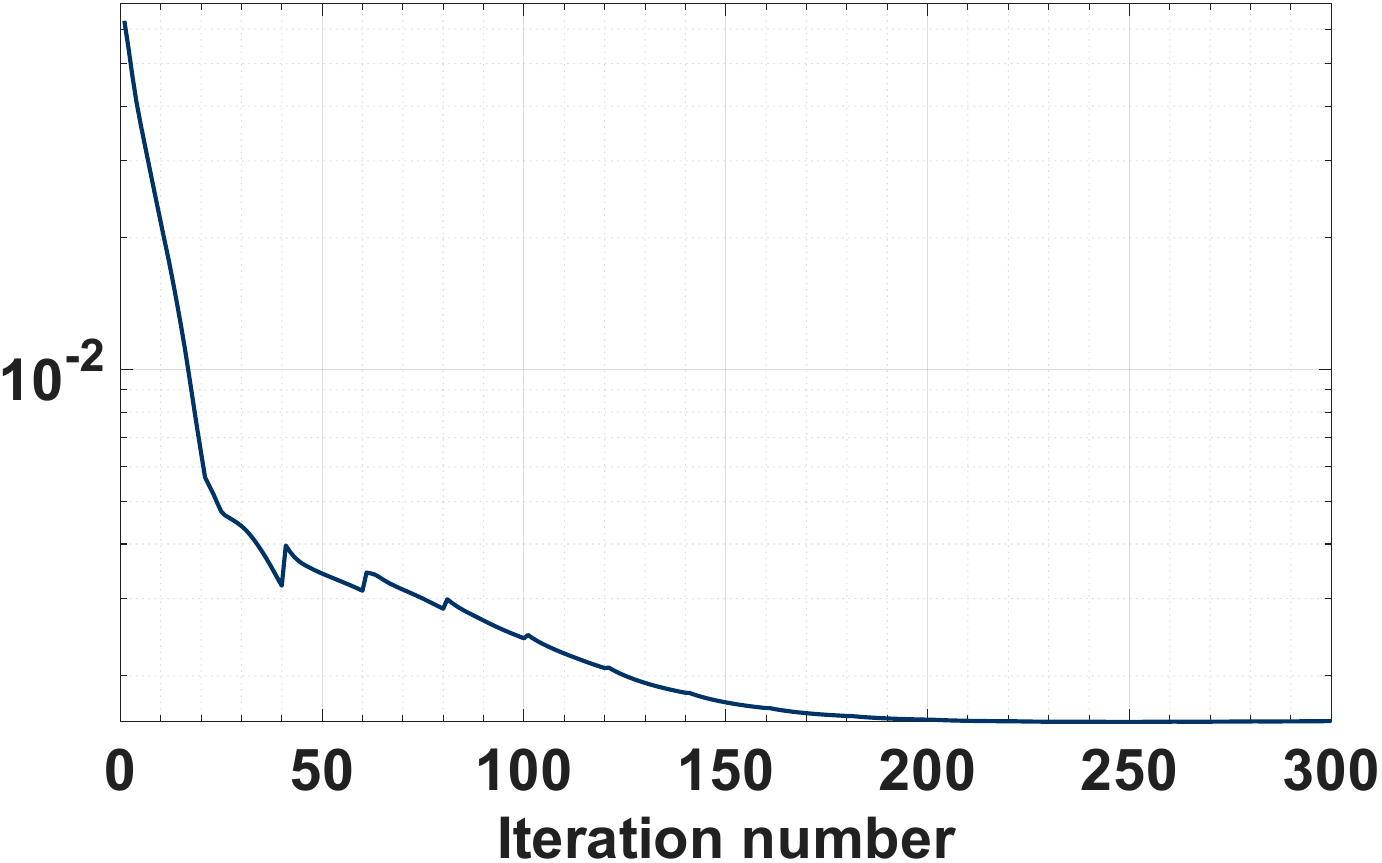}
		{\r\caption{Regularized-objective history}}
	\end{subfigure}\par\smallskip
	\begin{subfigure}[b]{0.33\textwidth}
		\centering
		\includegraphics[width=1.8 in]{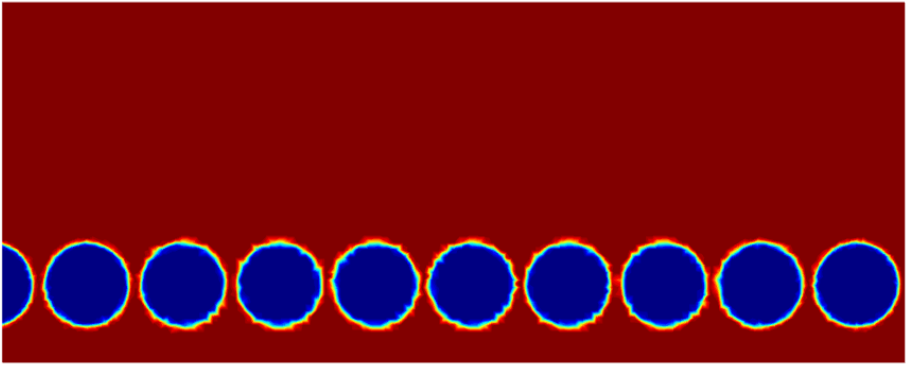}
		\caption{Initial shape}
	\end{subfigure}\hfill
	\begin{subfigure}[b]{0.33\textwidth}
		\centering
		\includegraphics[width=1.8 in]{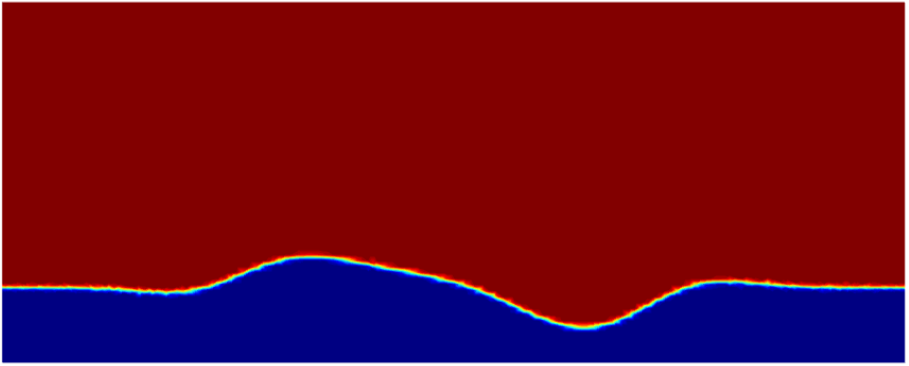}
		\caption{Reference shape}
	\end{subfigure}\hfill
	\begin{subfigure}[b]{0.33\textwidth}
		\centering
		\includegraphics[width=1.8 in]{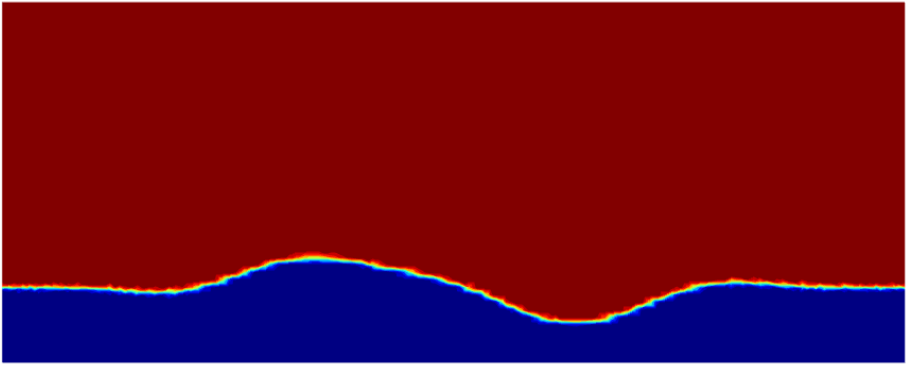}
		\caption{Reconstructed shape}
	\end{subfigure}\par\smallskip
	\begin{subfigure}[h]{0.33\textwidth}
		\centering
		\includegraphics[width=1.8 in]{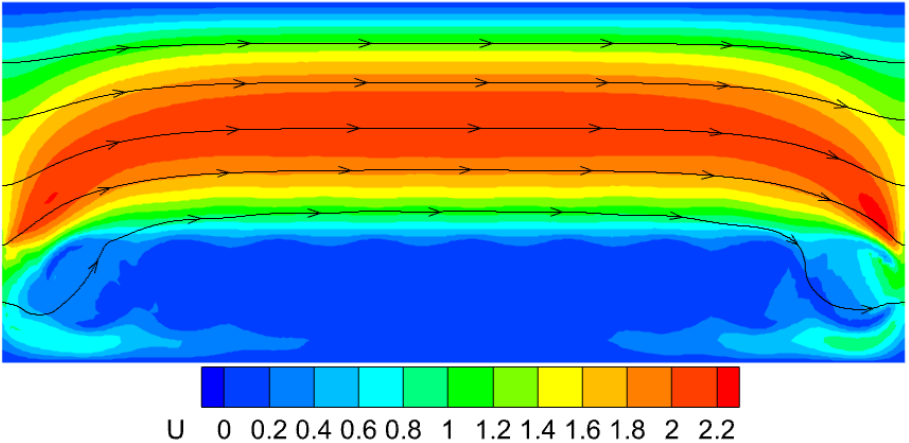}
		\caption{Velocity of initial shape} 	\end{subfigure}
	\begin{subfigure}[h]{0.33\textwidth}
		\centering
		\includegraphics[width=1.8 in]{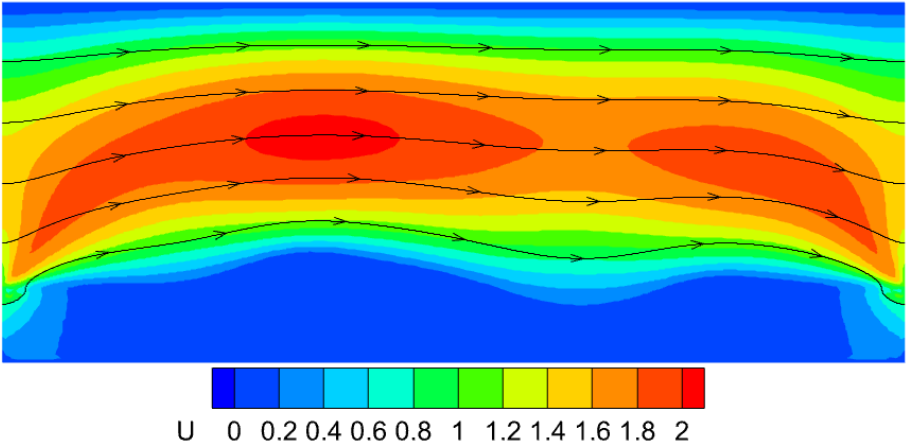}
		\caption{Velocity of reference shape} 	\end{subfigure}
	\begin{subfigure}[h]{0.33\textwidth}
		\centering
		\includegraphics[width=1.8 in]{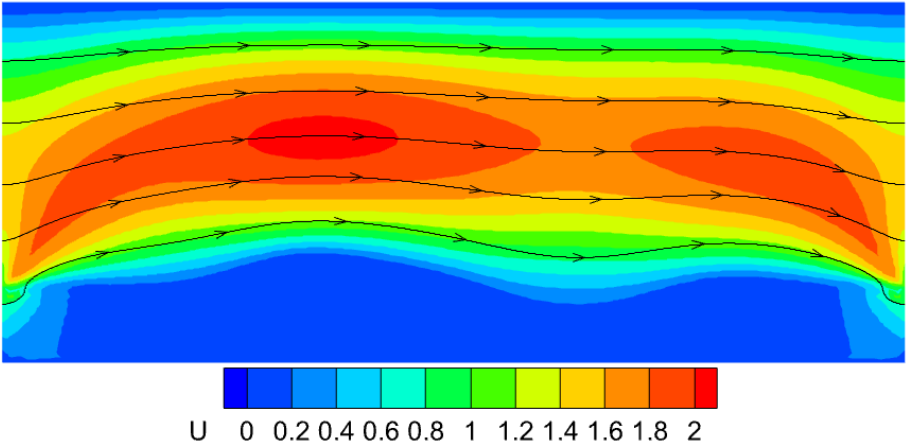} 
		\caption{Velocity of reconstructed shape}
	\end{subfigure}
	\caption{Reconstruction of an idealized riverbed interface from synthetic
		velocity observations.
		(a) {\r Physical motivation: river flow over a permeable sediment bed.}
		(b) Two-dimensional computational domain and boundary conditions.
		(c) History of the regularized objective.
		(d)--(f) initial, reference, and selected reconstructed phase configurations, respectively; and (g)--(i) the corresponding velocity magnitudes.}
	\label{fig:inverse-wave}
\end{figure}

{\r The interaction between river flow and seepage in a permeable bed
(Fig.~\ref{fig:inverse-wave}(a)) motivates this idealized
two-dimensional benchmark for riverbed interface reconstruction.
We use the prescribed interface in Fig.~\ref{fig:inverse-wave}(e)
as the reference bed. The synthetic observation field
$\mathbf{u}_{\mathrm{obs}}$ is obtained by solving the regularized
Stokes--Darcy system with the reference interface represented by a
signed-distance function.}
The computational domain is
$\Omega=[-0.25,1.25]\times[-0.30,0.30]$
(Fig.~\ref{fig:inverse-wave}(b)).
On the left and right boundaries, we prescribe
$\mathbf{u}=(1.5[1-(y/0.30)^2],0)^{\top}$,
while the horizontal walls satisfy $\mathbf{u}=\mathbf{0}$.


No volume constraint is imposed.
To incorporate the riverbed boundary prior, the reaction--diffusion
step imposes $\phi=0.05$ on $y=-0.30$ and $\phi=-0.05$ on $y=0.30$,
with homogeneous Neumann conditions on the vertical boundaries.
We set $\mu=0.1$, $K=10^{-4}$,
$\alpha=5$, $\varepsilon=10^{-2}$,
$\zeta=10^{-8}$, $\ell_s=0.03$,
$\Delta s=5\times10^{-4}$, $\kappa=0.5$, and $\tau=0.1$.
The phase-weight lower bound and perimeter weight are both $10^{-3}$,
and $N_{\rm re}=20$. We evaluate $300$ designs and report the one
with the smallest value of the regularized objective, consisting of the
velocity-tracking term and the perimeter term.
The observation and reconstruction meshes contain 24,158 and 15,806
triangles, respectively, with 12,330 and 8,104 vertices.
The initial porous phase is the intersection of $\Omega$ with ten
disks of radius $0.070$, centered at
$(-0.27+0.16i,-0.170)$, $i=0,\ldots,9$.

Figure~\ref{fig:inverse-wave} presents the history of the regularized
objective together with the initial, reference, and selected reconstructed
phase configurations and their velocity magnitudes.
The regularized objective decreases from $6.26573\times10^{-2}$
for the initial design to its smallest recorded value
$1.56984\times10^{-3}$ at iteration $248$; the latter consists of a
velocity-tracking term of $3.54530\times10^{-5}$ and a perimeter term of
$1.53439\times10^{-3}$. The ten initially separated porous inclusions merge
into a connected sediment layer. The selected interface reproduces the main
undulations of the reference bed, and the corresponding velocity-magnitude
field shows the same principal spatial pattern as the reference field.

\subsubsection{Reconstruction of multiple defects in a porous filtration
membrane from partial noisy observations}
\begin{figure}[!htbp]
			\begin{subfigure}[hb]{0.48\textwidth}
			\centering
			\includegraphics[width=2.8in]{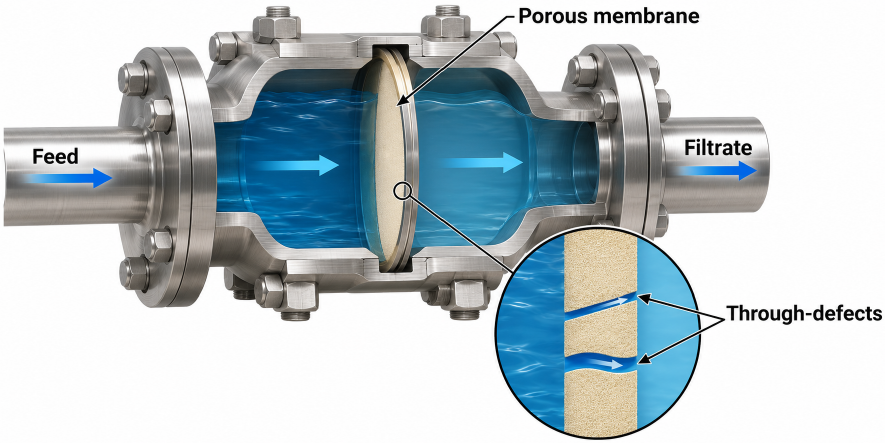}
			\caption{Physical motivation} 	\end{subfigure}
		\begin{subfigure}[hb]{0.48\textwidth}
			\centering
			\includegraphics[width=2.8 in]{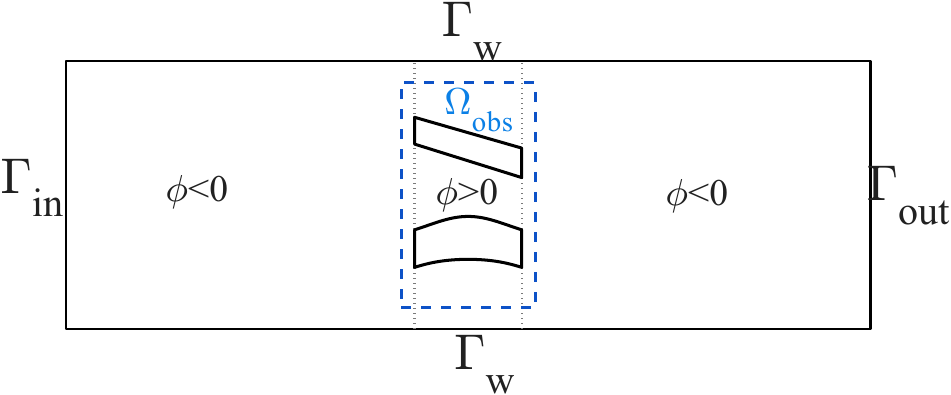}
			\caption{Domain and boundary conditions} 	\end{subfigure}
		\caption{{\r Physical motivation and computational setup for membrane-defect
			reconstruction, showing (a) an in-line membrane filter holder with
			idealized through-defects and (b) the computational domain,
			boundary conditions, and observation window.}}
		\label{fig:membrane-partial 1}
	\end{figure}

{\r The potential for flow to bypass a porous membrane through local
	breaches motivates this inverse problem
	(Fig.~\ref{fig:membrane-partial 1}(a)).
	We consider an idealized two-dimensional benchmark for reconstructing
	multiple through-defects from spatially restricted synthetic velocity
	observations.}
	The computational domain is $\Omega=(0,3)\times(0,1)$, and the nominal
	membrane support is $\Omega_m=[1.3,1.7]\times[0,1]$
	(Fig.~\ref{fig:membrane-partial 1}(b)).
	Using the nominal membrane location as prior information, we restrict
	the level-set update to $D=[1.2,1.8]\times[0.08,0.92]$ and keep the
	level-set degrees of freedom outside $D$ fixed.
	The left and right boundaries are the inlet $\Gamma_{\rm in}$ and
	outlet $\Gamma_{\rm out}$, respectively. No-slip conditions are
	imposed on the horizontal walls $\Gamma_w$.
    

The reference membrane contains a narrow tapered oblique defect in its
upper part and a wider curved defect in its lower part. With
$s=(x-1.3)/0.4$, their centerlines and half-openings are
\[
\begin{aligned}
c_u(s)&=0.74-0.12s,
& a_u(s)&=0.050+0.005s,\\
c_{\ell}(s)&=0.30+0.040\sin(\pi s),
& a_{\ell}(s)&=0.075-0.005\cos(2\pi s).
\end{aligned}
\]
The initial configuration consists of five straight candidate channels
centered at $y=0.14,0.32,0.50,0.68,$ and $0.86$, each with half-opening
$0.035$. Hence, the reconstruction requires both interface deformation
and elimination of three unsupported channels.

Velocity observations are available only in
$\Omega_{\rm obs}=[1.25,1.75]\times[0.08,0.92]$, which occupies
approximately $14\%$ of $\Omega$.
Define
$
g(y):=6y(1-y).
$
Two flux-compatible loading conditions are used to distinguish the
upper and lower defects. For the first loading, we prescribe
$
\mathbf{u}_d^{(1)}=(g(y),0)^{\top} \text{on }
\Gamma_{\rm in}\cup\Gamma_{\rm out}.
$
For the second loading, we prescribe
$\mathbf{u}_d^{(2)}
 =\bigl(g(y)[1+0.35\sin(2\pi y)],0\bigr)^{\top} \text{on } \Gamma_{\rm in}$,
$\mathbf{u}_d^{(2)}=(g(y),0)^{\top} \text{on } \Gamma_{\rm out}$.
Both profiles vanish at $y=0,1$ and are therefore compatible with the
no-slip data on the horizontal walls. Moreover,
$
\int_0^1 g(y)\,\mathrm{d}y
 =
\int_0^1 g(y)[1+0.35\sin(2\pi y)]\,\mathrm{d}y
$
so that each loading satisfies the zero-net-flux condition. The same
boundary profiles are used in the observation-generation and inverse
computations.

Let $\delta\geq0$ denote the relative noise level in the velocity
observations, measured in the $L^2(\Omega_{\rm obs})^d$ norm, with
$\delta=0$ corresponding to noiseless data.
In this example, we consider $\delta\in\{0,0.03,0.10\}$.
For each loading $m=1,2$, the observed velocity field
$\mathbf{u}_{{\rm obs},\delta}^{(m)}\in L^2(\Omega)^d$
is fixed independently of the design variable $\phi_h$.
The localized, normalized counterpart of the objective in
\eqref{eq:discrete-objectives} is
\begin{equation}
\label{eq:partial-two-load-objective}
\mathcal J_{2,h}^{\rm obs}(\phi_h)
=
\frac{1}{4}\sum_{m=1}^{2}
\frac{\displaystyle
\int_\Omega\chi_{\rm obs}
\left|
\mathbf{u}_h^{(m)}
-\mathbf{u}_{{\rm obs},\delta}^{(m)}
\right|^2\,\mathrm{d}x}
{\displaystyle
\int_\Omega\chi_{\rm obs}
\left|
\mathbf{u}_{{\rm obs},\delta}^{(m)}
\right|^2\,\mathrm{d}x}
+\gamma_{\rm per}\int_\Omega
\delta_\varepsilon(\phi_h)g_h(\phi_h)\,\mathrm{d}x .
\end{equation}
Here $\chi_{\rm obs}$ is a fixed smoothed characteristic function of
$\Omega_{\rm obs}$ with smoothing parameter $0.02$, satisfying
$0\leq\chi_{\rm obs}\leq1$.
For each loading $m=1,2$, assume that
\[
D_m:=
\int_\Omega\chi_{\rm obs}
\left|\mathbf{u}_{{\rm obs},\delta}^{(m)}\right|^2
\,\mathrm{d}x>0.
\]
Since $\chi_{\rm obs}$ and the observation fields are independent of
$\phi_h$, the normalization factors $D_m$ are fixed.
Their positivity, the boundedness of $\chi_{\rm obs}$, and the finite
number of loadings ensure that this extension preserves the
continuity and differentiability arguments for the single-load
functional.

For each loading $m=1,2$, the corresponding adjoint problem has the
form \eqref{eq:discrete-adjoint}, with right-hand side
\[
\ell_h^{(m)}(\mathbf{v}_h)
:=
\frac{1}{2D_m}
\int_\Omega\chi_{\rm obs}
\left(
\mathbf{u}_h^{(m)}
-\mathbf{u}_{{\rm obs},\delta}^{(m)}
\right)\cdot\mathbf{v}_h\,\mathrm{d}x .
\]
The design sensitivity is assembled by summing the two adjoint
contributions and adding the perimeter variation once.

Synthetic observations are generated on a fine mesh with $55{,}504$
triangles and projected onto an independent reconstruction mesh with
$23{,}862$ triangles, thereby mitigating the discretization inverse
crime. We consider the cases without noise, with $3\%$ noise, and with
$10\%$ noise. A fixed spatial Gaussian realization is scaled separately
for the two loadings in the $L^2(\Omega_{\rm obs})$ norm and retained
across the noisy tests, so that the comparison isolates the influence
of the noise amplitude. The parameters are $\mu=0.1$, $K=10^{-4}$, $\alpha=0.5$,
$\eta=10^{-3}$, $\varepsilon=0.04$, $\zeta=10^{-8}$, and
$\gamma_{\rm per}=5\times10^{-4}$. For the reaction--diffusion update,
we use $\Delta s=2\times10^{-4}$, $\kappa=0.5$,
$\tau=5\times10^{-4}$, and $\ell_s=2\times10^{-2}$. The reaction
displacement per update is limited to $0.3\varepsilon$. Fast-marching
reinitialization is applied after every update, with $N_{\max}=300$
and $\mathrm{tol}_{\mathcal M}=10^{-7}$.

\begin{figure*}[htbp]
		\centering
		
		\begin{tabular}{c c c c}
			\multirow{2}{*}[0.4cm]{
			\begin{subfigure}[hb]{0.30\textwidth}
				\centering
				\includegraphics[width=2.0 in]{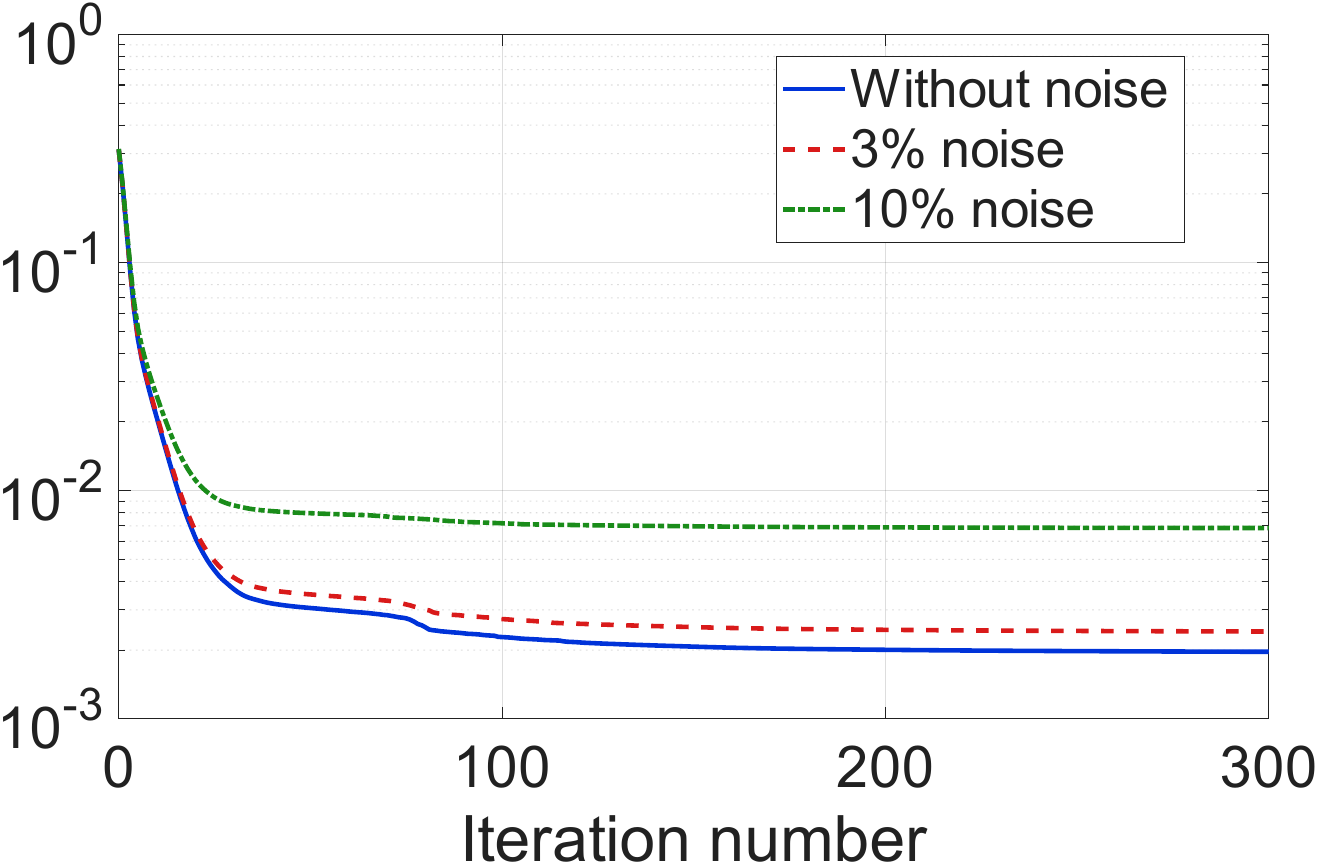}
				\caption{Regularized-objective histories} 	\end{subfigure}
			}
			&
				\begin{subfigure}[hb]{0.30\textwidth}
				\centering
				\includegraphics[width=1.9 in]{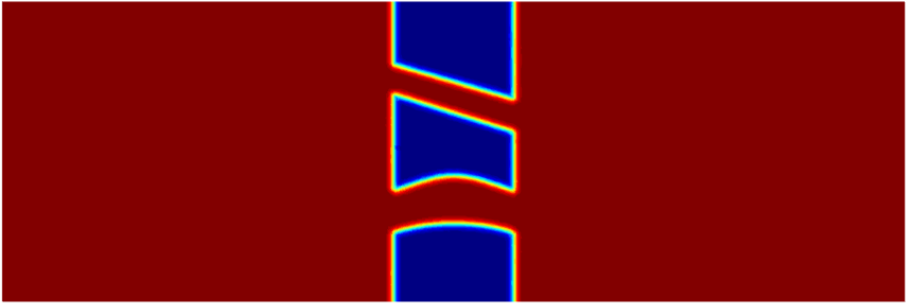}
				\caption{Reference shape}  	\end{subfigure}
			&
			\begin{subfigure}[hb]{0.30\textwidth}
				\centering
				\includegraphics[width=1.9 in]{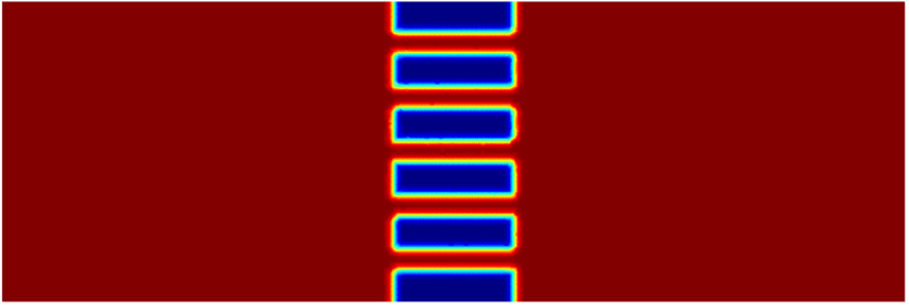}
				\caption{Initial shape}  	\end{subfigure}
				\\[0.1cm]
			
			&
			\begin{subfigure}[hb]{0.30\textwidth}
				\centering
				\includegraphics[width=1.9 in]{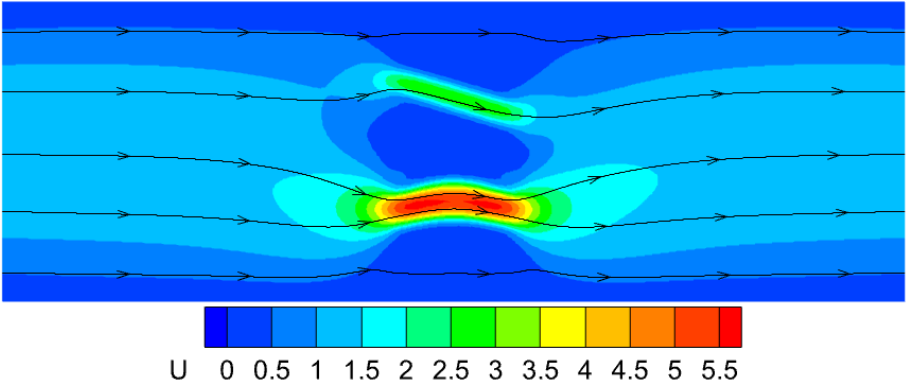}
				\caption{Velocity of reference shape}  	\end{subfigure}
			&
			\begin{subfigure}[hb]{0.30\textwidth}
				\centering
				\includegraphics[width=1.9 in]{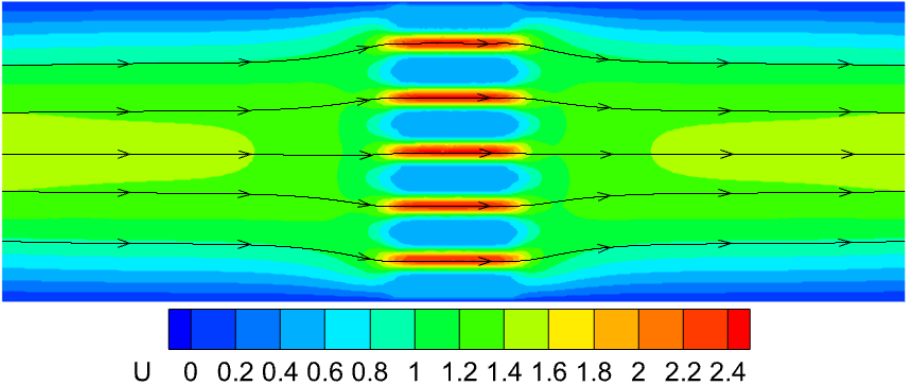}
				\caption{Velocity of initial shape} 	\end{subfigure}
				\\[0.1cm]
					\begin{subfigure}[hb]{0.30\textwidth}
					\centering
					\includegraphics[width=1.9 in]{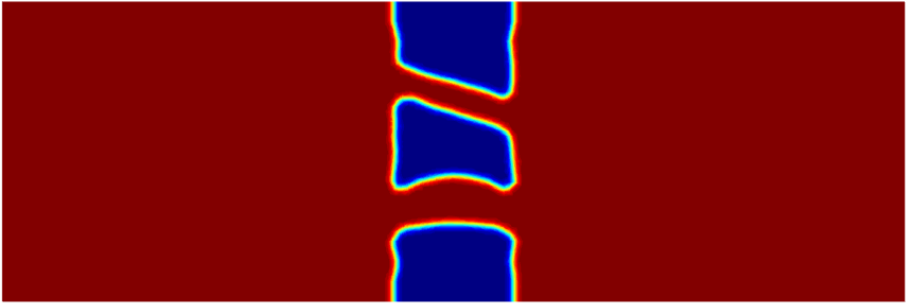}
					\caption{Reconstructed shape without noise} 	\end{subfigure}
			&
					\begin{subfigure}[hb]{0.30\textwidth}
					\centering
					\includegraphics[width=1.9 in]{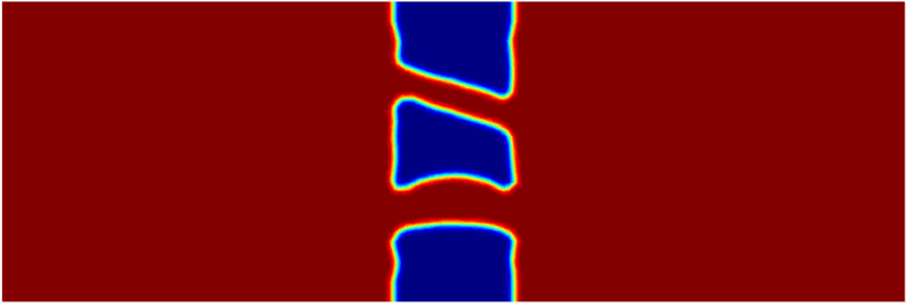}
					\caption{Reconstructed shape with $3\%$ noise}  	\end{subfigure}
			&
				\begin{subfigure}[hb]{0.30\textwidth}
					\centering
					\includegraphics[width=1.9 in]{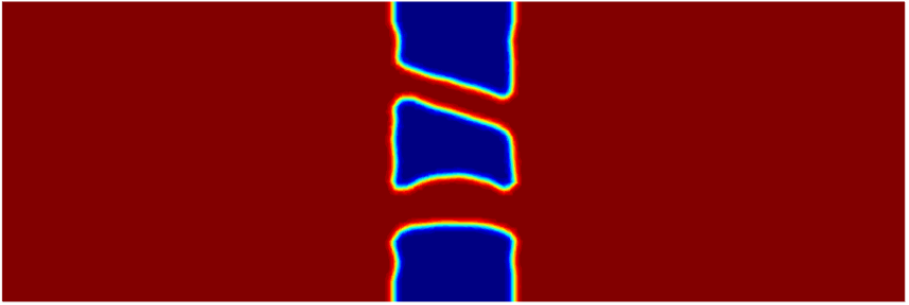}
					\caption{Reconstructed shape with $10\%$ noise} 	\end{subfigure}
			\\[0.1cm]
			
				\begin{subfigure}[h]{0.30\textwidth}
					\centering
					\includegraphics[width=1.9 in]{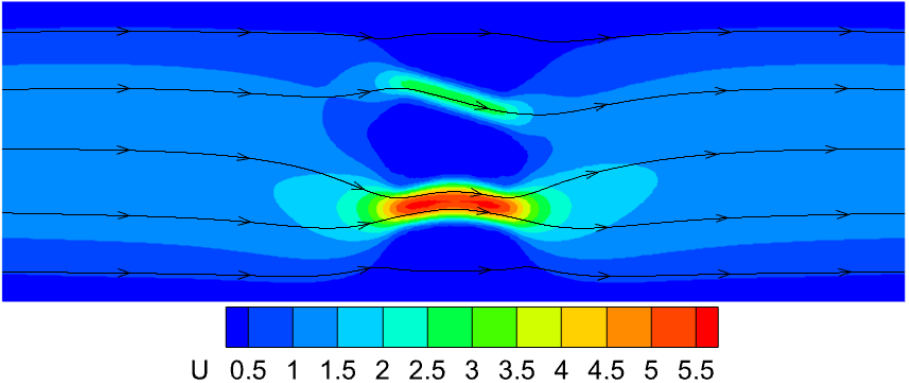}
					\caption{Velocity of reconstructed shape without noise} 	\end{subfigure}
			&
				\begin{subfigure}[h]{0.30\textwidth}
					\centering
					\includegraphics[width=1.9 in]{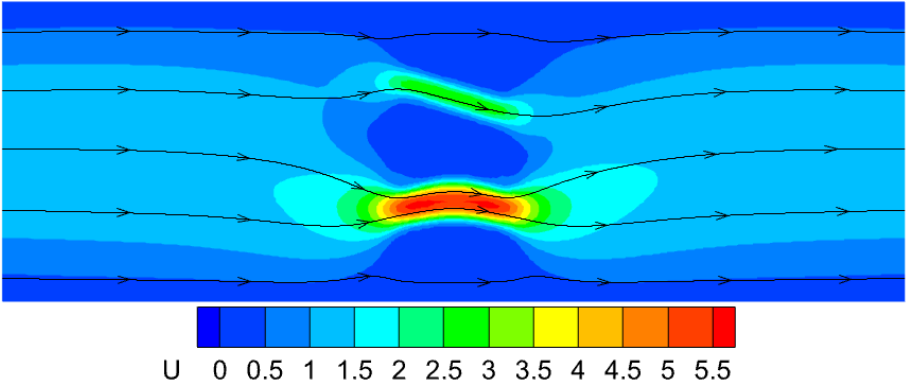}
					\caption{Velocity of reconstructed shape with 3\% noise}	\end{subfigure}
			&
				\begin{subfigure}[h]{0.30\textwidth}
					\centering
					\includegraphics[width=1.9 in]{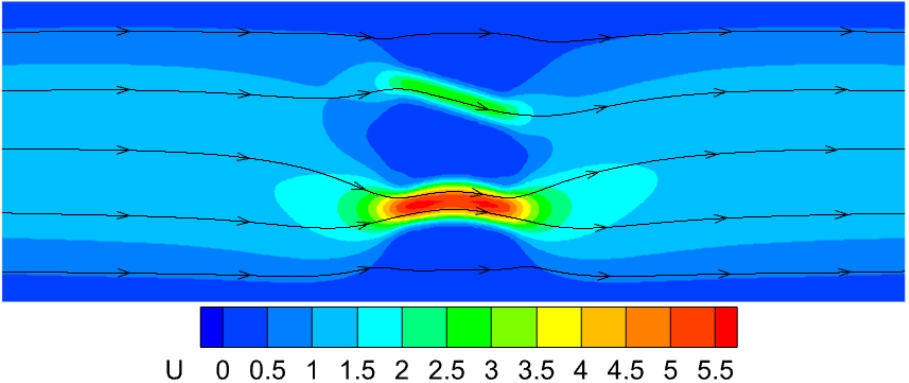}
					\caption{Velocity of reconstructed shape with 10\% noise}	\end{subfigure}
		\end{tabular}
\caption{
	Reconstruction of two through-defects from partial velocity
	observations:
	(a) regularized-objective histories;
	(b) reference configuration;
	(c) initial configuration;
	(d) and (e) velocity magnitudes of the reference and initial
	configurations;
	(f)--(h) reconstructions without noise, with $3\%$ noise, and with
	$10\%$ noise; and
	(i)--(k) corresponding velocity magnitudes.
}
\label{fig:membrane-partial}
	\end{figure*}

All three calculations reach the prescribed limit of $300$ updates;
the merit-change tolerance is not satisfied. The final regularized-objective
values are $1.966\times10^{-3}$, $2.412\times10^{-3}$, and
$6.854\times10^{-3}$ for the $0\%$, $3\%$, and $10\%$ noise cases,
respectively. The corresponding relative residuals in the observation
window are $3.04\%$, $4.26\%$, and $10.34\%$. As shown in
Fig.~\ref{fig:membrane-partial}, the five initial candidate channels are
reduced to two openings in all three cases. Their locations and main
orientations are similar to those of the reference defects, while the
interface details vary with the noise level.

\subsubsection{Three-dimensional reconstruction of a fixed porous particle}

We reconstruct a fixed porous inclusion in the rectangular channel
$\Omega=(0,3)\times(0,1)\times(0,1)$. The level-set update is restricted to
$D=(0.35,2.55)\times(0.08,0.92)\times(0.08,0.92)$, with its degrees of
freedom outside $D$ kept fixed. The initial inclusion is the sphere
\[
 \Omega_p^0=\{\boldsymbol{x}\in\Omega:
 |\boldsymbol{x}-(1.48,0.50,0.50)^{\top}|<0.36\}.
\]
The reference inclusion is a smooth horizontal three-dimensional heart. With
$\theta=10^\circ$ and
\[
\begin{aligned}
 \xi_1&=-\frac{\cos\theta(x-1.58)+\sin\theta(y-0.50)}{0.58},\\
 \xi_2&= \frac{-\sin\theta(x-1.58)+\cos\theta(y-0.50)}{0.42},
 &\xi_3&=\frac{z-0.50}{0.32},
\end{aligned}
\]
it is defined by
\[
 \Omega_p^\dagger=\left\{\boldsymbol{x}\in\Omega:
 (\xi_1^2+\xi_2^2-1)^3-\xi_2^2\xi_1^3+\xi_3^2<0\right\}.
\]

Three flux-compatible boundary loadings are used. On both
$\Gamma_{\rm in}=\{0\}\times(0,1)^2$ and
$\Gamma_{\rm out}=\{3\}\times(0,1)^2$, we prescribe
\[
 \mathbf u_d^{(m)}=
 16y(1-y)z(1-z)\bigl(1+b_y^{(m)}(2y-1)
 +b_z^{(m)}(2z-1),0,0\bigr)^{\top},
\]
where $(b_y^{(m)},b_z^{(m)})=(0,0),(0.6,0),(0,0.6)$ for
$m=1,2,3$, respectively; no-slip conditions are imposed on the remaining
walls. Synthetic full-domain velocity observations are generated on an
independent mesh. The functional used in this example is
\[
 \mathcal J_{2,h}^{(3)}(\phi_h)
 =\frac{1}{6}\sum_{m=1}^3\int_\Omega
 |\mathbf u_h^{(m)}-\mathbf u_{\rm obs}^{(m)}|^2\,\mathrm dx
 +\gamma_{\rm per}\int_\Omega
 \delta_\varepsilon(\phi_h)g_h(\phi_h)\,\mathrm dx.
\]
Accordingly, the design sensitivity is obtained by summing the three weighted
adjoint contributions and adding the perimeter variation once.

The observation and reconstruction meshes contain $26{,}808$ and $10{,}576$
tetrahedra, respectively. We set $\mu=1$, $K=10^{-3}$, $\alpha=1$,
$\eta=10^{-3}$, $\varepsilon=0.10$, $\zeta=10^{-8}$, and
$\gamma_{\rm per}=5\times10^{-6}$. For the reaction--diffusion update, we use
$\ell_s=0.06$, $\Delta s=0.02$, $\kappa=1$, and $\tau=10^{-5}$, and limit the
reaction displacement to $0.1\varepsilon$ per update. No volume constraint is
imposed. Fast-marching reinitialization is applied once, after update 40, and
the calculation is continued to the prescribed limit $N_{\max}=120$.

\begin{figure}[!htbp]
\centering
\begin{subfigure}[t]{0.50\textwidth}
  \centering\includegraphics[width=2.40in]{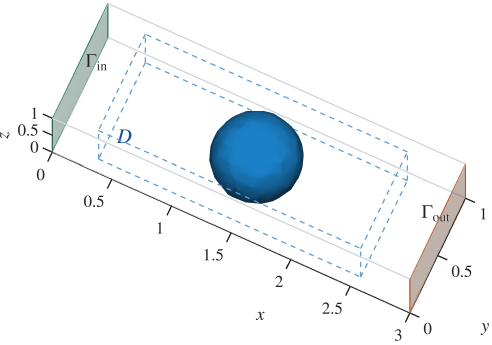}
  \caption{Domain, boundary conditions, and design region}
\end{subfigure}
\begin{subfigure}[t]{0.48\textwidth}
  \centering\includegraphics[width=2.0in]{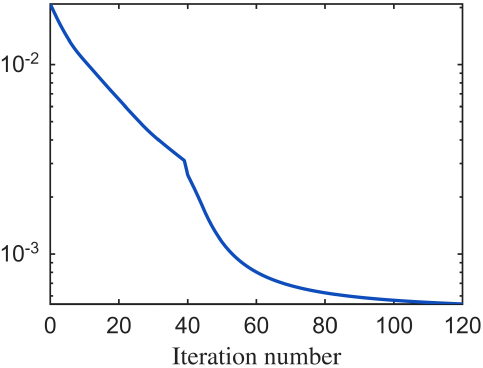}
  \caption{Regularized-objective history}
\end{subfigure}

\begin{subfigure}[t]{0.32\textwidth}
  {\b \centering\includegraphics[width=\linewidth]{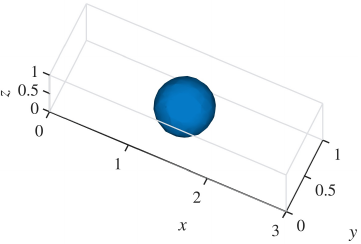}}
  \caption{Initial particle}
\end{subfigure}
\begin{subfigure}[t]{0.32\textwidth}
  {\b \centering\includegraphics[width=\linewidth]{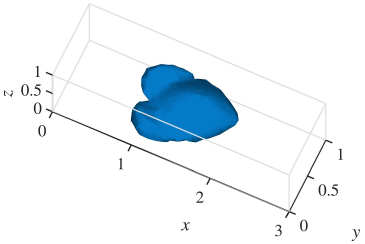}}
  \caption{Reference particle}
\end{subfigure}
\begin{subfigure}[t]{0.32\textwidth}
  {\b \centering\includegraphics[width=\linewidth]{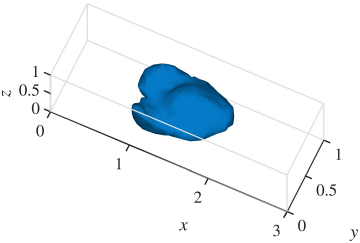}}
  \caption{Reconstructed particle}
\end{subfigure}

\begin{subfigure}[t]{0.32\textwidth}
  {\b \centering\includegraphics[width=\linewidth]{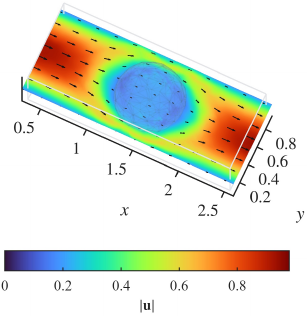}}
  \caption{Initial velocity magnitude}
\end{subfigure}
\begin{subfigure}[t]{0.32\textwidth}
  {\b \centering\includegraphics[width=\linewidth]{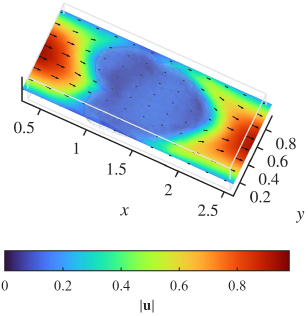}}
  \caption{Reference velocity magnitude}
\end{subfigure}
\begin{subfigure}[t]{0.32\textwidth}
  {\b \centering\includegraphics[width=\linewidth]{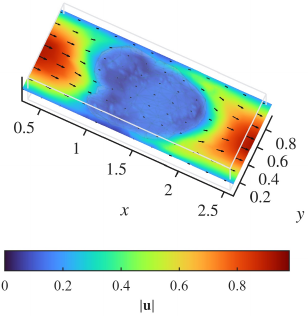}}
  \caption{Reconstructed velocity magnitude}
\end{subfigure}
\caption{Three-dimensional reconstruction of a fixed porous particle from
three full-domain noiseless velocity observations: (a) computational domain,
boundary conditions, and design region; (b) recorded regularized-objective
history; (c)--(e) initial, reference, and reconstructed interfaces in the
complete computational domain; and
(f)--(h) the corresponding velocity magnitudes on the oblique section
$z=0.2+0.6y$. The interface panels use a common view, and the velocity panels
use a common color scale.}
\label{fig:3d_inverse_heart}
\end{figure}

Figure~\ref{fig:3d_inverse_heart}(b) shows the recorded regularized objective,
which decreases from $2.0881\times10^{-2}$ to $5.4375\times10^{-4}$. The
aggregated relative velocity error decreases from $2.3288\times10^{-1}$ to
$3.7175\times10^{-2}$. Figures~\ref{fig:3d_inverse_heart}(c)--(e) show that the
reconstructed particle captures the displacement, horizontal orientation,
two lobes, central indentation, and tapered end of the reference particle,
although the final indentation remains shallower. The corresponding velocity
magnitudes on the common oblique section $z=0.2+0.6y$ are shown in
Figs.~\ref{fig:3d_inverse_heart}(f)--(h). The calculation terminates at the
prescribed update limit; the history is reported as a fixed-budget numerical
reconstruction.

\section{Conclusions}

The numerical results illustrate how changes in channel connectivity
and interface geometry can improve flow designs and support the
reconstruction of unknown fluid--porous configurations. They also
highlight the importance of distinguishing geometric effects from
those introduced by permeability and diffuse-interface regularization.
The present existence analysis is restricted to the specified
admissible classes and fixed regularization parameters.

Future work
will therefore address broader design classes and the consistency
of the diffuse formulation with the sharp-interface model.
Adaptive discretization and scalable solvers will be investigated
for larger three-dimensional problems. For inverse reconstruction,
further priorities include uncertainty quantification and validation
against experimental velocity measurements.

\appendix

\section{Compactness and existence of optimal designs}
\label{app:existence}

We collect the compactness and state-closedness results used in the
existence analysis, followed by the proof of
Theorem~\ref{thm:existence}.

\begin{lemma}[Compactness of admissible level-set functions]
\label{lem:admissible-compactness}
Let $j\in\{1,2\}$ and let $\{\phi_n\}\subset\mathcal A_j$.  Then there
are a subsequence, not relabeled, and $\phi^*\in\mathcal A_j$ such that
\[
 \phi_n\stackrel{*}{\rightharpoonup}\phi^*
   \text{ in }W^{1,\infty}(\Omega),
 \qquad
 \phi_n\to\phi^*\text{ uniformly on }\overline\Omega,
 \qquad
 \phi_n\to\phi^*\text{ strongly in }H^1(\Omega).
\]
Moreover,
\[
 \mathbf P_\tau(\phi_n)\to\mathbf P_\tau(\phi^*)
 \quad\text{strongly in }L^2(\Omega)^{d\times d}.
\]
\end{lemma}

\begin{proof}
The closed uniform $W^{1,\infty}$ bound gives, after extraction,
weak-$*$ convergence in $W^{1,\infty}$ and uniform convergence on
$\overline\Omega$. By the stability of viscosity solutions under locally uniform
convergence, the limit $\phi^*$ is again a viscosity solution
of $|\nabla\phi^*|=1$ in $\Omega$
\cite{CrandallLions1983,CrandallIshiiLions1992}.
Since Lipschitz viscosity solutions satisfy the
Eikonal equation almost everywhere, $|\nabla\phi_n|=|\nabla\phi^*|=1$ a.e. in $\Omega$.
The gradients converge weakly in $L^2(\Omega)^d$ and have the same norm
$|\Omega|^{1/2}$; uniform convexity of $L^2$ therefore yields strong
convergence of the gradients.  Together with uniform convergence of the
functions, this proves strong $H^1$ convergence.

If $j=1$, continuity of $\mathcal H_f$ and uniform convergence give
\[
 \int_\Omega\mathcal H_f(\phi^*)\dx
 =\lim_{n\to\infty}\int_\Omega\mathcal H_f(\phi_n)\dx=V_0.
\]
For $j=2$, no volume constraint is imposed.  The norm constraint is closed
under weak-$*$ convergence, hence $\phi^*\in\mathcal A_j$ in either case.
Finally, because the gradients have unit length,
$\mathbf P_\tau(\phi_n)
 =\mathbf I-\nabla\phi_n\otimes\nabla\phi_n$,
and $\|\mathbf P_\tau(\phi_n)-\mathbf P_\tau(\phi^*)\|_{L^2}
 \le2\|\nabla\phi_n-\nabla\phi^*\|_{L^2}\to0$.
\end{proof}

\begin{lemma}[Sequential closedness of the state system]
\label{lem:state-closedness}
Let $j\in\{1,2\}$ and let $\phi_n\in\mathcal A_j$ converge as in
Lemma~\ref{lem:admissible-compactness}, and let
$(\mathbf u_n,p_n)=S(\phi_n)$.  Then, after extraction,
\[
 \mathbf u_n\rightharpoonup\mathbf u^*
   \text{ in }H^1(\Omega)^d,
 \qquad
 \mathbf u_n\to\mathbf u^*
   \text{ in }L^2(\Omega)^d,
 \qquad
 p_n\rightharpoonup p^*
   \text{ in }L^2_0(\Omega),
\]
where $(\mathbf u^*,p^*)=S(\phi^*)$.
\end{lemma}
\begin{proof}
Lemma~\ref{lem:pde-wellposed} and the uniform $W^{1,\infty}$ bounds in
\eqref{eq:admissible-set} and \eqref{eq:inverse-admissible-set} imply that
$\{(\mathbf u_n,p_n)\}$ is uniformly
bounded in $H^1(\Omega)^d\times L^2_0(\Omega)$. Weak compactness and the compact embedding
$H^1(\Omega)^d\hookrightarrow L^2(\Omega)^d$ yield, after passing to a
further subsequence, the stated weak and strong convergences. The weak
closedness of $L^2_0(\Omega)$ yields $p^*\in L^2_0(\Omega)$. Moreover, all
states have trace $\mathbf u_d$; hence weak continuity of the trace map gives
$\mathbf u^*\in\mathcal U_{\mathbf u_d}$.

Uniform convergence of $\phi_n$ implies uniform convergence of
$\mathcal H_f(\phi_n)$, $\mathcal H_p(\phi_n)$, and
$\delta_\varepsilon(\phi_n)$.  For every fixed $\mathbf v\in\mathcal U$,
$\mathcal H_f(\phi_n)\mathcal E(\mathbf v)$ converges strongly in $L^2$,
whereas $\mathcal E(\mathbf u_n)$ converges weakly in $L^2$.  Hence
\[
 \int_\Omega2\mu\,\mathcal E(\mathbf u_n):\mathcal E(\mathbf v)
 \mathcal H_f(\phi_n)\dx
 \to
 \int_\Omega2\mu\,\mathcal E(\mathbf u^*):\mathcal E(\mathbf v)
 \mathcal H_f(\phi^*)\dx.
\]
Similarly, the strong $L^2$ convergence of $\mathbf u_n$ and the uniform
convergence of $\mathcal H_p(\phi_n)$ give
\[
 \int_\Omega\frac{\mu}{K}\mathbf u_n\cdot\mathbf v\,
 \mathcal H_p(\phi_n)\dx
 \to
 \int_\Omega\frac{\mu}{K}\mathbf u^*\cdot\mathbf v\,
 \mathcal H_p(\phi^*)\dx.
\]

For the Beavers--Joseph--Saffman (BJS) interface term, set
$\mathbf P_n:=\mathbf P_\tau(\phi_n)$,
$\mathbf P_*:=\mathbf P_\tau(\phi^*)$, and
$\delta_n:=\delta_\varepsilon(\phi_n)$.  By
Lemma~\ref{lem:admissible-compactness}, the uniform $L^\infty$ boundedness of
the projections, and interpolation,
\[
 \|\mathbf P_n-\mathbf P_*\|_{L^3}
 \le \|\mathbf P_n-\mathbf P_*\|_{L^2}^{2/3}
     \|\mathbf P_n-\mathbf P_*\|_{L^\infty}^{1/3}\to0.
\]
Since
$H^1(\Omega)^d\hookrightarrow L^6(\Omega)^d$ for $d\le3$ and
$\delta_n\to\delta_\varepsilon(\phi^*)$ uniformly, for every fixed
$\mathbf v\in\mathcal U$,
\[
 \delta_n\mathbf P_n\mathbf v
 \to\delta_\varepsilon(\phi^*)\mathbf P_*\mathbf v
 \quad\text{strongly in }L^2(\Omega)^d.
\]
Because $|\nabla\phi_n|=|\nabla\phi^*|=1$ almost everywhere and
$\mathbf P_n^\top\mathbf P_n=\mathbf P_n$, the strong $L^2$ convergence of
$\mathbf u_n$ and the preceding convergence imply
\[
 \begin{aligned}
 \int_\Omega\frac{\alpha\mu}{\sqrt K}
 (\mathbf P_n\mathbf u_n)\cdot(\mathbf P_n\mathbf v)
 \delta_n|\nabla\phi_n|\dx =\int_\Omega\frac{\alpha\mu}{\sqrt K}
 \mathbf u_n\cdot(\delta_n\mathbf P_n\mathbf v)\dx
 \to\int_\Omega\frac{\alpha\mu}{\sqrt K}
 (\mathbf P_*\mathbf u^*)\cdot(\mathbf P_*\mathbf v)
 \delta_\varepsilon(\phi^*)|\nabla\phi^*|\dx.
 \end{aligned}
\]

Finally, weak convergence gives
$b(p_n,\mathbf v)\to b(p^*,\mathbf v)$ and
$b(q,\mathbf u_n)\to b(q,\mathbf u^*)$.  Passing to the limit in
\eqref{eq:global-ls} and using uniqueness from
Lemma~\ref{lem:pde-wellposed} yield
$(\mathbf u^*,p^*)=S(\phi^*)$.
\end{proof}

\begin{proof}[Proof of Theorem~\ref{thm:existence}]
Fix $j\in\{1,2\}$ and let $\{\phi_n\}\subset\mathcal A_j$ be a minimizing
sequence.  By
Lemmas~\ref{lem:admissible-compactness} and \ref{lem:state-closedness}, after
passing to a subsequence,
\[
 \phi_n\to\phi^*\text{ uniformly},\qquad
 \mathbf u_n\rightharpoonup\mathbf u^*\text{ in }H^1(\Omega)^d,
 \qquad
 \mathbf u_n\to\mathbf u^*\text{ in }L^2(\Omega)^d,
\]
where $\phi^*\in\mathcal A_j$ and
$(\mathbf u^*,p^*)=S(\phi^*)$.

For $\mathcal J_1$, uniform convergence of
$\mathcal H_f(\phi_n)$ implies
$\sqrt{\mathcal H_f(\phi_n)}\to
\sqrt{\mathcal H_f(\phi^*)}$ uniformly.  Therefore
\[
 \sqrt{\mathcal H_f(\phi_n)}\mathcal E(\mathbf u_n)
 \rightharpoonup
 \sqrt{\mathcal H_f(\phi^*)}\mathcal E(\mathbf u^*)
 \quad\text{in }L^2,
\]
and weak lower semicontinuity of the squared $L^2$ norm gives
\[
 \mathcal J_1(\mathbf u^*,\phi^*)
 \le\liminf_{n\to\infty}\mathcal J_1(\mathbf u_n,\phi_n).
\]
For $\mathcal J_2$, the strong $L^2$ convergence of $\mathbf u_n$ gives
\[
 \frac12\int_\Omega|\mathbf u_n-\mathbf u_{\rm obs}|^2\dx
 \to
 \frac12\int_\Omega|\mathbf u^*-\mathbf u_{\rm obs}|^2\dx.
\]
Moreover, since $|\nabla\phi_n|=|\nabla\phi^*|=1$ almost everywhere and
$\delta_\varepsilon(\phi_n)\to\delta_\varepsilon(\phi^*)$ uniformly,
\[
 \mathcal P_\varepsilon(\phi_n)
 =\int_\Omega\delta_\varepsilon(\phi_n)\dx
 \to\int_\Omega\delta_\varepsilon(\phi^*)\dx
 =\mathcal P_\varepsilon(\phi^*).
\]
Consequently,
\[
 \mathcal J_2(\mathbf u_n,\phi_n)
 \to\mathcal J_2(\mathbf u^*,\phi^*).
\]
Thus the reduced objective is sequentially lower semicontinuous in either
case.  Applying the direct method of the calculus of variations proves
that $\phi^*$ is a minimizer.
\end{proof}

\bibliographystyle{unsrtnat}
\bibliography{ref}

\end{document}